%% file: American_basket.tex
\documentclass[10pt]{article}
\usepackage{graphicx,float,latexsym,amsmath,amsfonts,amssymb,subeqnarray,colordvi}
\usepackage[latin1]{inputenc}
\usepackage{a4wide}
\usepackage{color}
\usepackage{mathabx}
\usepackage{pgfplots}
\usepackage{rotating}

\definecolor{color1}{RGB}{0, 96, 0}
\definecolor{color2}{RGB}{255, 140, 0}
\definecolor{color3}{RGB}{0, 0, 192}
\definecolor{color4}{RGB}{192, 0, 0}

\def\R{\mathbb{R}}
\newcommand{\rT}{{\rm{T}}}

\providecommand*{\diff}{\@ifnextchar^{\DIfF}{\DIfF^{}}}
\def\DIfF^#1{\operatorname{d}\nolimits^{#1}\gobblespace}
\def\gobblespace{\futurelet\diffarg\opspace}
\def\opspace{\let\DiffSpace\!
	\ifx\diffarg(\let\DiffSpace\relax
	\else
	\ifx\diffarg[\let\DiffSpace\relax
	\else
	\ifx\diffarg\{\let\DiffSpace\relax
	\fi\fi\fi\DiffSpace}

\providecommand*{\pderiv}[3][]{\frac{\partial^{#1}#2}{{\partial #3}^{#1}}}

\title{Numerical valuation of American basket options\\ 
via partial differential complementarity problems}
\author{Karel~J.~in 't Hout
\footnote{Department of Mathematics,
University of Antwerp, Middelheimlaan 1, B-2020 Antwerp, Belgium.
\mbox{Email}: \texttt{\{karel.inthout,jacob.snoeijer\}@uantwerpen.be}.}
~and Jacob Snoeijer
\footnotemark[\value{footnote}]
}
\date{\today}

\begin{document}

\maketitle

\begin{abstract}
\noindent
We study the principal component analysis based approach introduced by Reisinger \& Wittum \cite{RW07}
and the comonotonic approach considered by Hanbali \& Linders \cite{HL19} for the approximation of American 
basket option values via multidimensional partial differential complementarity problems (PDCPs).
Both approximation approaches require the solution of just a limited number of low-dimensional PDCPs.
It is demonstrated by ample numerical experiments that they define approximations that lie close to each 
other. Next, an efficient discretisation of the pertinent PDCPs is presented that leads to a favourable 
convergence behaviour.
\end{abstract}

\vspace{0.2cm}\noindent

\normalsize

\setcounter{equation}{0}
\section{Introduction}\label{sec:introduction}

This paper deals with the valuation of American-style basket options.
Basket options constitute a popular type of financial derivatives and possess a payoff depending on a weighted average of different 
assets.
In general, exact valuation formulas for such options are not available in the literature in semi-closed analytic form.
Therefore, the development and analysis of efficient approximation methods for their fair values is of much importance.

In this paper we consider the valuation of American basket options through partial differential complementarity problems (PDCPs).
If $d$ denotes the number of different assets in the basket, then the pertinent PDCP is $d$-dimensional.
In this paper, we are interested in the situation where $d$ is medium or large, say $d\ge 5$.
It is well-known that this renders the application of standard discretisation methods for PDCPs impractical, due to the curse of 
dimensionality.
For European- and Bermudan-style basket options, leading to high-dimensional partial differential equations (PDEs), an effective 
approach has been introduced by Reisinger \& Wittum \cite{RW07} and next studied in e.g.~Reisinger \& Wissmann~\cite{RW15,RW17,RW18} 
and In~'t~Hout \& Snoeijer~\cite{HS21}. 
This approach is based on a principal component analysis (PCA) and yields an approximation formula for the value of the basket option 
that requires the solution of a limited number of only low-dimensional PDEs.
In the literature, an alternative useful approach has been investigated that employs the idea of comonotonicity.
For European basket options, this comonotonic approach has been developed notably by
Kaas et al.~\cite{KDG00}, 
Dhaene et al.~\cite{DDGKV02a,DDGKV02b}, 
Deelstra et al.~\cite{DLV04,DDV08} and 
Chen et al.~\cite{CDDV08,CDDLV15}.
Recently, an extension to American basket options has been presented by Hanbali \& Linders \cite{HL19}, who consider a comonotonic 
approximation formula that requires the solution of just two one-dimensional PDCPs.
In the present paper we shall study and compare the PCA-based and comonotonic approaches for the effective valuation of American 
basket options.
To our knowledge, this is the first paper where these two, different but related, approaches are jointly investigated.
In our subsequent analysis, we shall include also the (simpler) case of European basket options.

A \emph{European-style basket option} is a financial contract that gives the holder the right to buy or sell a prescribed weighted 
average of $d$ assets at a prescribed maturity date $T$ for a prescribed strike price $K$.
We assume in this paper the well-known Black--Scholes model.
Thus the asset prices $S_\tau^i$ ($i=1,2,\ldots,d$) evolve according to a multidimensional geometric Brownian motion, 
which is given (under the risk-neutral measure) by the system of stochastic differential equations (SDEs)
\begin{equation}\label{eq:SDE}
dS_\tau^i = r S_\tau^i d\tau + \sigma_i S_\tau^i dW_\tau^i \quad (0< \tau \le T,~1\le i\le d).
\end{equation}
Here $\tau$ is time, with $\tau=0$ representing the time of inception of the option, $r\ge 0$ is the given risk-free interest rate, $\sigma_i>0$ 
($i=1,2,\ldots,d$) are the given volatilities and $W_\tau^i$ ($i=1,2,\ldots,d$) is a multidimensional standard Brownian motion with given 
correlation matrix $(\rho_{ij})_{i,j=1}^d$.
Further, the initial asset prices $S_0^i>0$ ($i=1,2,\ldots,d$) are given.
In essentially all financial applications, the correlation matrix is full.

Let $u(s,t) = u(s_1,s_2,\ldots,s_d,t)$ be the fair value of a European basket option if at time  till maturity $t = T-\tau$ the $i$-th 
asset price equals $s_i$ ($i=1,2,\ldots,d$).
Financial mathematics theory yields that $u$ satisfies the $d$-dimensional time-dependent PDE
\begin{equation}\label{eq:PDE_s}
\pderiv{u}{t}(s,t) =
{\cal A}u(s,t) := 
\frac{1}{2} \sum_{i=1}^d \sum_{j=1}^d \sigma_i \sigma_j \rho_{ij} s_i s_j \frac{\partial^2 u}{\partial s_i \partial s_j}(s,t)
+ \sum_{i=1}^d r s_i \pderiv{u}{s_i}(s,t) - r u(s,t)
\end{equation}
whenever $(s, t) \in (0,\infty)^d \times (0,T]$.
The PDE \eqref{eq:PDE_s} is also satisfied if $s_i=0$ for any given $i$, thus at the boundary of the spatial domain.
At maturity time of the option its fair value is known and specified by the particular option contract.
If $\phi$ is the given payoff function of the option, then one has the initial condition
\begin{equation}\label{eq:IC_s}
u(s,0) = \phi(s)
\end{equation}
whenever $s \in (0,\infty)^d$.

An \emph{American-style basket option} is a financial contract that gives the holder the right to buy or sell a prescribed weighted 
average of $d$ assets for a prescribed strike price $K$ at any given single time up to and including a prescribed maturity time $T$.
The fair value function $u$ of an American basket option satisfies the (nonlinear) $d$-dimensional time-dependent PDCP
\begin{subeqnarray}\label{eq:PDCP_u}
&&u(s,t) \ge \phi (s),\\
\nonumber\\
&&\frac{\partial u} {\partial t} (s,t)  \ge {\cal A } u (s,t),\\
\nonumber\\
&&\left( u(s,t)-\phi (s) \right) \left( \frac{\partial u} {\partial t} (s,t) - {\cal A } u (s,t) \right) = 0
\end{subeqnarray}
\noindent
whenever $(s, t) \in (0,\infty)^d \times (0,T]$.
The PDCP \eqref{eq:PDCP_u} is provided with the same initial condition \eqref{eq:IC_s}.
Further, \eqref{eq:PDCP_u} also holds if $s_i=0$ for any given $i$.

In this paper we shall consider the class of basket put options.
These have a payoff function given by
\begin{equation}\label{eq:payoff}
\phi(s) = \max \left(K - \sum_{i=1}^d \omega_i s_i\,,\, 0\right)
\end{equation}
with prescribed weights $\omega_i>0$ ($i=1,2,\ldots,d$) such that $\sum_{i=1}^d \omega_i = 1$.

An outline of our paper is as follows.

Following Reisinger \& Wittum \cite{RW07}, in Section~\ref{sec:transformation} a convenient coordinate transformation is applied 
to the PDE \eqref{eq:PDE_s} for European basket options by means of a spectral decomposition of the covariance matrix.
This way, a $d$-dimensional time-dependent PDE for a transformed option value function is obtained in which each 
coefficient is directly proportional to one of the eigenvalues.
In Section~\ref{sec:approximation} this feature is exploited to derive a principal component analysis (PCA) based approximation. 
The key property of this approximation is that it is determined by just a limited number of one- and two-dimensional 
PDEs.
The presentation in Sections~\ref{sec:transformation} and~\ref{sec:approximation} follows largely that in \cite{HS21}.
In Section~\ref{sec:approximation_Am}, the PCA-based approximation approach is extended to American basket options.
This gives rise to an approximation that is defined by a limited number of one- and two-dimensional PDCPs.
In Section~\ref{sec:EU_discretisation}, an efficient discretisation of the one- and two-dimensional PDEs for European basket options 
is described, which employs finite differences on a nonuniform spatial grid followed by the Brian and Douglas alternating direction
implicit (ADI) scheme on a uniform temporal grid.
This discretisation is adapted in Section~\ref{sec:AM_discretisation} to the pertinent PDCPs for American basket options, where the
basic explicit payoff (EP) approach as well as the more advanced Ikonen--Toivanen (IT) splitting technique are considered.
Section~\ref{sec:comon} collects results from the literature on the comonotonic approach for valuing European and American basket 
options.
We consider the same comonotonic approximation as Hanbali \& Linders \cite{HL19}, which is determined 
by just two one-dimensional PDEs (for the European basket) or PDCPs (for the American basket). 
Section~\ref{sec:numericalexperiments} contains the main contribution of our paper. 
In this section we perform ample numerical experiments and obtain the positive result that the PCA-based and comonotonic approaches 
yield approximations to the option value that always lie close to each other for both European and American basket put options.
We next study in detail the error in the discretisation described in Section~\ref{sec:discretisation} for the PCA-based 
and comonotonic approximations and observe a favourable, near second-order convergence behaviour. 
The final Section~\ref{sec:conclusion} presents our conclusions and outlook.

\setcounter{equation}{0}
\section{PCA approximation approach}\label{sec:pcs}

\subsection{Coordinate transformation}\label{sec:transformation}
In this preliminary section we apply two subsequent coordinate transformations to the PDE \eqref{eq:PDE_s} for a European basket 
option.
We assume here that the elementary functions $\ln$, $\exp$, $\tan$, $\arctan$ are taken componentwise whenever their argument 
is a vector.

The covariance matrix $\Sigma = \left(\Sigma_{ij} \right) \in \R^{d \times d}$ is given by $\Sigma_{ij} = \sigma_i \rho_{ij} \sigma_j$ 
for $i,j = 1,2,\ldots,d$.
Let $\Lambda = {\rm diag} (\lambda_1,\lambda_2,\ldots,\lambda_d)$ denote a real diagonal matrix of eigenvalues of $\Sigma$ and 
$Q$ a real orthogonal matrix of eigenvectors of $\Sigma$ such that $\Sigma = Q \Lambda Q^\rT$.
Then, following \cite{RW07}, we apply the coordinate transformation
\begin{equation}\label{eq:transformationS2X}
x(s, t) = Q^\rT \left( \ln (s/K) - b(t) \right),
\end{equation}
where $b(t) = (b_1(t),b_2(t),\ldots,b_d(t))^\rT$ with $b_i(t) = (\tfrac{1}{2}\sigma_i^2-r)t$ for $1\le i\le d$.
Let the function $v$ be defined by
\begin{equation*}
u(s,t) = v(x(s,t),t).
\end{equation*}
An easy calculation yields that $v$ satisfies
\begin{equation}\label{eq:PDE_x}
\pderiv{v}{t}(x,t) = \frac{1}{2} \sum_{k=1}^d \lambda_k  \frac{\partial^2 v}{\partial x_k^2}(x,t) - r v(x,t)
\end{equation}
whenever $x \in \R^d$, $t\in (0,T]$.
Clearly, \eqref{eq:PDE_x} is a pure diffusion equation, without mixed derivative terms, and with a simple 
reaction term.
Following \cite{RW07}, we apply a second coordinate transformation, which maps the spatial domain $\R^d$ 
onto the $d$-dimensional open unit cube,
\begin{equation}\label{eq:transformationX2Y}
y(x) = \frac{1}{\pi} \arctan (x) + \frac{1}{2}.
\end{equation}
Let the function $w$ be defined by
\begin{equation*}
v(x,t) = w(y(x),t).
\end{equation*}
Then it is readily seen that 
\begin{equation}\label{eq:PDE_y}
\pderiv{w}{t}(y,t) =
{\cal B}w(y,t) :=
\sum_{k=1}^d \lambda_k  \left[ p(y_k) \frac{\partial^2 w}{\partial y_k^2}(y,t)
+ q(y_k) \frac{\partial w}{\partial y_k}(y,t) \right] - r w(y,t) 
\end{equation}
whenever $y \in (0,1)^d$, $t\in (0,T]$ with
\begin{equation*}
p(\eta) = \frac{1}{2\pi^2} \sin^4\! \left( \pi \eta \right),\quad
q(\eta) = \frac{1}{\pi} \sin^3\! \left( \pi \eta \right) \cos \left( \pi \eta \right)\quad
{\rm for}~\eta\in \R.
\end{equation*}
The PDE \eqref{eq:PDE_y} is a convection-diffusion-reaction equation without mixed derivatives.
Let $\psi$ denote the transform of the payoff function $\phi$, 
\begin{equation}\label{eq:psi}
\psi (y,t) = \phi \left( K \exp \left[ Qx +b(t) \right] \right) ~~{\rm with}~~  x = \tan \left[ \pi(y-\tfrac{1}{2}) \right]
\end{equation}
whenever $y \in (0,1)^d$, $t\in [0,T]$.
Then for \eqref{eq:PDE_y} one has the initial condition
\begin{equation}\label{eq:IC_y}
w(y,0) = \psi (y,0). 
\end{equation}
At the boundary $\partial D$ of the spatial domain $D=(0,1)^d$ we shall consider a Dirichlet condition.
As in \cite{HS21}, we make the minor assumption in this paper that each column of the matrix $Q$ satisfies 
one of the following two conditions: 
\begin{itemize}
\item[(a)] all its entries are strictly positive;
\item[(b)] it has both a strictly positive and a strictly negative entry.
\end{itemize}
For any given $k \in \{1,2,\ldots,d\}$ such that the $k$-th column of $Q$ satisfies condition (a) there holds 
\begin{equation}\label{eq:BC_y}
w(y,t) = K e^{-rt}
\end{equation}
whenever $y \in \partial D$ with $y_k = 0$ and $t\in (0,T]$.
On the complementary part of $\partial D$ a homogeneous Dirichlet condition is valid.
For a short proof of this result, see \cite{HS21}.

\subsection{PCA-based approximation for European basket option}\label{sec:approximation}
Assume the eigenvalues of the covariance matrix $\Sigma$ are ordered such that 
$\lambda_1 \geq \lambda_2 \geq \cdots \geq \lambda_d \geq 0$. 
In many financial applications it holds that $\lambda_1$ is dominant, that is, $\lambda_1$ is much larger than $\lambda_2$.
In view of this observation, Reisinger \& Wittum \cite{RW07} introduced a PCA-based approximation of the exact solution $w$ 
to the $d$-dimensional PDE \eqref{eq:PDE_y}.
To this purpose, regard $w$ also as a function of the eigenvalues and write $w(y,t;\lambda)$ with 
$\lambda = (\lambda_1,\lambda_2,\ldots,\lambda_d)^\rT$.
Let
\begin{equation*}
\widehat{\lambda} = (\lambda_1,0,\ldots,0)^\rT ~~ {\rm and} ~~ \delta \lambda = \lambda - \widehat{\lambda} = (0,\lambda_2,\ldots,\lambda_d)^\rT.
\end{equation*}
Under sufficient smoothness, a first-order Taylor expansion of $w$ at $\widehat{\lambda}$ yields
\begin{equation}\label{eq:TaylorExpansion}
w(y, t; \lambda) \approx w(y, t; \widehat{\lambda}) + \sum_{l=2}^d \delta \lambda_{l}\, \pderiv{w}{\lambda_{l}}(y, t; \widehat{\lambda}).
\end{equation}
The partial derivative ${\partial w}/{\partial \lambda_{l}}$ (for $2\le l \le d$) can be approximated by a forward finite difference,
\begin{equation}\label{eq:FiniteDifference}
\pderiv{w}{\lambda_{l}}(y, t; \widehat{\lambda}) \approx \frac{w(y, t; \widehat{\lambda} + \delta \lambda_l\, e_l) - w(y, t; \widehat{\lambda})}{\delta \lambda_l}\,,
\end{equation}
where $e_l$ denotes the $l$-th standard basis vector in $\R^d$.
From \eqref{eq:TaylorExpansion} and \eqref{eq:FiniteDifference}, it follows that
\begin{equation*}\label{eq:solutionApprox1}
w(y, t; \lambda) \approx w(y, t; \widehat{\lambda}) + \sum_{l=2}^d \left[ w(y, t; \widehat{\lambda} + \delta \lambda_{l}\, e_{l}) - w(y, t; \widehat{\lambda}) \right].
\end{equation*}
Write 
\begin{equation*}
w^{(1)}(y, t) = w(y, t; \widehat{\lambda}) ~~{\rm and}~~ w^{(1,\, l)}(y, t) = w(y, t; \widehat{\lambda} + \delta \lambda_{l}\, e_{l}). 
\end{equation*}
Then the {\it PCA-based approximation} reads
\begin{equation}\label{eq:solutionApprox2}
w(y, t) \approx {\widetilde w}(y, t) = w^{(1)}(y, t) + \sum_{l=2}^d \left[ w^{(1,\, l)}(y, t) - w^{(1)}(y, t) \right]
\end{equation}
whenever $y \in (0,1)^d$ and $t\in (0,T]$.
By construction, $w^{(1)}$ satisfies the PDE \eqref{eq:PDE_y} with $\lambda_k$ being set to zero for all $k\not= 1$, and $w^{(1,\, l)}$ 
satisfies \eqref{eq:PDE_y} with $\lambda_k$ being set to zero for all $k\not\in \{1,l\}$.
This is completed by the same initial and boundary conditions as for $w$, given above.
We write 
\[
{\widetilde u}(s,t) = {\widetilde w}(y(x(s,t)),t)
\] 
for the PCA-based approximation in the original coordinates.

In financial practice, one is often interested in the option value at inception in the single point $s=S_0$, where 
$S_0 = (S_0^1, S_0^2,\ldots, S_0^d)^\rT$ is the vector of initial (spot) asset prices.
Let 
\begin{equation*}
Y_0 = y(x(S_0,T))\in (0,1)^d 
\end{equation*}
denote the corresponding point in the $y$-domain.
Then $w^{(1)}(Y_0,T)$ can be acquired by solving a one-dimensional PDE on the line segment $L_1$ in the $y$-domain 
that is parallel to the $y_1$-axis and passes through $y=Y_0$.
Hence, $y_k$ can be fixed at the value $Y_{0,k}$ whenever $k\not= 1$.
Next, $w^{(1,\, l)}(Y_0,\! T)$ (for $2\le l \le d$) can be acquired by solving a two-dimensional PDE on the plane 
segment $P_l$ in the $y$-domain that is parallel to the $(y_1,y_l)$-plane and passes through $y=Y_0$.
Hence, in this case, $y_k$ can be fixed at the value $Y_{0,k}$ whenever $k\not\in \{1,l\}$.
Determining the PCA-based approximation ${\widetilde w} (Y_0,T) = {\widetilde u}(S_0,T)$ thus requires solving 
just 1 one-dimensional PDE and $d-1$ two-dimensional PDEs. 
This clearly constitutes a major computational advantage, compared to solving the full $d$-dimensional 
PDE at once whenever $d$ is medium or large. 
Notice further that the different terms in the approximation \eqref{eq:solutionApprox2} can be computed in parallel 
independently of each other.
Then the total computational cost equals that of solving just 1 two-dimensional PDE.

A rigorous error analysis of the PCA-based approximation relevant to European basket options has been given by 
Reisinger \& Wissmann \cite{RW17}.
In particular, under mild assumptions, these authors showed that $w-{\widetilde w} = {\cal O}\left( \lambda_2^2 \right)$ 
in the maximum norm.

\subsection{PCA-based approximation for American basket option}\label{sec:approximation_Am}
Applying the coordinate transformation from Section \ref{sec:transformation} to the PDCP \eqref{eq:PDCP_u} for the value function 
$u$ of an American basket option, directly yields the following PDCP for the transformed function~$w$,
\begin{subeqnarray}\label{eq:PDCP_w}
&&w(y,t) \ge \psi (y,t),\\
\nonumber\\
&&\pderiv{w}{t}(y,t) \ge {\cal B}w(y,t),\\
\nonumber\\
&&\left( w(y,t)-\psi (y,t) \right) \left( \pderiv{w}{t}(y,t) - {\cal B}w(y,t) \right) = 0
\end{subeqnarray}
\noindent
whenever $y \in (0,1)^d$, $t\in (0,T]$ with function $\psi$ defined by \eqref{eq:psi}.
As for European options, a Dirichlet condition is taken at the boundary of the spatial domain $D$.
For any given $k \in \{1,2,\ldots,d\}$ such that the entries of the $k$-th column of $Q$ are all strictly positive there holds
\begin{equation}\label{eq:BC_ynew}
w(y,t) = K 
\end{equation}
whenever $y \in \partial D$ with $y_k = 0$ and $t\in (0,T]$.
Notice that, compared to \eqref{eq:BC_y}, the discount factor $\exp(-rt)$ is absent in \eqref{eq:BC_ynew}.
On the complementary part of $\partial D$, a homogeneous Dirichlet condition is valid.

The PCA-based approximation for the American basket option value function $w$ is given by \eqref{eq:solutionApprox2},
where by definition $w^{(1)}$ satisfies the PDCP \eqref{eq:PDCP_w} with $\lambda_k$ being set to zero for all $k\not= 1$, 
and $w^{(1,\, l)}$ satisfies \eqref{eq:PDCP_w} with $\lambda_k$ being set to zero for all $k\not\in \{1,l\}$.

\setcounter{equation}{0}
\section{Discretisation}\label{sec:discretisation}

\subsection{Discretisation for European basket option}\label{sec:EU_discretisation}
To numerically obtain the values $w^{(1)}(Y_0,T)$ and $w^{(1,\, l)}(Y_0,T)$ (for $2\le l \le d$) in the approximation 
${\widetilde w}(Y_0,T)$ of $w(Y_0,T)$ for a European basket option, we adopt the finite difference discretisation of 
the pertinent one- and two-dimensional PDEs on a Cartesian nonuniform smooth spatial grid constructed in~\cite{HS21}.

Let $m\ge 1$ be the number of (unidirectional) spatial grid points in the interval $(0,1)$.  
Semidiscretisation of the PDE for $w^{(1,\, l)}$ on the plane segment $P_l$ as described in~\cite{HS21} yields a system of 
ordinary differential equations (ODEs) of the form
\begin{equation}\label{eq:ODE_y}
W^\prime(t) = \left( A_1 +  A_l \right) W(t) \,+\,g(t)
\end{equation}
for $t\in (0,T]$.
Here $W(t)$ is a vector of dimension $m^2$ and $A_1$, $A_l$ are given $m^2 \times m^2$ matrices that are tridiagonal 
(possibly up to permutation) and commute and correspond to, respectively, the first and the $l$-th spatial direction.
Next, $g(t)=g_1(t)+g_l(t)$ is a given vector of dimension $m^2$, which is obtained from the Dirichlet boundary condition 
stated at the end of Section~\ref{sec:transformation}.
The ODE system \eqref{eq:ODE_y} is provided with an initial condition 
\[
W(0) = W_0,
\] 
where the vector $W_0$ is determined by $\psi(\cdot,0)$ on $P_l$, with function $\psi$ defined by \eqref{eq:psi}.

Since the payoff function $\phi$ given by \eqref{eq:payoff} is continuous but not everywhere differentiable, the same holds 
for $\psi$.
It is well-known that this nonsmoothness can have an adverse impact on the convergence of the spatial discretisation.
To mitigate this, we apply cell averaging near the points of nonsmoothness of $\psi(\cdot,0)$ in defining the initial vector 
$W_0$, see e.g.~\cite{H17}.

For the temporal discretisation of the ODE system \eqref{eq:ODE_y}, a common Alternating Direction Implicit (ADI) method is used.
Let a step size $\Delta t = T/N$ with integer $N\ge 1$ be given and define temporal grid points $t_n = n \Delta t$ for $n=0,1,\ldots,N$.
Then the familiar second-order Brian and Douglas ADI scheme for two-dimensional PDEs yields approximations $W_n\approx W(t_n)$ 
that are successively defined for $n=1,2,\ldots,N$ by
\begin{equation}\label{eq:Douglas}
\left\{
\begin{array}{l}
Z_0 = W_{n-1}+\Delta t \left( A_1 +  A_l \right) W_{n-1} \,+\, \Delta t\, g(t_{n-1}), \\\\
Z_1 = Z_0+\tfrac{1}{2}\Delta t\, A_1 (Z_1-W_{n-1}) \,+\, \tfrac{1}{2} \Delta t\, (g_1(t_n)\,-\, g_1(t_{n-1})), \\\\
Z_2 = Z_1+\tfrac{1}{2}\Delta t\, A_l\, (Z_2-W_{n-1}) \,+\, \tfrac{1}{2} \Delta t\, (g_l(t_n)\,-\, g_l(t_{n-1})), \\\\
W_n = Z_2.
\end{array}\right.
\end{equation}
The two linear systems in each time step can be solved very efficiently by employing a priori $LU$ factorisations of the pertinent 
two matrices.
The number of floating-point operations per time step is then directly proportional to the number of spatial grid points $m^2$, 
which is optimal.

As for the spatial discretisation, also the convergence of the temporal discretisation can be adversely affected by the nonsmooth 
payoff function.
To alleviate this, we apply backward Euler damping at the initial time, also known as Rannacher time stepping, that 
is, the first time step is replaced by two half steps with the backward Euler method, see e.g.~\cite{H17}.

Discretisation of the PDE for $w^{(1)}$ on the line segment $L_1$ is done analogously to the above.
Then a semidiscrete system 
\begin{equation}\label{eq:ODE_y1}
W^\prime(t) = A_1  W(t) \,+\, g_1(t)
\end{equation}
is obtained with $W(t)$ and $g_1(t)$ vectors of dimension $m$ and $A_1$ an $m \times m$ tridiagonal matrix.
Temporal discretisation is performed by the Crank--Nicolson scheme with backward Euler damping.
Recall that the Crank--Nicolson scheme can be regarded as a special case of the Brian and Douglas scheme, which 
is seen upon setting $A_l$ and $g_l$ both equal to zero in \eqref{eq:Douglas}.

Favourable rigorous stability results for the spatial and temporal discretisations discussed in this section have been 
proved in~\cite{HS21}.

\subsection{Discretisation for American basket option}\label{sec:AM_discretisation}
Semidiscretisation of the pertinent one- and two-dimensional PDCPs in the case of American basket options follows along the same 
lines as described in Section \ref{sec:EU_discretisation} for the corresponding PDEs in the case of European basket options.
The relevant boundary condition \eqref{eq:BC_ynew} is now independent of time, and hence, the same holds for $g$.
Semidiscretisation of the PDCP for $w^{(1,\, l)}$ on the plane segment $P_l$ yields 
\begin{subeqnarray}\label{eq:semiPDCP_w}
&&W(t) \ge \Psi (t),\\
\nonumber\\
&&W^\prime(t) \ge \left( A_1 +  A_l \right) W(t) \,+\,g,\\
\nonumber\\
&&\left( W(t) - \Psi (t) \right)^\rT \left( W^\prime(t) - \left( A_1 +  A_l \right) W(t) \,-\,g \right) = 0
\end{subeqnarray}
\noindent
for $t\in (0,T]$ and $W(0) = W_0$.
Here $\Psi(t)$ is a vector of dimension $m^2$ that is determined by the function $\psi(\cdot,t)$ on $P_l$.
Inequalities for vectors are to be understood componentwise.

For the temporal discretisation of the semidiscrete PDCP \eqref{eq:semiPDCP_w} we consider two adaptations of the Brian 
and Douglas ADI scheme \eqref{eq:Douglas}.
They both generate successive approximations $\widehat{W}_n$ to $W(t_n)$ for $n=1,2,\ldots,N$ with $\widehat{W}_0 = W_0$.

The first adaptation is elementary and follows the so-called explicit payoff (EP) approach,
\begin{equation}\label{eq:Douglas_EP}
\left\{
\begin{array}{l}
Z_0 = \widehat{W}_{n-1}+\Delta t \left( A_1 +  A_l \right) \widehat{W}_{n-1} \,+\, \Delta t\, g, \phantom{\,+\, \Delta t\, \widehat{\mu}_{n-1}}\\\\
Z_1 = Z_0+\tfrac{1}{2}\Delta t\, A_1 (Z_1-\widehat{W}_{n-1}), \\\\
Z_2 = Z_1+\tfrac{1}{2}\Delta t\, A_l\, (Z_2-\widehat{W}_{n-1}), \\\\
\widebar{W}_n = Z_2,\\\\
\widehat{W}_n = \max \{ \widebar{W}_n , \Psi_n \}.
\end{array}\right.
\end{equation}
Here $\Psi_n = \Psi(t_n)$ and the maximum of two vectors is to be taken componentwise.
The adaptation \eqref{eq:Douglas_EP} can be regarded as first carrying out a time step by ignoring the American 
constraint and next applying this constraint explicitly. 

The second adaptation is more advanced and employs the Ikonen--Toivanen (IT) splitting technique~\cite{HH15,IT04,IT09},
\begin{equation}\label{eq:Douglas_IT}
\left\{
\begin{array}{l}
Z_0 = \widehat{W}_{n-1}+\Delta t \left( A_1 +  A_l \right) \widehat{W}_{n-1} \,+\, \Delta t\, g \,+\, \Delta t\, \widehat{\mu}_{n-1}, \\\\
Z_1 = Z_0+\tfrac{1}{2}\Delta t\, A_1 (Z_1-\widehat{W}_{n-1}), \\\\
Z_2 = Z_1+\tfrac{1}{2}\Delta t\, A_l\, (Z_2-\widehat{W}_{n-1}), \\\\
\widebar{W}_n = Z_2,\\\\
\widehat{W}_n = \max\left\{\widebar{W}_n- \Delta t\, \widehat{\mu}_{n-1}\,,\,\Psi_n\right\},\\\\
\widehat{\mu}_n~ = \max\left\{0\,,\,\widehat{\mu}_{n-1}+(\Psi_n-\widebar{W}_n)/\Delta t \right\}
\end{array}\right.
\end{equation}
with $\widehat{\mu}_0 =0$.
The auxiliary vector $\widehat{\mu}_n$ is often called a Lagrange multiplier.
For a useful interpretation of this adaptation we refer to \cite{HV17}.
The vector $\widehat{W}_n$ and the auxiliary vector $\widehat{\mu}_n$ are computed in two parts.
In the first part, an intermediate vector $\widebar{W}_n$ is computed.
In the second part, $\widebar{W}_n$ and $\widehat{\mu}_{n-1}$ are updated to $\widehat{W}_n$ and $\widehat{\mu}_n$ 
by a certain simple, explicit formula.

The obtained accuracy for the adaptation by the IT approach is generally better than by the EP approach, 
see e.g.~\cite{H17,HV17} and also Section~\ref{sec:numericalexperiments} below.
A virtue of both adaptations \eqref{eq:Douglas_EP}, \eqref{eq:Douglas_IT} is that the computational cost per time step 
is essentially the same as that for \eqref{eq:Douglas}.

\setcounter{equation}{0}
\section{Comonotonic approach}\label{sec:comon}
In a variety of papers in the literature, the concept of comonotonicity has been employed for arriving at 
efficiently computable approximations as well as upper and lower bounds for option values.
For European-style basket options, relevant references to the comonotonic approach are, notably,
Kaas et al.~\cite{KDG00}, 
Dhaene et al.~\cite{DDGKV02a,DDGKV02b}, 
Deelstra et al.~\cite{DLV04,DDV08} and
Chen et al.~\cite{CDDV08,CDDLV15}.
Recently, an extension to American-style basket options has been considered by Hanbali \& Linders \cite{HL19}.
In this section we review results obtained with the comonotonic approach and applied in loc.~cit.
Here the assumption has been made that the payoff function $\phi$ is convex, which is satisfied by \eqref{eq:payoff}, 
and that all correlations in the SDE system \eqref{eq:SDE} are nonnegative.

It follows from \cite{KDG00} that an upper bound for the European basket option value function $u$ is acquired 
by setting all correlations in \eqref{eq:SDE} equal to one, i.e., $\rho_{ij}=1$ for all $i,j$.
Denote this upper bound by $u^{\rm up}$.
Consider the same coordinate transformations as in Section~\ref{sec:transformation} and denote the obtained
transformed functions by $v^{\rm up}$ and $w^{\rm up}$. 
The pertinent covariance matrix $\Sigma^{\rm up} = \left(\sigma_i \sigma_j \right)$ has single nonzero eigenvalue 
$\lambda^{\rm up} = \sum_{i=1}^d \sigma_i^2$.
Hence, the function $v^{\rm up}$ satisfies the one-dimensional PDE
\begin{equation}\label{eq:PDE_x_up}
\pderiv{v^{\rm up}}{t}(x,t) = \frac{1}{2} \lambda^{\rm up} \frac{\partial^2 v^{\rm up}}{\partial x_1^2}(x,t) - r v^{\rm up}(x,t)
\end{equation}
whenever $x \in \R^d$, $t\in (0,T]$.
Next, the function $w^{\rm up}$ satisfies the one-dimensional PDE
\begin{equation}\label{eq:PDE_y_up}
\pderiv{w^{\rm up}}{t}(y,t) =
{\cal B^{\rm up}}w^{\rm up}(y,t) :=
\lambda^{\rm up} \left[ p(y_1) \frac{\partial^2 w^{\rm up}}{\partial y_1^2}(y,t)
+ q(y_1) \frac{\partial w^{\rm up}}{\partial y_1}(y,t) \right] - r w^{\rm up}(y,t) 
\end{equation}
whenever $y \in (0,1)^d$, $t\in (0,T]$.
The same initial and boundary conditions apply as in Section~\ref{sec:transformation}, using the pertinent function
$\psi^{\rm up}$.

It turns out that the upper bound above is, in general, rather crude.
In the comonotonic approach, accurate lower bounds for the European basket option value have been derived, however.
We consider here the lower bound chosen in \cite{HL19}, which has been motivated by results obtained in \cite{KDG00,DLV04}.
Let $\nu_i \in (0,1]$ be given by
\begin{equation}
\nu_i = \frac{\sum_{j=1}^d \omega_j S^j_0 \rho_{ij} \sigma_j}{\sqrt{\sum_{j=1}^d 
\sum_{k=1}^d \omega_j \omega_k S^j_0 S^k_0 \rho_{jk} \sigma_j \sigma_k}} \quad{\rm for}~1\le i \le d.
\end{equation}
The lower bound is acquired upon replacing the volatility $\sigma_i$ by $\nu_i \sigma_i$ for $1\le i \le d$
and subsequently setting in \eqref{eq:SDE} all correlations equal to one.
Denote this bound by $u^{\rm low}$ and the corresponding transformed functions by $v^{\rm low}$ and $w^{\rm low}$.
Then, with $\lambda^{\rm low} = \sum_{i=1}^d (\nu_i\sigma_i)^2$, the function $v^{\rm low}$ satisfies the one-dimensional PDE
\begin{equation}\label{eq:PDE_x_low}
\pderiv{v^{\rm low}}{t}(x,t) = \frac{1}{2} \lambda^{\rm low} \frac{\partial^2 v^{\rm low}}{\partial x_1^2}(x,t) - r v^{\rm low}(x,t)
\end{equation}
whenever $x \in \R^d$, $t\in (0,T]$.
Next, the function $w^{\rm low}$ satisfies the one-dimensional PDE
\begin{equation}\label{eq:PDE_y_low}
\pderiv{w^{\rm low}}{t}(y,t) =
{\cal B^{\rm low}}w^{\rm low}(y,t) :=
\lambda^{\rm low} \left[ p(y_1) \frac{\partial^2 w^{\rm low}}{\partial y_1^2}(y,t)
+ q(y_1) \frac{\partial w^{\rm low}}{\partial y_1}(y,t) \right] - r w^{\rm low}(y,t) 
\end{equation}
whenever $y \in (0,1)^d$, $t\in (0,T]$.
The same initial and boundary conditions apply as in Section~\ref{sec:transformation}, using the pertinent function
$\psi^{\rm low}$.

Clearly, the comonotonic upper as well as lower bound can be viewed as obtained upon replacing in the PDE 
\eqref{eq:PDE_s} the covariance matrix $\Sigma$ by a certain matrix of rank one.
For the lower bound, this rank-one matrix is given by $\Sigma^{\rm low} = \xi\, \xi^\rT$ with (eigen)vector 
$\xi = (\nu_1\sigma_1,\nu_2\sigma_2,\ldots,\nu_d\sigma_d)^\rT$ and single nonzero eigenvalue 
$\lambda^{\rm low} = \xi^\rT \xi$.

Based on a result by Vyncke et al.~\cite{VGD04}, a specific linear combination of the comonotonic lower and upper 
bounds has been considered in \cite{HL19}, which approximates the value of a European basket option.
This {\it comonotonic approximation reads 
\begin{equation}\label{eq:com_w}
u^{\rm app}(S_0,T) = z u^{\rm low}(S_0,T) + (1-z) u^{\rm up}(S_0,T),
\end{equation} 
}
where $z \ge 0$ is given by
\begin{equation*}
z = \frac{c - b}{c - a}
\end{equation*}
with
\begin{equation*}
\begin{split}
a &= \sum_{i=1}^d\sum_{j=1}^d \omega_i \omega_j S_0^i S_0^j \left( e^{\nu_i \nu_j \sigma_i \sigma_j T} - 1 \right),\\
b &= \sum_{i=1}^d\sum_{j=1}^d \omega_i \omega_j S_0^i S_0^j \left( e^{\rho_{ij}   \sigma_i \sigma_j T} - 1 \right), \\
c &= \sum_{i=1}^d\sum_{j=1}^d \omega_i \omega_j S_0^i S_0^j \left( e^{            \sigma_i \sigma_j T} - 1 \right).
\end{split}
\end{equation*}

In \cite{HL19} the authors next proposed \eqref{eq:com_w} as an approximation to the value of an American
basket option, where $u^{\rm low}$ and $u^{\rm up}$ are now defined via the solutions $w^{\rm low}$ and 
$w^{\rm up}$ to the PDCP \eqref{eq:PDCP_w} with ${\cal B}$ replaced by ${\cal B^{\rm low}}$ and ${\cal B^{\rm up}}$, 
respectively, and function $\psi$ replaced by $\psi^{\rm low}$ and $\psi^{\rm up}$, respectively.
We remark that, to our knowledge, it is an open question in the literature at present whether these functions 
$u^{\rm low}$ and $u^{\rm up}$ form actual lower and upper bounds for the American basket option value.

For the numerical solution of the pertinent PDEs and PDCPs, in \cite{HL19} a finite difference method was applied 
in space and the explicit Euler method in time, with the EP approach for American basket options.
In the following, we shall employ the spatial and temporal discretisations described in Section~\ref{sec:discretisation}.
In particular this allows for much less time steps than is required, in view of stability, by the explicit Euler method.

\setcounter{equation}{0}
\section{Numerical experiments}\label{sec:numericalexperiments}

In this section we perform ample numerical experiments.
Our main aims are to determine whether the PCA-based and comonotonic approaches define approximations to European and 
American basket put option values that lie close to each other, and next, to gain insight into the error of the 
discretisations described in Section~\ref{sec:discretisation} in computing these approximations. 
 
We consider two parts of experiments, depending on the parameter sets chosen for the basket option and underlying 
asset price model.
In the first part we choose the same six parameter sets A--F as considered in \cite{HS21}.
In the second part we shall select parameter sets similar to those in~\cite{HL19}.

Commencing with the first part, Set~A is taken from Reisinger \& Wittum \cite{RW07}. 
Here $d=5$, $K=1$, $T=1$, $r=0.05$ and
\[
(\rho_{ij})_{i,j=1}^d = 
\begin{pmatrix}
1.00 & 0.79 & 0.82 & 0.91 & 0.84\\
0.79 & 1.00 & 0.73 & 0.80 & 0.76\\
0.82 & 0.73 & 1.00 & 0.77 & 0.72\\
0.91 & 0.80 & 0.77 & 1.00 & 0.90\\
0.84 & 0.76 & 0.72 & 0.90 & 1.00
\end{pmatrix},\qquad~~
\]
\[
(\sigma_i)_{i=1}^d = 
\begin{pmatrix}
0.518 & 0.648 & 0.623 & 0.570 & 0.530
\end{pmatrix},
\]
\[
(\omega_i)_{i=1}^d = 
\begin{pmatrix}
0.381 & 0.065 & 0.057 & 0.270 & 0.227
\end{pmatrix}.
\]
The corresponding covariance matrix $\Sigma$ has eigenvalues 
\[
(\lambda_i)_{i=1}^d = \begin{pmatrix} 1.4089 & 0.1124 & 0.1006 & 0.0388 & 0.0213 \end{pmatrix}
\]
and it is clear that $\lambda_1$ is dominant.

Sets~B and C are obtained from Jain \& Oosterlee \cite{JO15} and possess dimensions $d=10$ and $d=15$, respectively. 
Here $K=40$, $T=1$, $r=0.06$ and $\rho_{ij}=0.25$, $\sigma_i = 0.20$, $\omega_i = 1/d$ for $1\le i\not= j \le d$.
Sets B and C have $\lambda_1 = 0.13$ and $\lambda_1 = 0.18$, respectively, and $\lambda_2 = \ldots = \lambda_d = 0.03$.
Hence, $\lambda_1$ is also dominant for these parameter sets.

Sets D, E, F possess dimensions $d=5, 10, 15$, respectively, where $K=100$, $T=1$, $r=0.04$ and 
$\rho_{ij}=\exp(-\mu |i-j|)$, $\sigma_i = 0.30$, $\omega_i = 1/d$ for $1\le i, j \le d$ with $\mu=0.0413$.
The relevant correlation structure has been considered in for example Reisinger \& Wissmann~\cite{RW15} and yields
eigenvalues that decrease rapidly.
Sets D, E, F have in particular 
\begin{equation*}
(\lambda_1~\lambda_2~\lambda_3) = 
\begin{pmatrix} 0.4218 & 0.0180 & 0.0053\end{pmatrix},~ 
\begin{pmatrix} 0.7897 & 0.0647 & 0.0187\end{pmatrix},~ 
\begin{pmatrix} 1.1126 & 0.1337 & 0.0402\end{pmatrix},
\end{equation*}
respectively.

It can be verified that for all Sets~A--F the pertinent matrix of eigenvectors $Q$ satisfies the assumption from
Section~\ref{sec:transformation}.

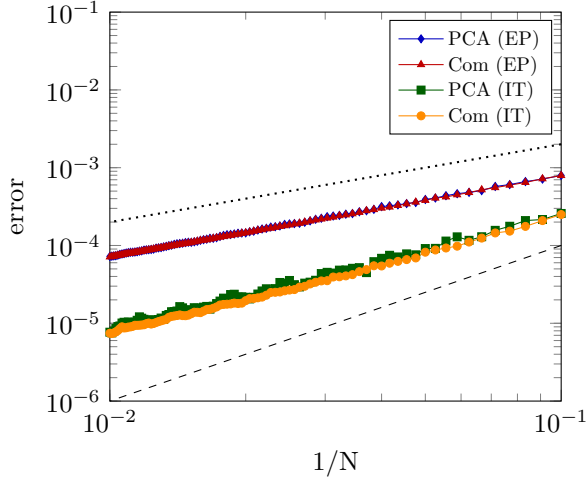
\begin{figure}[h!]
	\centering
	\input{timeerrors-BSamND-PutOnAvg-set3-5D-r=0.05-yDomain-non-uniform-N=10to100.tikz}
	\caption{Error with respect to the semidiscrete values for ${\widetilde u}(S_0,T)$ and $u^{\rm app}(S_0,T)$  
	if $m=100$. Two reference lines included for first-order convergence (dotted) and second-order convergence 
	(dashed).}
	\label{fig:timeerrors}
\end{figure}

Our first numerical experiment concerns the two adaptations of the temporal discretisation scheme to PDCPs by the 
EP and IT approaches as described in Section~\ref{sec:AM_discretisation} for American-style options.
Consider Set A and $S_0=(K,K,\ldots,K)^\rT$. 
For a fixed number of spatial grid points, given by $m=100$, we study the absolute error in the two pertinent 
discretisations of the PCA-based and comonotonic approximations ${\widetilde u}(S_0,T)$ and $u^{\rm app}(S_0,T)$ 
in function of the number of time steps $N=10,11,12,\ldots,100$.
Figure~\ref{fig:timeerrors} displays the obtained errors with respect to the values computed for a large number 
of time steps, $N=1000$. 
Note that these errors do not contain the error due to spatial discretisation, but only due to the 
temporal discretisation.
Figure~\ref{fig:timeerrors} clearly illustrates that, in the PCA-based as well as the comonotonic case, the IT 
approach yields a (much) smaller error than the EP approach for any given $N$.
Further, the observed order of convergence for IT is approximately 1.5, whereas for EP it is only approximately 
1.0.
The better performance of IT compared to EP is well-known in the literature, see e.g.~\cite{H17,IT09,HV17}.
Accordingly, in the following, we shall always apply the IT approach.

Let $S_0=(K,K,\ldots,K)^\rT$ as above.
Table~\ref{tbl:refsEU} displays our reference values for the PCA-based and comonotonic approximations 
${\widetilde u}(S_0,T)$ and $u^{\rm app}(S_0,T)$, respectively, as well as the lower bound $u^{\rm low}(S_0,T)$ 
for the European basket put option.
These values have been obtained by applying the PDE discretisation from Section~\ref{sec:discretisation} with 
$m=N=1000$ spatial and temporal grid points.
Clearly, the positive result holds that, for each given set, the two approximations and the lower bound lie close 
to each other.

\begin{table}[h!]
	\centering
	\begin{tabular}{c|c|c|c}
		Set & ${\widetilde u}(S_0,T)$ & $u^{\rm app}(S_0,T)$ & $u^{\rm low}(S_0,T)$ \\
		\hline
		A & 0.17577 & 0.17583 & 0.17577 \\
		B & 0.83257 & 0.84125 & 0.83942\\
		C & 0.77065 & 0.78083 & 0.77955\\
		D & 9.46550 & 9.46570 & 9.46523\\
		E & 9.10039 & 9.10128 & 9.09974\\
		F & 8.76358 & 8.76554 & 8.76255
	\end{tabular}
	\caption{Reference values ${\widetilde u}(S_0,T)$, $u^{\rm app}(S_0,T)$, $u^{\rm low}(S_0,T)$ for 
	European basket put option.}
	\label{tbl:refsEU}
\end{table}

Similarly, Table~\ref{tbl:refsAM} shows our reference values for ${\widetilde u}(S_0,T)$, $u^{\rm app}(S_0,T)$, 
$u^{\rm low}(S_0,T)$ for the American basket put option.
These values have been obtained by applying the PDCP discretisation from Section~\ref{sec:discretisation}
and $m=N=1000$.
We find the favourable result that also in the American case, for each given set, the PCA-based and comonotonic 
approximations lie close to each other.
Recall that, at present, it is not clear whether $u^{\rm low}(S_0,T)$ forms an actual lower bound in this case.
 
\begin{table}[h!]
	\centering
	\begin{tabular}{c|c|c|c}
		Set & ${\widetilde u}(S_0,T)$ & $u^{\rm app}(S_0,T)$ & $u^{\rm low}(S_0,T)$ \\
		\hline
		A & 0.18110 & 0.18120 & 0.18114\\
		B & 1.07928 & 1.08615 & 1.08431\\
		C & 1.01641 & 1.02435 & 1.02306\\
		D & 9.86176 & 9.86206 & 9.86159\\
		E & 9.49645 & 9.49774 & 9.49620\\
		F & 9.15935 & 9.16219 & 9.15920 
	\end{tabular}
	\caption{Reference values ${\widetilde u}(S_0,T)$, $u^{\rm app}(S_0,T)$, $u^{\rm low}(S_0,T)$ 
	for American basket put option.}
	\label{tbl:refsAM}
\end{table}

We next study, for European and American basket put options and Sets A--F, the absolute error in the discretisation 
described in Section \ref{sec:discretisation} of the PCA-based and comonotonic approximations ${\widetilde u}(S_0,T)$ 
and $u^{\rm app}(S_0,T)$ in function of $m=N=10,11,12,\ldots,100$.
To determine the error of the discretisation for the PCA-based and comonotonic approximations, the corresponding 
reference values from Tables~\ref{tbl:refsEU} and \ref{tbl:refsAM} are used.

\begin{figure}
	\centering
	\input{errors-BSeuND-PutOnAvg-set3-5D-r=0.05-yDomain-non-uniform-M=10to100.tikz} ~
	\input{errors-BSamND-PutOnAvg-set3-5D-r=0.05-yDomain-non-uniform-M=10to100.tikz} \\
	\input{errors-BSeuND-PutOnAvg-set2-10D-r=0.06-yDomain-non-uniform-M=10to100.tikz} ~
	\input{errors-BSamND-PutOnAvg-set2-10D-r=0.06-yDomain-non-uniform-M=10to100.tikz} \\
	\input{errors-BSeuND-PutOnAvg-set2-15D-r=0.06-yDomain-non-uniform-M=10to100.tikz} ~
	\input{errors-BSamND-PutOnAvg-set2-15D-r=0.06-yDomain-non-uniform-M=10to100.tikz} \\
	\caption{Discretisation error for ${\widetilde u}(S_0,T)$ and $u^{\rm app}(S_0,T)$ in cases A (top), B (middle) 
	and C (bottom).
	Left: European basket option. Right: American basket option. 
	Reference line (dashed) included for second-order convergence.}
	\label{fig:errors-PutOnAvg-setABC}
\end{figure}
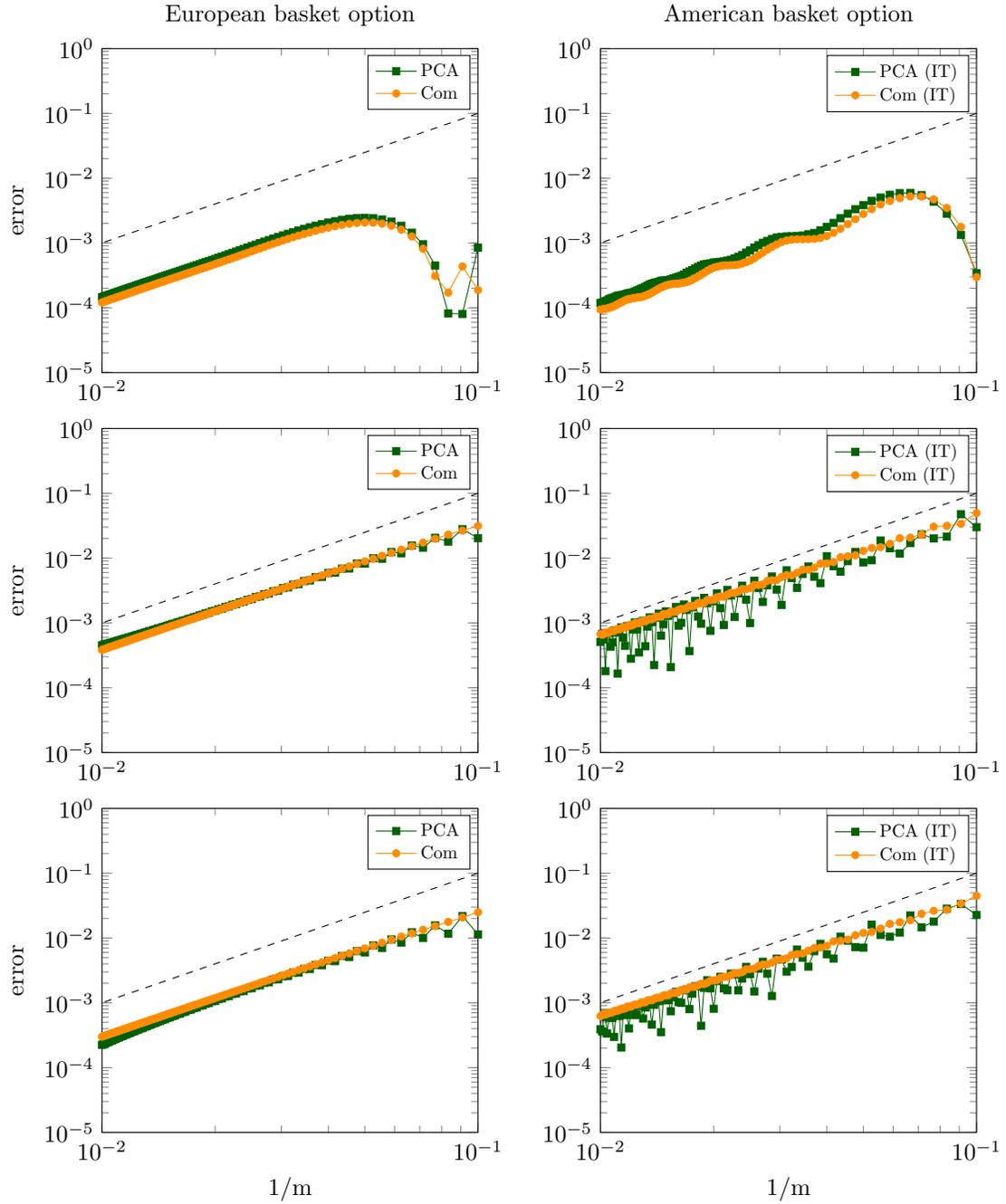

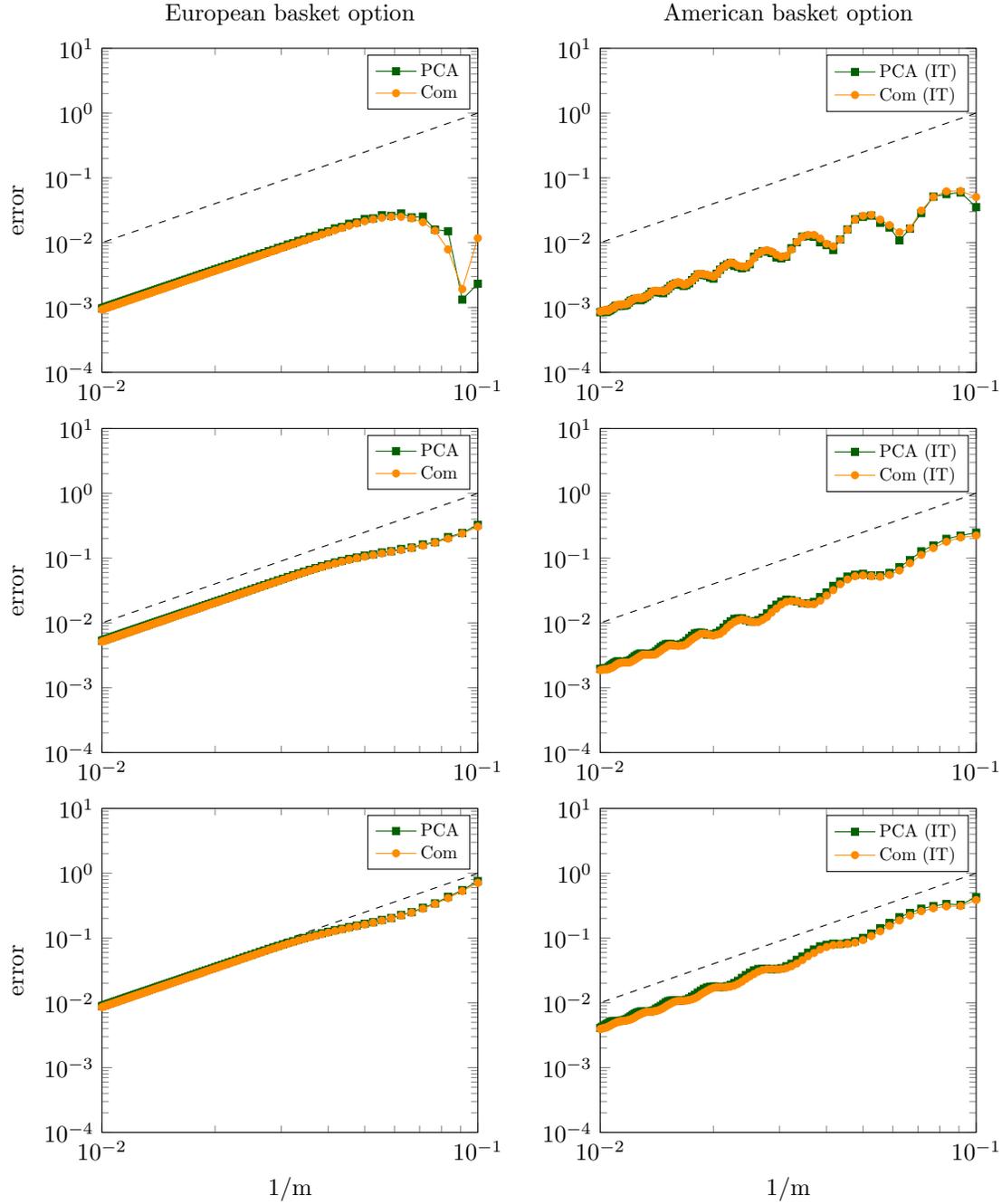
\begin{figure}
	\centering
	\input{errors-BSeuND-PutOnAvg-set5-phi=0.0413-5D-r=0.04-yDomain-non-uniform-M=10to100.tikz} ~
	\input{errors-BSamND-PutOnAvg-set5-phi=0.0413-5D-r=0.04-yDomain-non-uniform-M=10to100.tikz} \\
	\input{errors-BSeuND-PutOnAvg-set5-phi=0.0413-10D-r=0.04-yDomain-non-uniform-M=10to100.tikz} ~
	\input{errors-BSamND-PutOnAvg-set5-phi=0.0413-10D-r=0.04-yDomain-non-uniform-M=10to100.tikz} \\
	\input{errors-BSeuND-PutOnAvg-set5-phi=0.0413-15D-r=0.04-yDomain-non-uniform-M=10to100.tikz} ~
	\input{errors-BSamND-PutOnAvg-set5-phi=0.0413-15D-r=0.04-yDomain-non-uniform-M=10to100.tikz}
	\caption{Discretisation error for ${\widetilde u}(S_0,T)$ and $u^{\rm app}(S_0,T)$ in cases D (top), E (middle) 
	and F (bottom). 
	Left: European basket option. Right: American basket option. 
	Reference line (dashed) included for second-order convergence.}
	\label{fig:errors-PutOnAvg-setDEF}
\end{figure}

Figures~\ref{fig:errors-PutOnAvg-setABC} and~\ref{fig:errors-PutOnAvg-setDEF} display for Sets A, B, C and D, E, F, 
respectively, the absolute error in the discretisation of ${\widetilde u}(S_0,T)$ and $u^{\rm app}(S_0,T)$ 
versus $1/m$, where the left column concerns the European option and the right column the American option.

As a main observation, Figures~\ref{fig:errors-PutOnAvg-setABC} and~\ref{fig:errors-PutOnAvg-setDEF} clearly indicate (near) 
second-order convergence of the discretisation error in all cases, that is, for all Sets A--F, for both the European and 
American basket options, and for both the PCA-based and comonotonic approximations.
This is a very favourable result. 
Additional experiments indicate that the error stems essentially from the spatial discretisation (and not the 
temporal discretisation).

For the European option and Sets~A and D, we remark that the error drop in the (less important) region $m\le 20$ 
corresponds to a change of sign.
Besides this, in the case of the European basket option, the behaviour of the discretisation error is always seen to be 
regular.

For the American option, it is found that the discretisation error often behaves somewhat less regular, with oscillations
occurring.
A similar phenomenon has recently been observed and studied in~\cite{HS21} for Bermudan basket options and is attributed 
to the spatial nonsmoothness of the exact option value function at the early exercise boundary.

In the following we consider the second part of experiments and choose parameter sets inspired by those from \cite{HL19}.
Here a basket put option with $d=8$ equally weighted underlying assets is taken and $S_0=(40,40,\ldots,40)^\rT$. 
Next, the strike $K\in \{35,40,45\}$ and the maturity time $T\in\{0.5,1,2\}$.
For the interest rate we choose\footnote{In \cite{HL19} the rate $r=0.01$ is taken, but then American option values 
are often close to their European counterpart, which is less interesting.} $r=0.05$ and the volatilities 
are given by 
\[
(\sigma_i)_{i=1}^8 = 
\begin{pmatrix}
\sigma_1 & 0.6 & 0.1 & 0.9 & 0.3 & 0.7 & 0.8 & 0.2
\end{pmatrix}
\]
with $\sigma_1 \in \{0.3,0.9\}$.
We select correlation $\rho_{ij}=0.8$ for all $i\not= j$.
Then, for the pertinent two covariance matrices, the first eigenvalue is dominant.
In particular, there holds
\begin{equation*}
	\begin{split}
	\sigma_1=0.3: \quad\quad &(\lambda_i)_{i=1}^8 = \begin{pmatrix} 2.1398 & 0.1461 & 0.1101 & 0.0796 & \ldots \end{pmatrix},\\
	\sigma_1=0.9: \quad\quad &(\lambda_i)_{i=1}^8 = \begin{pmatrix} 2.7299 & 0.1620 & 0.1396 & 0.1076 & \ldots \end{pmatrix}.
	\end{split}
\end{equation*}
Further, the relevant matrices of eigenvectors $Q$ satisfy the assumption from Section~\ref{sec:transformation}.

Tables~\ref{tbl:refsEU4HL} and \ref{tbl:refsAM4HL} show our reference values for ${\widetilde u}(S_0,T)$, 
$u^{\rm app}(S_0,T)$, $u^{\rm low}(S_0,T)$ for the European and American basket put option, respectively, which have 
been obtained in the same way as above.
Again, we find the favourable result that, for each given parameter set and each given (European or American)
option, these three values lie close to each other.

\begin{figure}
	\centering
	\input{BlackScholesNDEuropean-K=40-T=0.5-vol=0.3-rho=0.8-r=0.05.tikz} ~
	\input{BlackScholesNDAmerican-K=40-T=0.5-vol=0.3-rho=0.8-r=0.05.tikz} \\ \bigskip
	\input{BlackScholesNDEuropean-K=40-T=1-vol=0.3-rho=0.8-r=0.05.tikz} ~
	\input{BlackScholesNDAmerican-K=40-T=1-vol=0.3-rho=0.8-r=0.05.tikz} \\ \bigskip
	\input{BlackScholesNDEuropean-K=40-T=2-vol=0.3-rho=0.8-r=0.05.tikz} ~
	\input{BlackScholesNDAmerican-K=40-T=2-vol=0.3-rho=0.8-r=0.05.tikz}
	\caption{Discretisation error for ${\widetilde u}(S_0,T)$ and $u^{\rm app}(S_0,T)$ 
	in cases $T=0.5$ (top), $T=1$ (middle) and $T=2$ (bottom)
	where $K = 40$, $\sigma_1 = 0.3$.
	Left: European basket option. Right: American basket option. 
	Reference line (dashed) included for second-order convergence.}	
	\label{fig:errors-PutOnAvg-setHL-vol3}
\end{figure}
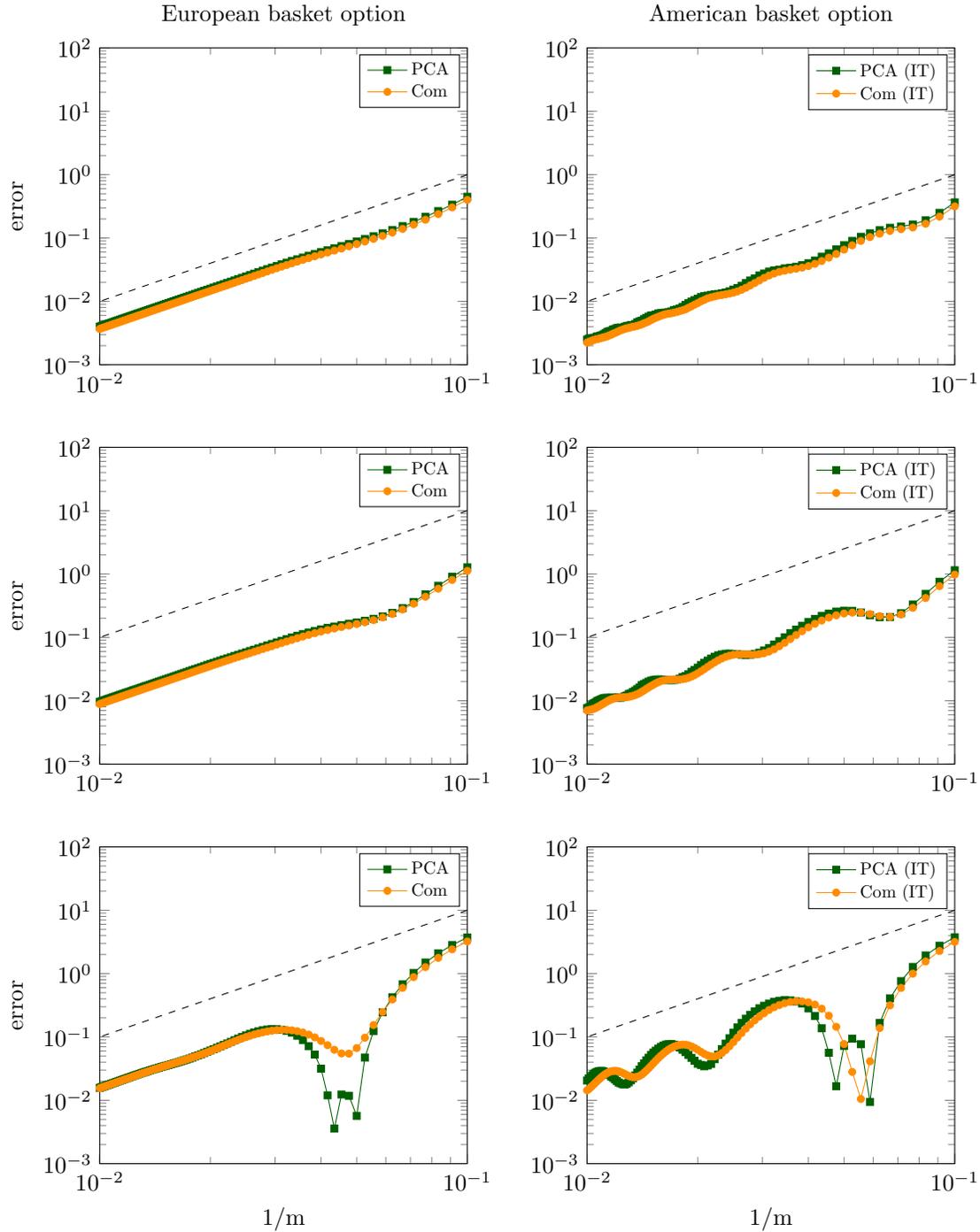

Figure~\ref{fig:errors-PutOnAvg-setHL-vol3} displays, analogously to Figures~\ref{fig:errors-PutOnAvg-setABC} 
and~\ref{fig:errors-PutOnAvg-setDEF}, the absolute error in the discretisation of ${\widetilde u}(S_0,T)$ 
and $u^{\rm app}(S_0,T)$ for the (representative) three parameter sets given by $T\in\{0.5,1,2\}$, 
$K=40$, $\sigma_1=0.3$.
The outcomes again indicate a favourable, second-order convergence result.
The regularity of the error behaviour is seen to decrease as the maturity time $T$ increases.
We note that for $T=2$ this behaviour is partly explained from a (near) vanishing error when $m\approx 20$.

\setcounter{equation}{0}
\section{Conclusions}\label{sec:conclusion}
In this paper we have studied an extension of the PCA-based approach by Reisinger \& Wittum \cite{RW07} 
to valuate American basket options. 
This approximation approach is highly effective, as the numerical solution of only a limited 
number of low-dimensional PDCPs is required.
In addition, we have considered the comonotonic approach, which was developed for basket options notably 
in 
\cite{KDG00,DDGKV02a,DDGKV02b,DLV04,DDV08,CDDV08,CDDLV15,VGD04}.
We have studied the comonotonic approximation formula for American basket option values recently examined
in Hanbali \& Linders \cite{HL19}.
This comonotonic approach is also highly effective, since it requires the numerical solution of just two 
one-dimensional PDCPs.

For the discretisation of the pertinent PDCPs, we applied finite differences on a nonuniform spatial grid 
followed by the Brian and Douglas ADI scheme on a uniform temporal grid and selected the Ikonen--Toivanen 
(IT) technique~\cite{HH15,IT04,IT09} to efficiently handle the complementarity problem in each time step.

As a first main result, we find in ample numerical experiments that the PCA-based and comonotonic 
approaches always yield approximations to the value of an American (as well as European) basket option 
that lie close to each other.

As a next main result, we observe near second-order convergence of the discretisation error in all 
numerical experiments for both the PCA-based and comonotonic approaches for American (as well as 
European) basket options.

At this moment it is still open which (if any) of the two approaches, PCA-based or 
comonotonic, is to be preferred for the approximate valuation of American basket options on $d\ge 5$
assets.
In particular, whereas in our experiments the two approaches always define approximations that lie close 
to each other, it is not clear at present which approach (if any) generally yields the smallest error 
with respect to the exact option value.
A further investigation into the PCA-based and comonotonic approaches, both experimental 
and analytical, will be the subject of future research.

\clearpage
\vfill

\clearpage
\vfill

\clearpage
\begin{table}
	\centering
	\begin{tabular}{c|c|c|r|r|r}
		$T$	& $K$ & $\sigma_1$ & ${\widetilde u}(S_0,T)$ & $u^{\rm app}(S_0,T)$ & $u^{\rm low}(S_0,T)$ \\ 
		\hline
				& 35	& 0.3	&  2.13020	&  2.13271	&  2.12954	\\
				& 		& 0.9	&  2.74982	&  2.75307	&  2.74963	\\
			0.5 & 40	& 0.3	&  4.40336	&  4.40715	&  4.40328	\\
				&   	& 0.9	&  5.14582	&  5.15003	&  5.14595	\\
				& 45	& 0.3	&  7.45442	&  7.45827	&  7.45427	\\
				&   	& 0.9	&  8.21316	&  8.21738	&  8.21313	\\ \hline
				& 35	& 0.3	&  3.35805	&  3.36599	&  3.35620	\\
				&   	& 0.9	&  4.23834	&  4.24750	&  4.23731	\\
			1 	& 40	& 0.3	&  5.78199	&  5.79261	&  5.78114	\\
				&   	& 0.9	&  6.79656	&  6.80770	&  6.79599	\\
				& 45	& 0.3	&  8.75406	&  8.76551	&  8.75329	\\
				&   	& 0.9	&  9.82315	&  9.83486	&  9.82235	\\ \hline
				& 35	& 0.3	&  4.71159	&  4.73545	&  4.70532	\\
				&   	& 0.9	&  5.89254	&  5.91682	&  5.88742	\\
			2 	& 40	& 0.3	&  7.20593	&  7.23607	&  7.20149	\\
				&   	& 0.9	&  8.54494	&  8.57378	&  8.54048	\\
				& 45	& 0.3	& 10.08246	& 10.11611	& 10.07862	\\
				&   	& 0.9	& 11.51843	& 11.54974	& 11.51371
	\end{tabular}
	\caption{Reference values ${\widetilde u}(S_0,T)$, $u^{\rm app}(S_0,T)$, $u^{\rm low}(S_0,T)$ 
	for European basket put option.}
	\label{tbl:refsEU4HL}
\end{table}

\begin{table}
	\centering
	\begin{tabular}{c|c|c|r|r|r}
		$T$	& $K$	& $\sigma_1$ & ${\widetilde u}(S_0,T)$ & $u^{\rm app}(S_0,T)$ & $u^{\rm low}(S_0,T)$ \\ 
		\hline
		& 35	& 0.3	&  2.17006	&  2.17293	&  2.16973	\\
		& 		& 0.9	&  2.79440	&  2.79840	&  2.79494	\\
		0.5 & 40	& 0.3	&  4.50018	&  4.50506	&  4.50118	\\
		&   	& 0.9	&  5.24177	&  5.24795	&  5.24387	\\
		& 45	& 0.3	&  7.64424	&  7.65063	&  7.64670	\\
		&   	& 0.9	&  8.38729	&  8.39562	&  8.39142	\\ \hline
		& 35	& 0.3	&  3.48012	&  3.48874	&  3.47879	\\
		&   	& 0.9	&  4.37236	&  4.38280	&  4.37246	\\
		1 	& 40	& 0.3	&  6.01652	&  6.02870	&  6.01717	\\
		&   	& 0.9	&  7.03281	&  7.04676	&  7.03498	\\
		& 45	& 0.3	&  9.14612	&  9.16072	&  9.14867	\\
		&   	& 0.9	& 10.19561	& 10.21256	& 10.20013	\\ \hline
		& 35	& 0.3	&  5.06452	&  5.08982	&  5.05865	\\
		&   	& 0.9	&  6.27930	&  6.30536	&  6.27500	\\
		2 	& 40	& 0.3	&  7.78521	&  7.81748	&  7.78222	\\
		&   	& 0.9	&  9.14045	&  9.17258	&  9.13855	\\
		& 45	& 0.3	& 10.94634	& 10.98327	& 10.94585	\\
		&   	& 0.9	& 12.36710	& 12.40399	& 12.36770
	\end{tabular}
	\caption{Reference values ${\widetilde u}(S_0,T)$, $u^{\rm app}(S_0,T)$, $u^{\rm low}(S_0,T)$
	for American basket put option.}
	\label{tbl:refsAM4HL}
\end{table}
\vfill

\clearpage
\vfill

\clearpage
\bibliographystyle{unsrt}
\bibliography{references}
 
\end{document}

%% file: timeerrors-BSamND-PutOnAvg-set3-5D-r=0.05-yDomain-non-uniform-N=10to100.tikz
%
%
%
\begin{tikzpicture}

\begin{axis}[%
width=6cm,
scale only axis,
xmode=log,
xmin=0.01,
xmax=0.1,
xminorticks=true,
xlabel={1/N},
ymode=log,
ymin=1e-6,
ymax=1e-1,
yminorticks=true,
ylabel={error},
legend style={legend cell align=left, align=left, nodes={scale=0.75, transform shape}}
]
\addplot [color=color3, mark=diamond*, mark options={solid, color3, scale=0.75}]
  table[row sep=crcr]{%
0.1	0.00079528706290799\\
0.0909090909090909	0.00071310662846974\\
0.0833333333333333	0.000663360780490535\\
0.0769230769230769	0.000606829157437067\\
0.0714285714285714	0.000573865819471192\\
0.0666666666666667	0.000506919154436847\\
0.0625	0.000481401040883123\\
0.0588235294117647	0.000463563571800829\\
0.0555555555555556	0.000436937686083855\\
0.0526315789473684	0.000407746659341818\\
0.05	0.000390344112803837\\
0.0476190476190476	0.000356724042617107\\
0.0454545454545455	0.000345923011594362\\
0.0434782608695652	0.000332902856891804\\
0.0416666666666667	0.000323666681182083\\
0.04	0.000316334269957697\\
0.0384615384615385	0.000291254423119253\\
0.037037037037037	0.000279519482762947\\
0.0357142857142857	0.000271543576653571\\
0.0344827586206897	0.00025735046049713\\
0.0333333333333333	0.000252583499005921\\
0.032258064516129	0.000243685141417621\\
0.03125	0.000239157896393255\\
0.0303030303030303	0.000234646683204803\\
0.0294117647058824	0.000223952744953154\\
0.0285714285714286	0.000217242259313327\\
0.0277777777777778	0.000207298837760161\\
0.027027027027027	0.000198895134149307\\
0.0263157894736842	0.000193054238280366\\
0.0256410256410256	0.000191478251943156\\
0.025	0.000189476254087023\\
0.024390243902439	0.000186521071605134\\
0.0238095238095238	0.000178571723267157\\
0.0232558139534884	0.000174374711863323\\
0.0227272727272727	0.000172672424221004\\
0.0222222222222222	0.000166304122693361\\
0.0217391304347826	0.000162906321290568\\
0.0212765957446809	0.000158131532889599\\
0.0208333333333333	0.000152673808528681\\
0.0204081632653061	0.000150172417072242\\
0.02	0.000145128822070961\\
0.0196078431372549	0.000144656842762464\\
0.0192307692307692	0.000143542241510197\\
0.0188679245283019	0.000140657639672082\\
0.0185185185185185	0.00013858306365655\\
0.0181818181818182	0.000137589694865242\\
0.0178571428571429	0.000134581033686576\\
0.0175438596491228	0.000130870085568641\\
0.0172413793103448	0.000126033085282806\\
0.0169491525423729	0.000123768899814941\\
0.0166666666666667	0.000123647115531034\\
0.0163934426229508	0.000120488259754081\\
0.0161290322580645	0.000118798641286111\\
0.0158730158730159	0.000116434402125581\\
0.015625	0.000113309161815678\\
0.0153846153846154	0.000111547527951283\\
0.0151515151515152	0.000109817621135599\\
0.0149253731343284	0.000110127587980291\\
0.0147058823529412	0.000109023286944676\\
0.0144927536231884	0.000106965797514708\\
0.0142857142857143	0.000105945011315917\\
0.0140845070422535	0.00010519958238317\\
0.0138888888888889	0.000103277838060634\\
0.0136986301369863	0.000100850599105279\\
0.0135135135135135	9.94898512524933e-05\\
0.0133333333333333	9.73330109317605e-05\\
0.0131578947368421	9.59717791645109e-05\\
0.012987012987013	9.4078579234963e-05\\
0.0128205128205128	9.32349706597724e-05\\
0.0126582278481013	9.23567807553094e-05\\
0.0125	8.97872730536431e-05\\
0.0123456790123457	8.88645266372556e-05\\
0.0121951219512195	8.83845902038272e-05\\
0.0120481927710843	8.78555102479806e-05\\
0.0119047619047619	8.71940916522262e-05\\
0.0117647058823529	8.53414302844913e-05\\
0.0116279069767442	8.45574035726326e-05\\
0.0114942528735632	8.32965481954995e-05\\
0.0113636363636364	8.23711878402689e-05\\
0.0112359550561798	8.18954345126466e-05\\
0.0111111111111111	8.18542875664074e-05\\
0.010989010989011	8.09621307156227e-05\\
0.0108695652173913	7.99800437191189e-05\\
0.010752688172043	7.87895863759436e-05\\
0.0106382978723404	7.77852117672118e-05\\
0.0105263157894737	7.67662955517834e-05\\
0.0104166666666667	7.50485911987508e-05\\
0.0103092783505155	7.38963791498126e-05\\
0.0102040816326531	7.34826702704361e-05\\
0.0101010101010101	7.27659098611055e-05\\
0.01	7.22603145251188e-05\\
};
\addlegendentry{PCA (EP)}

\addplot [color=color4, mark=triangle*, mark options={solid, color4, scale=0.75}]
  table[row sep=crcr]{%
0.1	0.000792889469385771\\
0.0909090909090909	0.000720423818121152\\
0.0833333333333333	0.000646296247044986\\
0.0769230769230769	0.000596680855246251\\
0.0714285714285714	0.000554320501098687\\
0.0666666666666667	0.000520365858747651\\
0.0625	0.000479444139460189\\
0.0588235294117647	0.000448322583771837\\
0.0555555555555556	0.000420988366887587\\
0.0526315789473684	0.000406837584243896\\
0.05	0.000380461876600985\\
0.0476190476190476	0.000361319177968694\\
0.0454545454545455	0.000342794056227558\\
0.0434782608695652	0.000327604693299566\\
0.0416666666666667	0.000311964602822956\\
0.04	0.00030235954871255\\
0.0384615384615385	0.000289340389976717\\
0.037037037037037	0.000280141701058084\\
0.0357142857142857	0.000266526688101509\\
0.0344827586206897	0.000255520542923104\\
0.0333333333333333	0.000247318976067029\\
0.032258064516129	0.000241387725102543\\
0.03125	0.000231856272964293\\
0.0303030303030303	0.000225053832268701\\
0.0294117647058824	0.000218674152614229\\
0.0285714285714286	0.000212941181441662\\
0.0277777777777778	0.000206266810125977\\
0.027027027027027	0.000199394707124406\\
0.0263157894736842	0.00019309735254025\\
0.0256410256410256	0.000189097805569777\\
0.025	0.000183581147819162\\
0.024390243902439	0.000179466237586939\\
0.0238095238095238	0.000176250465303868\\
0.0232558139534884	0.000171330606436065\\
0.0227272727272727	0.000166724931545181\\
0.0222222222222222	0.000163631180596729\\
0.0217391304347826	0.000159417534083439\\
0.0212765957446809	0.000156349942060124\\
0.0208333333333333	0.000152634505613147\\
0.0204081632653061	0.000148775920438243\\
0.02	0.000145682268491293\\
0.0196078431372549	0.000143140662776597\\
0.0192307692307692	0.000139796024430316\\
0.0188679245283019	0.00013752757214508\\
0.0185185185185185	0.000135126088894599\\
0.0181818181818182	0.000132704171007886\\
0.0178571428571429	0.000130395231333014\\
0.0175438596491228	0.000127965585517698\\
0.0172413793103448	0.000125681795651367\\
0.0169491525423729	0.000123554208362187\\
0.0166666666666667	0.000120712152114233\\
0.0163934426229508	0.000118460535785986\\
0.0161290322580645	0.000117018409114239\\
0.0158730158730159	0.00011560167187788\\
0.015625	0.000113326672326941\\
0.0153846153846154	0.000111161861403974\\
0.0151515151515152	0.000109084767924844\\
0.0149253731343284	0.000107559133938256\\
0.0147058823529412	0.000106207030813515\\
0.0144927536231884	0.000104748782180086\\
0.0142857142857143	0.00010353057142834\\
0.0140845070422535	0.000102233313182626\\
0.0138888888888889	0.000100421851952776\\
0.0136986301369863	9.90230225427491e-05\\
0.0135135135135135	9.78134082382509e-05\\
0.0133333333333333	9.62967189426234e-05\\
0.0131578947368421	9.467659526044e-05\\
0.012987012987013	9.33404028253815e-05\\
0.0128205128205128	9.20758756932927e-05\\
0.0126582278481013	9.10738651134224e-05\\
0.0125	8.98871915205046e-05\\
0.0123456790123457	8.84794569886438e-05\\
0.0121951219512195	8.75855952105487e-05\\
0.0120481927710843	8.64876626881783e-05\\
0.0119047619047619	8.52797675213457e-05\\
0.0117647058823529	8.4287533467009e-05\\
0.0116279069767442	8.32589991631161e-05\\
0.0114942528735632	8.22251528233009e-05\\
0.0113636363636364	8.14110221444009e-05\\
0.0112359550561798	8.04588674584161e-05\\
0.0111111111111111	7.97374165031062e-05\\
0.010989010989011	7.9246602270211e-05\\
0.0108695652173913	7.82731520359847e-05\\
0.010752688172043	7.71264655828052e-05\\
0.0106382978723404	7.63145285724576e-05\\
0.0105263157894737	7.52541333788048e-05\\
0.0104166666666667	7.44122795322821e-05\\
0.0103092783505155	7.36069405974649e-05\\
0.0102040816326531	7.27120844971063e-05\\
0.0101010101010101	7.1989628853214e-05\\
0.01	7.1305574692565e-05\\
};
\addlegendentry{Com (EP)}

\addplot [color=color1, mark=square*, mark options={solid, color1, scale=0.75}]
  table[row sep=crcr]{%
0.1	0.000258011862309587\\
0.0909090909090909	0.00021570867377238\\
0.0833333333333333	0.000209047924472877\\
0.0769230769230769	0.00017677043037026\\
0.0714285714285714	0.000156571206955614\\
0.0666666666666667	0.000129624977861731\\
0.0625	0.000116606332579261\\
0.0588235294117647	0.000129164454163055\\
0.0555555555555556	0.000113991128299629\\
0.0526315789473684	9.18427350184436e-05\\
0.05	9.13459622438051e-05\\
0.0476190476190476	7.59136703508689e-05\\
0.0454545454545455	7.80437154185021e-05\\
0.0434782608695652	7.12814293837016e-05\\
0.0416666666666667	7.48308684880294e-05\\
0.04	6.84277183026494e-05\\
0.0384615384615385	6.19811777715362e-05\\
0.037037037037037	4.48011325203668e-05\\
0.0357142857142857	5.19148400240366e-05\\
0.0344827586206897	5.09241115363435e-05\\
0.0333333333333333	4.84704307434358e-05\\
0.032258064516129	4.8225876250696e-05\\
0.03125	4.42191055265995e-05\\
0.0303030303030303	4.51840638967804e-05\\
0.0294117647058824	4.35635376409838e-05\\
0.0285714285714286	3.57400166298705e-05\\
0.0277777777777778	3.37838082275499e-05\\
0.027027027027027	3.28181371242464e-05\\
0.0263157894736842	2.93330658043822e-05\\
0.0256410256410256	3.13633943131997e-05\\
0.025	3.55144911448335e-05\\
0.024390243902439	3.26291115476662e-05\\
0.0238095238095238	3.34717092275172e-05\\
0.0232558139534884	2.79141886361522e-05\\
0.0227272727272727	2.57819343540966e-05\\
0.0222222222222222	2.77847400685372e-05\\
0.0217391304347826	2.60390236364594e-05\\
0.0212765957446809	2.37765697004566e-05\\
0.0208333333333333	2.15777626265956e-05\\
0.0204081632653061	2.0559883129212e-05\\
0.02	2.16072098050135e-05\\
0.0196078431372549	2.18973216271068e-05\\
0.0192307692307692	2.17369232710096e-05\\
0.0188679245283019	2.36433906707356e-05\\
0.0185185185185185	2.34073696685211e-05\\
0.0181818181818182	1.98771315443602e-05\\
0.0178571428571429	2.10560906110102e-05\\
0.0175438596491228	1.93020575238734e-05\\
0.0172413793103448	1.79266112052268e-05\\
0.0169491525423729	1.65964773251037e-05\\
0.0166666666666667	1.51997398917803e-05\\
0.0163934426229508	1.64767092947216e-05\\
0.0161290322580645	1.61280904545091e-05\\
0.0158730158730159	1.54304067366551e-05\\
0.015625	1.56209801589757e-05\\
0.0153846153846154	1.58641914768387e-05\\
0.0151515151515152	1.48382748845211e-05\\
0.0149253731343284	1.47364706926667e-05\\
0.0147058823529412	1.48606017424946e-05\\
0.0144927536231884	1.58058054892385e-05\\
0.0142857142857143	1.6379631987945e-05\\
0.0140845070422535	1.43754639913563e-05\\
0.0138888888888889	1.41419034308532e-05\\
0.0136986301369863	1.40850973145701e-05\\
0.0135135135135135	1.29712861610953e-05\\
0.0133333333333333	1.26080149702723e-05\\
0.0131578947368421	1.12324247111828e-05\\
0.012987012987013	1.16761209507932e-05\\
0.0128205128205128	1.11818640883554e-05\\
0.0126582278481013	1.07943527817567e-05\\
0.0125	1.10060303969428e-05\\
0.0123456790123457	1.1091843011457e-05\\
0.0121951219512195	1.05818124543633e-05\\
0.0120481927710843	1.04713203855078e-05\\
0.0119047619047619	1.12565769007922e-05\\
0.0117647058823529	1.1635844231167e-05\\
0.0116279069767442	1.21091493735237e-05\\
0.0114942528735632	1.09287151224546e-05\\
0.0113636363636364	1.03116148637283e-05\\
0.0112359550561798	9.86259729637839e-06\\
0.0111111111111111	9.53034826353671e-06\\
0.010989010989011	1.03953088197162e-05\\
0.0108695652173913	1.02150558101421e-05\\
0.010752688172043	1.03122569393432e-05\\
0.0106382978723404	1.0073844459918e-05\\
0.0105263157894737	9.24240799812792e-06\\
0.0104166666666667	9.05973778780877e-06\\
0.0103092783505155	8.56201359827891e-06\\
0.0102040816326531	7.92786610501195e-06\\
0.0101010101010101	7.52390966801686e-06\\
0.01	7.68033982179173e-06\\
};
\addlegendentry{PCA (IT)}

\addplot [color=color2, mark=*, mark options={solid, color2, scale=0.75}]
  table[row sep=crcr]{%
0.1	0.000249971825147993\\
0.0909090909090909	0.000208863833123935\\
0.0833333333333333	0.000176023915587104\\
0.0769230769230769	0.000152869057844313\\
0.0714285714285714	0.000145002248846754\\
0.0666666666666667	0.000124455482339608\\
0.0625	0.000110698735200476\\
0.0588235294117647	9.80063804059128e-05\\
0.0555555555555556	9.28317178048588e-05\\
0.0526315789473684	8.75450468020567e-05\\
0.05	8.19993937554442e-05\\
0.0476190476190476	6.91573023687275e-05\\
0.0454545454545455	6.58602730165825e-05\\
0.0434782608695652	6.24794417636576e-05\\
0.0416666666666667	5.91683166362711e-05\\
0.04	5.45186359009464e-05\\
0.0384615384615385	5.48469024255482e-05\\
0.037037037037037	4.89874835213344e-05\\
0.0357142857142857	4.61354271842906e-05\\
0.0344827586206897	4.21954594779506e-05\\
0.0333333333333333	4.06706796919898e-05\\
0.032258064516129	3.96002448098121e-05\\
0.03125	3.93367580196202e-05\\
0.0303030303030303	3.52660078070943e-05\\
0.0294117647058824	3.50403399898835e-05\\
0.0285714285714286	3.36845925363316e-05\\
0.0277777777777778	3.20219004030509e-05\\
0.027027027027027	2.9891823668915e-05\\
0.0263157894736842	2.90879835455338e-05\\
0.0256410256410256	2.7119455413227e-05\\
0.025	2.67392022139379e-05\\
0.024390243902439	2.61467661210602e-05\\
0.0238095238095238	2.52466945056817e-05\\
0.0232558139534884	2.50202074219552e-05\\
0.0227272727272727	2.48797874863871e-05\\
0.0222222222222222	2.29463376939243e-05\\
0.0217391304347826	2.16996483194398e-05\\
0.0212765957446809	2.12490320157521e-05\\
0.0208333333333333	2.08206927339305e-05\\
0.0204081632653061	2.0754387553279e-05\\
0.02	1.95280632362405e-05\\
0.0196078431372549	1.8151030048541e-05\\
0.0192307692307692	1.79334829354394e-05\\
0.0188679245283019	1.81848082950742e-05\\
0.0185185185185185	1.77278468319886e-05\\
0.0181818181818182	1.75954698318492e-05\\
0.0178571428571429	1.7476757482815e-05\\
0.0175438596491228	1.71726493131963e-05\\
0.0172413793103448	1.62170779427484e-05\\
0.0169491525423729	1.55899813120142e-05\\
0.0166666666666667	1.51275006216312e-05\\
0.0163934426229508	1.49956153367992e-05\\
0.0161290322580645	1.43224077372162e-05\\
0.0158730158730159	1.36979168176565e-05\\
0.015625	1.38506656094106e-05\\
0.0153846153846154	1.38850516070155e-05\\
0.0151515151515152	1.32723169051774e-05\\
0.0149253731343284	1.27270718227879e-05\\
0.0147058823529412	1.25058799102518e-05\\
0.0144927536231884	1.25469329350647e-05\\
0.0142857142857143	1.28148419465401e-05\\
0.0140845070422535	1.25692142089184e-05\\
0.0138888888888889	1.24080402174598e-05\\
0.0136986301369863	1.22595922926461e-05\\
0.0135135135135135	1.19356038222329e-05\\
0.0133333333333333	1.11945225601262e-05\\
0.0131578947368421	1.09984351348569e-05\\
0.012987012987013	1.0754332666868e-05\\
0.0128205128205128	1.0519457313285e-05\\
0.0126582278481013	1.02222336366697e-05\\
0.0125	9.93489970482764e-06\\
0.0123456790123457	9.76824232154505e-06\\
0.0121951219512195	1.00732379394453e-05\\
0.0120481927710843	9.88938747673029e-06\\
0.0119047619047619	9.57357002726034e-06\\
0.0117647058823529	9.5043294334185e-06\\
0.0116279069767442	9.44558921711636e-06\\
0.0114942528735632	9.41661808739891e-06\\
0.0113636363636364	9.43233311256897e-06\\
0.0112359550561798	9.13167272600934e-06\\
0.0111111111111111	8.9867322557724e-06\\
0.010989010989011	8.93784772706963e-06\\
0.0108695652173913	8.75761818813769e-06\\
0.010752688172043	8.70195948804442e-06\\
0.0106382978723404	9.01027949923305e-06\\
0.0105263157894737	8.55177151609654e-06\\
0.0104166666666667	8.07265345323382e-06\\
0.0103092783505155	7.73970701059556e-06\\
0.0102040816326531	7.51578545565845e-06\\
0.0101010101010101	7.3941243493858e-06\\
0.01	7.4521031102337e-06\\
};
\addlegendentry{Com (IT)}

\addplot [dotted, thick]
  table[row sep=crcr]{%
0.1	2e-3\\
0.01	2e-04\\
};

\addplot [dashed]
  table[row sep=crcr]{%
0.1	1e-4\\
0.01	1e-06\\
};

\end{axis}
\end{tikzpicture}%

%% file: errors-BSeuND-PutOnAvg-set3-5D-r=0.05-yDomain-non-uniform-M=10to100.tikz
%
%
%
\begin{tikzpicture}

\begin{axis}[%
width=0.37\textwidth,
scale only axis,
xmode=log,
xmin=0.01,
xmax=0.1,
xminorticks=true,
ymode=log,
ymin=1e-5,
ymax=1e-0,
yminorticks=true,
ylabel={error},
title={European basket option},
legend cell align={left},
legend style={nodes={scale=0.75, transform shape}}
]
\addplot [color=color1, mark=square*, mark options={solid, color1, scale=0.75}]
  table[row sep=crcr]{%
0.1	0.000847414078197245\\
0.0909090909090909	8.03109257746515e-05\\
0.0833333333333333	8.1666282706977e-05\\
0.0769230769230769	0.000449639125642914\\
0.0714285714285714	0.000954125302136066\\
0.0666666666666667	0.00143842287480661\\
0.0625	0.001841923170986\\
0.0588235294117647	0.00213037395852367\\
0.0555555555555556	0.00230221647476475\\
0.0526315789473684	0.00240366292318883\\
0.05	0.00243613013294441\\
0.0476190476190476	0.00241545408247576\\
0.0454545454545455	0.00235170132907225\\
0.0434782608695652	0.00226689118394474\\
0.0416666666666667	0.00216501687804702\\
0.04	0.00206110954058081\\
0.0384615384615385	0.00195104077926953\\
0.037037037037037	0.00184620253674922\\
0.0357142857142857	0.00174303202016871\\
0.0344827586206897	0.00164588237984656\\
0.0333333333333333	0.00155303068367177\\
0.032258064516129	0.0014656881476669\\
0.03125	0.00138437896638047\\
0.0303030303030303	0.00130704372426807\\
0.0294117647058824	0.00123585325322986\\
0.0285714285714286	0.00116809328577935\\
0.0277777777777778	0.0011059213688952\\
0.027027027027027	0.00104674720044984\\
0.0263157894736842	0.0009922305220392\\
0.0256410256410256	0.000941280749555434\\
0.025	0.000893786935057317\\
0.024390243902439	0.000849754540666825\\
0.0238095238095238	0.000808533920238469\\
0.0232558139534884	0.000770739765472567\\
0.0227272727272727	0.000735128361831844\\
0.0222222222222222	0.000702427951292245\\
0.0217391304347826	0.000671352345901405\\
0.0212765957446809	0.000642862987288856\\
0.0208333333333333	0.000616049365099125\\
0.0204081632653061	0.000591053618692189\\
0.02	0.000567717351462854\\
0.0196078431372549	0.00054580082092015\\
0.0192307692307692	0.000525374306740489\\
0.0188679245283019	0.00050595993396868\\
0.0185185185185185	0.000487942684225473\\
0.0181818181818182	0.000470652351356737\\
0.0178571428571429	0.000454476208205462\\
0.0175438596491228	0.000439050580074568\\
0.0172413793103448	0.00042444974354916\\
0.0169491525423729	0.000410574062293334\\
0.0166666666666667	0.000397338720293433\\
0.0163934426229508	0.00038481674595095\\
0.0161290322580645	0.000372739748230877\\
0.0158730158730159	0.000361423460075283\\
0.015625	0.000350375764192445\\
0.0153846153846154	0.000339994300668861\\
0.0151515151515152	0.000329947222568711\\
0.0149253731343284	0.000320402055815905\\
0.0147058823529412	0.000311242705696907\\
0.0144927536231884	0.000302432211161413\\
0.0142857142857143	0.000294073395175343\\
0.0140845070422535	0.000285942991487209\\
0.0138888888888889	0.000278283276775093\\
0.0136986301369863	0.000270739392218111\\
0.0135135135135135	0.000263635480745789\\
0.0133333333333333	0.000256704411237152\\
0.0131578947368421	0.000250125586357147\\
0.012987012987013	0.000243761401545334\\
0.0128205128205128	0.000237676437711581\\
0.0126582278481013	0.000231791143636445\\
0.0125	0.000226075647110913\\
0.0123456790123457	0.000220641322236398\\
0.0121951219512195	0.000215303371498166\\
0.0120481927710843	0.00021024318650531\\
0.0119047619047619	0.000205288007009835\\
0.0117647058823529	0.00020055465196997\\
0.0116279069767442	0.000195971252587523\\
0.0114942528735632	0.000191520501066217\\
0.0113636363636364	0.000187248775367577\\
0.0112359550561798	0.000183076566902446\\
0.0111111111111111	0.000179102215293575\\
0.010989010989011	0.000175181157201498\\
0.0108695652173913	0.000171453499560298\\
0.010752688172043	0.000167792234311098\\
0.0106382978723404	0.000164271714331005\\
0.0105263157894737	0.000160852098785913\\
0.0104166666666667	0.000157530712378978\\
0.0103092783505155	0.000154329098982747\\
0.0102040816326531	0.000151194246740927\\
0.0101010101010101	0.000148197535117384\\
0.01	0.000145234757171347\\
};
\addlegendentry{PCA}

\addplot [color=color2, mark=*, mark options={solid, color2, scale=0.75}]
  table[row sep=crcr]{%
0.1	0.000187601564665629\\
0.0909090909090909	0.000434434956055524\\
0.0833333333333333	0.000171522123783241\\
0.0769230769230769	0.000309000517879271\\
0.0714285714285714	0.000815476662548748\\
0.0666666666666667	0.0012600202677252\\
0.0625	0.00160460131846935\\
0.0588235294117647	0.00184286434903139\\
0.0555555555555556	0.00198745886612528\\
0.0526315789473684	0.00205361243652996\\
0.05	0.00206318386838711\\
0.0476190476190476	0.00203057659760764\\
0.0454545454545455	0.00197199343730278\\
0.0434782608695652	0.00189580996399158\\
0.0416666666666667	0.00181119707181443\\
0.04	0.00172174558012586\\
0.0384615384615385	0.001632088390081\\
0.037037037037037	0.00154357829664431\\
0.0357142857142857	0.0014583082354877\\
0.0344827586206897	0.00137644124253436\\
0.0333333333333333	0.00129884388289506\\
0.032258064516129	0.00122530952542782\\
0.03125	0.00115612241489207\\
0.0303030303030303	0.00109099690034786\\
0.0294117647058824	0.00103005571532033\\
0.0285714285714286	0.000973002195606781\\
0.0277777777777778	0.000919825448781386\\
0.027027027027027	0.000870236038752997\\
0.0263157894736842	0.000824195113625054\\
0.0256410256410256	0.000781381266971359\\
0.025	0.000741718743901748\\
0.024390243902439	0.000704895341372552\\
0.0238095238095238	0.000670827935304513\\
0.0232558139534884	0.000639185896745526\\
0.0227272727272727	0.000609894118244553\\
0.0222222222222222	0.000582649212505587\\
0.0217391304347826	0.000557379296780458\\
0.0212765957446809	0.000533797986800061\\
0.0208333333333333	0.000511859580969248\\
0.0204081632653061	0.000491313799348997\\
0.02	0.000472123856172901\\
0.0196078431372549	0.000454071488977581\\
0.0192307692307692	0.00043715029195307\\
0.0188679245283019	0.00042116980645171\\
0.0185185185185185	0.000406123457725283\\
0.0181818181818182	0.000391866069386709\\
0.0178571428571429	0.000378394944784088\\
0.0175438596491228	0.000365593204472187\\
0.0172413793103448	0.000353455896666494\\
0.0169491525423729	0.000341897387933193\\
0.0166666666666667	0.00033091124575485\\
0.0163934426229508	0.000320428122877997\\
0.0161290322580645	0.00031044142970077\\
0.0158730158730159	0.000300900363393142\\
0.015625	0.000291794538211598\\
0.0153846153846154	0.000283082902297893\\
0.0151515151515152	0.000274756076938026\\
0.0149253731343284	0.00026678315749365\\
0.0147058823529412	0.000259150432086491\\
0.0144927536231884	0.000251834384453725\\
0.0142857142857143	0.000244822687049534\\
0.0140845070422535	0.000238096279652794\\
0.0138888888888889	0.000231640259487104\\
0.0136986301369863	0.000225441673198268\\
0.0135135135135135	0.000219486237393179\\
0.0133333333333333	0.000213762883644414\\
0.0131578947368421	0.000208256632577564\\
0.012987012987013	0.000202961379485139\\
0.0128205128205128	0.00019786177587286\\
0.0126582278481013	0.000192952672297436\\
0.0125	0.000188219641663823\\
0.0123456790123457	0.000183660605859615\\
0.0121951219512195	0.00017926056722864\\
0.0120481927710843	0.000175018279232192\\
0.0119047619047619	0.000170920158345272\\
0.0117647058823529	0.000166966645553823\\
0.0116279069767442	0.000163143562650458\\
0.0114942528735632	0.000159452442828695\\
0.0113636363636364	0.00015588054797458\\
0.0112359550561798	0.000152429707708179\\
0.0111111111111111	0.000149087075730858\\
0.010989010989011	0.000145855806194833\\
0.0108695652173913	0.000142723761351049\\
0.010752688172043	0.00013969408115197\\
0.0106382978723404	0.000136754865390776\\
0.0105263157894737	0.000133910406295901\\
0.0104166666666667	0.000131149032185196\\
0.0103092783505155	0.000128474993644712\\
0.0102040816326531	0.000125877059770463\\
0.0101010101010101	0.000123360317877075\\
0.01	0.000120913579081727\\
};
\addlegendentry{Com}

\addplot [dashed]
  table[row sep=crcr]{%
0.1	1e-1\\
0.01	1e-03\\
};
\end{axis}
\end{tikzpicture}%

%% file: errors-BSamND-PutOnAvg-set3-5D-r=0.05-yDomain-non-uniform-M=10to100.tikz
%
%
%
\begin{tikzpicture}

\begin{axis}[%
width=0.37\textwidth,
scale only axis,
xmode=log,
xmin=0.01,
xmax=0.1,
xminorticks=true,
ymode=log,
ymin=1e-5,
ymax=1e-0,
yminorticks=true,
title={American basket option},
legend style={nodes={scale=0.75, transform shape}},
legend cell align={left}
]
\addplot [color=color1, mark=square*, mark options={solid, color1, scale=0.75}]
  table[row sep=crcr]{%
0.1	0.000342031634473744\\
0.0909090909090909	0.00133145540257762\\
0.0833333333333333	0.0028195679857646\\
0.0769230769230769	0.00433913159083393\\
0.0714285714285714	0.00547888857838522\\
0.0666666666666667	0.00589916729022316\\
0.0625	0.00587639383062513\\
0.0588235294117647	0.00553537515782254\\
0.0555555555555556	0.00501540177194104\\
0.0526315789473684	0.00445721609067762\\
0.05	0.00383942049170041\\
0.0476190476190476	0.00331371219854401\\
0.0454545454545455	0.00282421158976093\\
0.0434782608695652	0.00240719243694812\\
0.0416666666666667	0.00203791224397729\\
0.04	0.00177286349524575\\
0.0384615384615385	0.00156837766263146\\
0.037037037037037	0.00142666004441883\\
0.0357142857142857	0.00135993506540755\\
0.0344827586206897	0.00130304462919012\\
0.0333333333333333	0.00127856120676803\\
0.032258064516129	0.00126682024692701\\
0.03125	0.00126133816564436\\
0.0303030303030303	0.00123615335141938\\
0.0294117647058824	0.0012159481668822\\
0.0285714285714286	0.00114964165746981\\
0.0277777777777778	0.00107395276890235\\
0.027027027027027	0.00100695225194389\\
0.0263157894736842	0.00092162730475398\\
0.0256410256410256	0.000847728743651271\\
0.025	0.000766419383938372\\
0.024390243902439	0.000710411420170753\\
0.0238095238095238	0.000650165315268542\\
0.0232558139534884	0.00060970939945279\\
0.0227272727272727	0.000577288750340965\\
0.0222222222222222	0.000549818729671525\\
0.0217391304347826	0.0005342346327811\\
0.0212765957446809	0.0005186972994147\\
0.0208333333333333	0.000511102940875363\\
0.0204081632653061	0.000503383567086707\\
0.02	0.000499886427438756\\
0.0196078431372549	0.000483126414426605\\
0.0192307692307692	0.000476509495015431\\
0.0188679245283019	0.000467441382693251\\
0.0185185185185185	0.000448174646512117\\
0.0181818181818182	0.000426438505353177\\
0.0178571428571429	0.000403424913009154\\
0.0175438596491228	0.000381068186267902\\
0.0172413793103448	0.00036312823617704\\
0.0169491525423729	0.000347369852193274\\
0.0166666666666667	0.000327161980240631\\
0.0163934426229508	0.000310365415309793\\
0.0161290322580645	0.000300537689153058\\
0.0158730158730159	0.000291479707952902\\
0.015625	0.00028319995983403\\
0.0153846153846154	0.000274964979893849\\
0.0151515151515152	0.000268609859043423\\
0.0149253731343284	0.000263958092669686\\
0.0147058823529412	0.000262645166646563\\
0.0144927536231884	0.000261391791578747\\
0.0142857142857143	0.000252435859950773\\
0.0140845070422535	0.000246940251567879\\
0.0138888888888889	0.000241859454426341\\
0.0136986301369863	0.000234528054193328\\
0.0135135135135135	0.000229548829169229\\
0.0133333333333333	0.000218841814787235\\
0.0131578947368421	0.000210141879368742\\
0.012987012987013	0.000202094934598579\\
0.0128205128205128	0.000196764521120718\\
0.0126582278481013	0.000189338031969022\\
0.0125	0.000183716336211087\\
0.0123456790123457	0.00017813656200627\\
0.0121951219512195	0.000172891479313098\\
0.0120481927710843	0.000168911521352788\\
0.0119047619047619	0.000168351866038086\\
0.0117647058823529	0.00016581977658478\\
0.0116279069767442	0.000160923308323346\\
0.0114942528735632	0.000158539314655287\\
0.0113636363636364	0.00015665891029587\\
0.0112359550561798	0.0001543760597148\\
0.0111111111111111	0.000153058941694012\\
0.010989010989011	0.000148830610319256\\
0.0108695652173913	0.000144117163545537\\
0.010752688172043	0.000141629940918764\\
0.0106382978723404	0.00013880856488574\\
0.0105263157894737	0.00013509969363168\\
0.0104166666666667	0.000129513693009342\\
0.0103092783505155	0.000126864127855369\\
0.0102040816326531	0.000122526815087953\\
0.0101010101010101	0.000120219383940245\\
0.01	0.000118715213449028\\
};
\addlegendentry{PCA (IT)}

\addplot [color=color2, mark=*, mark options={solid, color2, scale=0.75}]
  table[row sep=crcr]{%
0.1	0.000296881718446829\\
0.0909090909090909	0.00176786993673647\\
0.0833333333333333	0.00346466451391994\\
0.0769230769230769	0.00470629302170583\\
0.0714285714285714	0.00520721534698526\\
0.0666666666666667	0.00521377007935317\\
0.0625	0.00491999587164016\\
0.0588235294117647	0.00444680951080345\\
0.0555555555555556	0.003910189328654\\
0.0526315789473684	0.00333040578457924\\
0.05	0.00279727785274692\\
0.0476190476190476	0.0023307396013009\\
0.0454545454545455	0.00195897841971832\\
0.0434782608695652	0.00165676501923989\\
0.0416666666666667	0.00143883187273516\\
0.04	0.00128364801480552\\
0.0384615384615385	0.00119280324965912\\
0.037037037037037	0.00115236002278657\\
0.0357142857142857	0.00114420921909744\\
0.0344827586206897	0.00113865489915899\\
0.0333333333333333	0.00114180369358963\\
0.032258064516129	0.0011246184139396\\
0.03125	0.00109124329302837\\
0.0303030303030303	0.00103003159578058\\
0.0294117647058824	0.000948546703660069\\
0.0285714285714286	0.000869824765615868\\
0.0277777777777778	0.000787728536698867\\
0.027027027027027	0.000708449822624135\\
0.0263157894736842	0.00064361033317864\\
0.0256410256410256	0.000581962046993928\\
0.025	0.000534134029113287\\
0.024390243902439	0.000504998416222147\\
0.0238095238095238	0.000480848757675734\\
0.0232558139534884	0.00046529211246088\\
0.0227272727272727	0.000456848565940443\\
0.0222222222222222	0.000453959759947165\\
0.0217391304347826	0.000450995432988183\\
0.0212765957446809	0.000451123129861924\\
0.0208333333333333	0.000441383462750772\\
0.0204081632653061	0.000429467446744464\\
0.02	0.000415869514006878\\
0.0196078431372549	0.000394416600096786\\
0.0192307692307692	0.000371829182169825\\
0.0188679245283019	0.000346247630000934\\
0.0185185185185185	0.000323085092252623\\
0.0181818181818182	0.000303969636236384\\
0.0178571428571429	0.000287312940500079\\
0.0175438596491228	0.000270454965146771\\
0.0172413793103448	0.000257580248633232\\
0.0169491525423729	0.000249616542324194\\
0.0166666666666667	0.000244536367967202\\
0.0163934426229508	0.000240477596239852\\
0.0161290322580645	0.000236736318526387\\
0.0158730158730159	0.000234475050885558\\
0.015625	0.000233432754337348\\
0.0153846153846154	0.000231248636665432\\
0.0151515151515152	0.000226360178664342\\
0.0149253731343284	0.00022065356503076\\
0.0147058823529412	0.000212915361445049\\
0.0144927536231884	0.00020553818073149\\
0.0142857142857143	0.000196767177080032\\
0.0140845070422535	0.000188062925378985\\
0.0138888888888889	0.00017852694858711\\
0.0136986301369863	0.000170653790107916\\
0.0135135135135135	0.000163713757744849\\
0.0133333333333333	0.000158248306674352\\
0.0131578947368421	0.000153729368675332\\
0.012987012987013	0.000149419024998598\\
0.0128205128205128	0.000146799357391619\\
0.0126582278481013	0.00014476108286704\\
0.0125	0.000144105474184919\\
0.0123456790123457	0.000141849595228227\\
0.0121951219512195	0.000140429004332393\\
0.0120481927710843	0.000138025662165836\\
0.0119047619047619	0.000136410070264104\\
0.0117647058823529	0.000133960099818087\\
0.0116279069767442	0.000130039892282996\\
0.0114942528735632	0.000125814393851148\\
0.0113636363636364	0.000121281965834219\\
0.0112359550561798	0.000117925432134791\\
0.0111111111111111	0.000113690082801654\\
0.010989010989011	0.000109978829738039\\
0.0108695652173913	0.00010635643605364\\
0.010752688172043	0.000103424740586\\
0.0106382978723404	0.000101587851759205\\
0.0105263157894737	9.94066623268841e-05\\
0.0104166666666667	9.80099716214045e-05\\
0.0103092783505155	9.6595849316844e-05\\
0.0102040816326531	9.54340441687018e-05\\
0.0101010101010101	9.44981848332582e-05\\
0.01	9.37579345651751e-05\\
};
\addlegendentry{Com (IT)}

\addplot [dashed]
  table[row sep=crcr]{%
0.1	1e-1\\
0.01	1e-03\\
};
\end{axis}
\end{tikzpicture}%

%% file: errors-BSeuND-PutOnAvg-set2-10D-r=0.06-yDomain-non-uniform-M=10to100.tikz
%
%
%
\begin{tikzpicture}

\begin{axis}[%
width=0.37\textwidth,
scale only axis,
xmode=log,
xmin=0.01,
xmax=0.1,
xminorticks=true,
ymode=log,
ymin=1e-5,
ymax=1e-0,
yminorticks=true,
ylabel={error},
legend cell align={left},
legend style={nodes={scale=0.75, transform shape}}
]
\addplot [color=color1, mark=square*, mark options={solid, color1, scale=0.75}]
  table[row sep=crcr]{%
0.1	0.0203018278276195\\
0.0909090909090909	0.0279359719045927\\
0.0833333333333333	0.0179158601182613\\
0.0769230769230769	0.0205485454436769\\
0.0714285714285714	0.0144427088096181\\
0.0666666666666667	0.0156326346127843\\
0.0625	0.0118424471438877\\
0.0588235294117647	0.0123714738038169\\
0.0555555555555556	0.0097591079905357\\
0.0526315789473684	0.00992296534208115\\
0.05	0.00821751934831205\\
0.0476190476190476	0.00825508123812746\\
0.0454545454545455	0.00690693510689955\\
0.0434782608695652	0.00687123407978163\\
0.0416666666666667	0.00596250648857799\\
0.04	0.00588846486129946\\
0.0384615384615385	0.00513625250566407\\
0.037037037037037	0.00506938567227266\\
0.0357142857142857	0.00452180131878732\\
0.0344827586206897	0.00442228150237156\\
0.0333333333333333	0.00396890357504931\\
0.032258064516129	0.0038867048866571\\
0.03125	0.00353562290524267\\
0.0303030303030303	0.00344635214679867\\
0.0294117647058824	0.00315192574490286\\
0.0285714285714286	0.0030813047443059\\
0.0277777777777778	0.00283926123877765\\
0.027027027027027	0.00276744664902151\\
0.0263157894736842	0.00256884679751124\\
0.0256410256410256	0.00250596835454053\\
0.025	0.00233295610288398\\
0.024390243902439	0.00227163570914779\\
0.0238095238095238	0.00213163557390428\\
0.0232558139534884	0.00208017835484176\\
0.0227272727272727	0.00195544790262847\\
0.0222222222222222	0.00190594968834967\\
0.0217391304347826	0.00180041214759719\\
0.0212765957446809	0.00175551525486262\\
0.0208333333333333	0.00166328305250596\\
0.0204081632653061	0.00162242937874058\\
0.02	0.00154142265418578\\
0.0196078431372549	0.00150455731318433\\
0.0192307692307692	0.00143347579894759\\
0.0188679245283019	0.00140083717091299\\
0.0185185185185185	0.00133786339152153\\
0.0181818181818182	0.00130652171042789\\
0.0178571428571429	0.00124997900499946\\
0.0175438596491228	0.00122260066564384\\
0.0172413793103448	0.00117260385516327\\
0.0169491525423729	0.00114715565370227\\
0.0166666666666667	0.00110117541965071\\
0.0163934426229508	0.00107829031976625\\
0.0161290322580645	0.00103807215192842\\
0.0158730158730159	0.00101627880784627\\
0.015625	0.000979836606479267\\
0.0153846153846154	0.000959052287547446\\
0.0151515151515152	0.000926584103077333\\
0.0149253731343284	0.000907928868475683\\
0.0147058823529412	0.000878139562840885\\
0.0144927536231884	0.000860434012619149\\
0.0142857142857143	0.000834057280636347\\
0.0140845070422535	0.000817365605930931\\
0.0138888888888889	0.000792378401266669\\
0.0136986301369863	0.000777424069941302\\
0.0135135135135135	0.000755178163716175\\
0.0133333333333333	0.000740664723493101\\
0.0131578947368421	0.000719646906195415\\
0.012987012987013	0.000707270012235361\\
0.0128205128205128	0.000687918562114698\\
0.0126582278481013	0.000675551283882192\\
0.0125	0.000657547979343032\\
0.0123456790123457	0.000646783752112245\\
0.0121951219512195	0.000630193474940755\\
0.0120481927710843	0.00061942189746389\\
0.0119047619047619	0.000604165058534778\\
0.0117647058823529	0.000594384297672002\\
0.0116279069767442	0.000580259177379472\\
0.0114942528735632	0.000570802971261664\\
0.0113636363636364	0.000557657148762325\\
0.0112359550561798	0.000549050297313891\\
0.0111111111111111	0.000536669065587581\\
0.010989010989011	0.000528288608747207\\
0.0108695652173913	0.000517286869486044\\
0.010752688172043	0.000509471213516299\\
0.0106382978723404	0.000498670360084663\\
0.0105263157894737	0.000491262686514116\\
0.0104166666666667	0.00048142037604948\\
0.0103092783505155	0.000474394019646995\\
0.0102040816326531	0.000464969266554149\\
0.0101010101010101	0.000458539731005514\\
0.01	0.000449779005556872\\
};
\addlegendentry{PCA}

\addplot [color=color2, mark=*, mark options={solid, color2, scale=0.75}]
  table[row sep=crcr]{%
0.1	0.0313103692839237\\
0.0909090909090909	0.0265708900113523\\
0.0833333333333333	0.022975568677526\\
0.0769230769230769	0.0198152743278951\\
0.0714285714285714	0.0173714050139888\\
0.0666666666666667	0.0152078965735831\\
0.0625	0.0135297670319666\\
0.0588235294117647	0.0120228268485102\\
0.0555555555555556	0.0108453467739891\\
0.0526315789473684	0.00977330624179795\\
0.05	0.00890942716271415\\
0.0476190476190476	0.00810111757930188\\
0.0454545454545455	0.00743754362562232\\
0.0434782608695652	0.0068153635424566\\
0.0416666666666667	0.00629952453857652\\
0.04	0.00581325521301079\\
0.0384615384615385	0.00540530917880455\\
0.037037037037037	0.00501739811009649\\
0.0357142857142857	0.00468799673471609\\
0.0344827586206897	0.00437388589489884\\
0.0333333333333333	0.00410447615332588\\
0.032258064516129	0.00384624356233532\\
0.03125	0.00362303242830153\\
0.0303030303030303	0.00340854185345796\\
0.0294117647058824	0.00322145130555707\\
0.0285714285714286	0.00304118335895442\\
0.0277777777777778	0.00288304953756779\\
0.027027027027027	0.00273014837055019\\
0.0263157894736842	0.00259508216447801\\
0.0256410256410256	0.00246436106283565\\
0.025	0.00234821258151419\\
0.024390243902439	0.00223546033472255\\
0.0238095238095238	0.00213481881838118\\
0.0232558139534884	0.00203701681808488\\
0.0227272727272727	0.00194920629682305\\
0.0222222222222222	0.00186373865949552\\
0.0217391304347826	0.00178674876881357\\
0.0212765957446809	0.00171165744111079\\
0.0208333333333333	0.00164369733158476\\
0.0204081632653061	0.00157739116625932\\
0.02	0.0015171628128684\\
0.0196078431372549	0.00145827541924692\\
0.0192307692307692	0.00140462296910437\\
0.0188679245283019	0.001352144550132\\
0.0185185185185185	0.00130413672293783\\
0.0181818181818182	0.00125712573167014\\
0.0178571428571429	0.00121403021825361\\
0.0175438596491228	0.00117177392791556\\
0.0172413793103448	0.00113290466733262\\
0.0169491525423729	0.00109478621986592\\
0.0166666666666667	0.00105963979420087\\
0.0163934426229508	0.00102511721822918\\
0.0161290322580645	0.000993214802668785\\
0.0158730158730159	0.000961877552022927\\
0.015625	0.000932832869779854\\
0.0153846153846154	0.00090427584088526\\
0.0151515151515152	0.000877771201209021\\
0.0149253731343284	0.000851689276265111\\
0.0147058823529412	0.000827417875146064\\
0.0144927536231884	0.000803530846435474\\
0.0142857142857143	0.000781267033932442\\
0.0140845070422535	0.000759327693589107\\
0.0138888888888889	0.000738843286423374\\
0.0136986301369863	0.000718660057917542\\
0.0135135135135135	0.000699773620255195\\
0.0133333333333333	0.000681149485035171\\
0.0131578947368421	0.00066370443292163\\
0.012987012987013	0.000646492770299156\\
0.0128205128205128	0.000630336188819247\\
0.0126582278481013	0.000614394038178867\\
0.0125	0.000599413231083479\\
0.0123456790123457	0.000584616048456166\\
0.0121951219512195	0.000570690825530074\\
0.0120481927710843	0.00055693916605537\\
0.0119047619047619	0.00054397630416525\\
0.0117647058823529	0.000531165240804765\\
0.0116279069767442	0.000519079624128449\\
0.0114942528735632	0.000507132299793933\\
0.0113636363636364	0.000495841447241907\\
0.0112359550561798	0.00048467816715736\\
0.0111111111111111	0.000474120795733723\\
0.010989010989011	0.000463674012521453\\
0.0108695652173913	0.000453781785327845\\
0.010752688172043	0.000443995407353848\\
0.0106382978723404	0.000434716971355775\\
0.0105263157894737	0.000425531203549712\\
0.0104166666666667	0.000416816632651251\\
0.0103092783505155	0.000408188281268429\\
0.0102040816326531	0.000399990284583374\\
0.0101010101010101	0.000391871825075718\\
0.01	0.000384154659481939\\
};
\addlegendentry{Com}

\addplot [dashed]
  table[row sep=crcr]{%
0.1	1e-1\\
0.01	1e-03\\
};

\end{axis}
\end{tikzpicture}%

%% file: errors-BSamND-PutOnAvg-set2-10D-r=0.06-yDomain-non-uniform-M=10to100.tikz
%
%
%
\begin{tikzpicture}

\begin{axis}[%
width=0.37\textwidth,
scale only axis,
xmode=log,
xmin=0.01,
xmax=0.1,
xminorticks=true,
ymode=log,
ymin=1e-5,
ymax=1e-0,
yminorticks=true,
legend cell align={left},
legend style={nodes={scale=0.75, transform shape}}
]
\addplot [color=color1, mark=square*, mark options={solid, color1, scale=0.75}]
  table[row sep=crcr]{%
0.1	0.0299790834363509\\
0.0909090909090909	0.0474340193986511\\
0.0833333333333333	0.0213952689593049\\
0.0769230769230769	0.0200763777915358\\
0.0714285714285714	0.0232608049816625\\
0.0666666666666667	0.01703267525235\\
0.0625	0.0117361409110051\\
0.0588235294117647	0.0141567222834005\\
0.0555555555555556	0.0186341442841811\\
0.0526315789473684	0.00925632652894115\\
0.05	0.00857132193576104\\
0.0476190476190476	0.0123953833862991\\
0.0454545454545455	0.00897357344304672\\
0.0434782608695652	0.00617841004372921\\
0.0416666666666667	0.00749060892904985\\
0.04	0.0106921306214338\\
0.0384615384615385	0.00411083510740817\\
0.037037037037037	0.00516344732539742\\
0.0357142857142857	0.00740847840642056\\
0.0344827586206897	0.00569155297859836\\
0.0333333333333333	0.00347741437788818\\
0.032258064516129	0.0049334100800742\\
0.03125	0.00647986556530555\\
0.0303030303030303	0.00189538354444574\\
0.0294117647058824	0.00327177002805246\\
0.0285714285714286	0.00519936624403528\\
0.0277777777777778	0.00381527824800809\\
0.027027027027027	0.00210856592519981\\
0.0263157894736842	0.00344962293488948\\
0.0256410256410256	0.00446541304595893\\
0.025	0.000994838010966781\\
0.024390243902439	0.00228308121684506\\
0.0238095238095238	0.00375449378535175\\
0.0232558139534884	0.00285470165284174\\
0.0227272727272727	0.00123319680259826\\
0.0222222222222222	0.00264661783230458\\
0.0217391304347826	0.00321910416972493\\
0.0212765957446809	0.000926566461711076\\
0.0208333333333333	0.00169017959043538\\
0.0204081632653061	0.00284966662036323\\
0.02	0.00215869527778967\\
0.0196078431372549	0.000756241487118059\\
0.0192307692307692	0.00202035667808409\\
0.0188679245283019	0.00243760162907991\\
0.0185185185185185	0.000973750264173567\\
0.0181818181818182	0.0012547301952186\\
0.0178571428571429	0.00218968253617935\\
0.0175438596491228	0.00174285757347659\\
0.0172413793103448	0.00036664566181388\\
0.0169491525423729	0.00158469639137526\\
0.0166666666666667	0.00192180166150591\\
0.0163934426229508	0.00100352445788632\\
0.0161290322580645	0.000905608949598236\\
0.0158730158730159	0.00175938482061455\\
0.015625	0.00141025804067563\\
0.0153846153846154	0.000207896734726765\\
0.0151515151515152	0.00128706630340458\\
0.0149253731343284	0.00150504530228468\\
0.0147058823529412	0.000952470872416233\\
0.0144927536231884	0.000634671895261452\\
0.0142857142857143	0.00139110901596151\\
0.0140845070422535	0.0012068301761774\\
0.0138888888888889	0.00022376714651573\\
0.0136986301369863	0.00102510743888318\\
0.0135135135135135	0.00121756569752063\\
0.0133333333333333	0.000878812041835486\\
0.0131578947368421	0.000436789034816698\\
0.012987012987013	0.00107989081982063\\
0.0128205128205128	0.00103551644800093\\
0.0126582278481013	0.000350414362646267\\
0.0125	0.000789292943073283\\
0.0123456790123457	0.00101659062861259\\
0.0121951219512195	0.000791692860285576\\
0.0120481927710843	0.000280588472532539\\
0.0119047619047619	0.000878179534331336\\
0.0117647058823529	0.000891410658339531\\
0.0116279069767442	0.000442219097396945\\
0.0114942528735632	0.000591623210032477\\
0.0113636363636364	0.00085537939833813\\
0.0112359550561798	0.00074650342472582\\
0.0111111111111111	0.000164805419694769\\
0.010989010989011	0.000709930013976123\\
0.0108695652173913	0.000762425278332701\\
0.010752688172043	0.000501763256181009\\
0.0106382978723404	0.000428616504988044\\
0.0105263157894737	0.000693197497880771\\
0.0104166666666667	0.000694140292950784\\
0.0103092783505155	0.000179887166726012\\
0.0102040816326531	0.000551799043107648\\
0.0101010101010101	0.000680781688509491\\
0.01	0.000508985652112992\\
};
\addlegendentry{PCA (IT)}

\addplot [color=color2, mark=*, mark options={solid, color2, scale=0.75}]
  table[row sep=crcr]{%
0.1	0.049520146543151\\
0.0909090909090909	0.0337740862779861\\
0.0833333333333333	0.0313204172655195\\
0.0769230769230769	0.0305117487568072\\
0.0714285714285714	0.023008598435464\\
0.0666666666666667	0.0206944013875143\\
0.0625	0.0202539031166116\\
0.0588235294117647	0.0166286987663531\\
0.0555555555555556	0.0147271006442593\\
0.0526315789473684	0.0143157718168401\\
0.05	0.0129027407323379\\
0.0476190476190476	0.0110537944304525\\
0.0454545454545455	0.0106504238237675\\
0.0434782608695652	0.0102908841642075\\
0.0416666666666667	0.00869330188192663\\
0.04	0.00826621552706275\\
0.0384615384615385	0.00821810485857877\\
0.037037037037037	0.00709091464181499\\
0.0357142857142857	0.00661065399921612\\
0.0344827586206897	0.00662474901061572\\
0.0333333333333333	0.00593148986453018\\
0.032258064516129	0.00542473732905835\\
0.03125	0.00546458679667716\\
0.0303030303030303	0.00509165479426898\\
0.0294117647058824	0.00455254332759103\\
0.0285714285714286	0.00455266852840208\\
0.0277777777777778	0.004459062969703\\
0.027027027027027	0.00392373478850838\\
0.0263157894736842	0.00385228406028171\\
0.0256410256410256	0.0038423431189516\\
0.025	0.00343038653727268\\
0.024390243902439	0.00331019368298091\\
0.0238095238095238	0.00333139005331629\\
0.0232558139534884	0.00303867114557166\\
0.0227272727272727	0.00288649309225919\\
0.0222222222222222	0.00290657235128466\\
0.0217391304347826	0.00272907256328803\\
0.0212765957446809	0.00253493582179365\\
0.0208333333333333	0.00255251053207828\\
0.0204081632653061	0.00248466155260929\\
0.02	0.0022638492461291\\
0.0196078431372549	0.00224862841767792\\
0.0192307692307692	0.00225491257033106\\
0.0188679245283019	0.00204411864093212\\
0.0185185185185185	0.00200061531623286\\
0.0181818181818182	0.00203094418683425\\
0.0178571428571429	0.00185988968771267\\
0.0175438596491228	0.00180241394435243\\
0.0172413793103448	0.00182648415343944\\
0.0169491525423729	0.00170521676762325\\
0.0166666666666667	0.00163734694449147\\
0.0163934426229508	0.0016424709860221\\
0.0161290322580645	0.0015896989705233\\
0.0158730158730159	0.00148540694740484\\
0.015625	0.00148716568788254\\
0.0153846153846154	0.00149530487351779\\
0.0151515151515152	0.00136096191040047\\
0.0149253731343284	0.00135870260487936\\
0.0147058823529412	0.00137084708172619\\
0.0144927536231884	0.00126510111789568\\
0.0142857142857143	0.0012460333355333\\
0.0140845070422535	0.0012547682998123\\
0.0138888888888889	0.00118722990156717\\
0.0136986301369863	0.00113956938361559\\
0.0135135135135135	0.00115736289572332\\
0.0133333333333333	0.00111338871437949\\
0.0131578947368421	0.0010555872786675\\
0.012987012987013	0.00107040031743422\\
0.0128205128205128	0.00104875894804546\\
0.0126582278481013	0.000985008641175122\\
0.0125	0.000987192768564515\\
0.0123456790123457	0.000989828712861796\\
0.0121951219512195	0.00092465453500945\\
0.0120481927710843	0.000910360417991418\\
0.0119047619047619	0.000927194276114651\\
0.0117647058823529	0.000870619492031599\\
0.0116279069767442	0.000847836804011459\\
0.0114942528735632	0.000862812684511605\\
0.0113636363636364	0.000822751652650666\\
0.0112359550561798	0.000795105211061076\\
0.0111111111111111	0.000804083522091936\\
0.010989010989011	0.000784536918104317\\
0.0108695652173913	0.000748522712528565\\
0.010752688172043	0.000749719151177386\\
0.0106382978723404	0.000751707215799868\\
0.0105263157894737	0.000705253284062124\\
0.0104166666666667	0.000702071694582918\\
0.0103092783505155	0.000712227495414108\\
0.0102040816326531	0.000667699454122594\\
0.0101010101010101	0.000658671542822065\\
0.01	0.000670841739900174\\
};
\addlegendentry{Com (IT)}

\addplot [dashed]
  table[row sep=crcr]{%
0.1	1e-1\\
0.01	1e-03\\
};

\end{axis}
\end{tikzpicture}%

%% file: errors-BSeuND-PutOnAvg-set2-15D-r=0.06-yDomain-non-uniform-M=10to100.tikz
%
%
%
\begin{tikzpicture}

\begin{axis}[%
width=0.37\textwidth,
scale only axis,
xmode=log,
xmin=0.01,
xmax=0.1,
xminorticks=true,
xlabel={1/m},
ymode=log,
ymin=1e-5,
ymax=1e-0,
yminorticks=true,
ylabel={error},
legend cell align={left},
legend style={nodes={scale=0.75, transform shape}}
]
\addplot [color=color1, mark=square*, mark options={solid, color1, scale=0.75}]
  table[row sep=crcr]{%
0.1	0.011372671668017\\
0.0909090909090909	0.0218171060453025\\
0.0833333333333333	0.0116591976555542\\
0.0769230769230769	0.0156948312741144\\
0.0714285714285714	0.0100420994310745\\
0.0666666666666667	0.0122650019611489\\
0.0625	0.00851232107141098\\
0.0588235294117647	0.00951791487189291\\
0.0555555555555556	0.00703405254532008\\
0.0526315789473684	0.00767947480220821\\
0.05	0.00599643610018175\\
0.0476190476190476	0.00628495490387637\\
0.0454545454545455	0.00508888910355743\\
0.0434782608695652	0.00528436996994097\\
0.0416666666666667	0.00440320903883107\\
0.04	0.00445762273414696\\
0.0384615384615385	0.00379427537210675\\
0.037037037037037	0.003844794296331\\
0.0357142857142857	0.00332236318114054\\
0.0344827586206897	0.00332730563288974\\
0.0333333333333333	0.00292369558531647\\
0.032258064516129	0.00292466418075399\\
0.03125	0.00258808693198476\\
0.0303030303030303	0.00256826100770458\\
0.0294117647058824	0.00230962195280937\\
0.0285714285714286	0.00228386066249053\\
0.0277777777777778	0.00206718824940622\\
0.027027027027027	0.00204047563766496\\
0.0263157894736842	0.00186153383281196\\
0.0256410256410256	0.00183009400183531\\
0.025	0.00168175512788471\\
0.024390243902439	0.0016566331006922\\
0.0238095238095238	0.00152367848038326\\
0.0232558139534884	0.0014984589986673\\
0.0227272727272727	0.00139177361884246\\
0.0222222222222222	0.00136722438870995\\
0.0217391304347826	0.00126774345621294\\
0.0212765957446809	0.00124768031068279\\
0.0208333333333333	0.00116722312824924\\
0.0204081632653061	0.00114312015628626\\
0.02	0.0010715312037276\\
0.0196078431372549	0.00105259784073586\\
0.0192307692307692	0.000989798419146593\\
0.0188679245283019	0.000967992261469819\\
0.0185185185185185	0.000913613312882977\\
0.0181818181818182	0.00089721649106711\\
0.0178571428571429	0.000847168777151297\\
0.0175438596491228	0.000828842431122134\\
0.0172413793103448	0.000786233460585595\\
0.0169491525423729	0.000772224451406012\\
0.0166666666666667	0.000732017248712191\\
0.0163934426229508	0.000715907681522965\\
0.0161290322580645	0.000682290778948769\\
0.0158730158730159	0.000669188486350225\\
0.015625	0.000637062886534201\\
0.0153846153846154	0.000623381976963655\\
0.0151515151515152	0.000596059854346964\\
0.0149253731343284	0.000583846416915335\\
0.0147058823529412	0.000557825169653414\\
0.0144927536231884	0.000546796478379519\\
0.0142857142857143	0.000524345884974187\\
0.0140845070422535	0.000512647375841224\\
0.0138888888888889	0.000491338063686841\\
0.0136986301369863	0.0004820254213459\\
0.0135135135135135	0.000462850997806186\\
0.0133333333333333	0.000452582663317158\\
0.0131578947368421	0.000434957607658326\\
0.012987012987013	0.000426548987512931\\
0.0128205128205128	0.000410219069831963\\
0.0126582278481013	0.000401602661699085\\
0.0125	0.000386808395214544\\
0.0123456790123457	0.000379066465317823\\
0.0121951219512195	0.000365211535600762\\
0.0120481927710843	0.000357995211917927\\
0.0119047619047619	0.000344963038002355\\
0.0117647058823529	0.000338003877488835\\
0.0116279069767442	0.000326421999753967\\
0.0114942528735632	0.000319932752480101\\
0.0113636363636364	0.000308689221334868\\
0.0112359550561798	0.000302506539232694\\
0.0111111111111111	0.000292489560030362\\
0.010989010989011	0.000286650390434851\\
0.0108695652173913	0.000277136573356263\\
0.010752688172043	0.000271668038045347\\
0.0106382978723404	0.000262693547658399\\
0.0105263157894737	0.00025748024888872\\
0.0104166666666667	0.000249147541581451\\
0.0103092783505155	0.000244434761354384\\
0.0102040816326531	0.0002365839112487\\
0.0101010101010101	0.000231621198057308\\
0.01	0.000224436308295606\\
};
\addlegendentry{PCA}

\addplot [color=color2, mark=*, mark options={solid, color2, scale=0.75}]
  table[row sep=crcr]{%
0.1	0.0250105583483559\\
0.0909090909090909	0.0207336148361136\\
0.0833333333333333	0.0177287421430543\\
0.0769230769230769	0.015288042606274\\
0.0714285714285714	0.0134253235034151\\
0.0666666666666667	0.0118016490662118\\
0.0625	0.0105333554866942\\
0.0588235294117647	0.00938574230529987\\
0.0555555555555556	0.00845938650429878\\
0.0526315789473684	0.00761166872683838\\
0.05	0.00692488434749861\\
0.0476190476190476	0.00629816577858255\\
0.0454545454545455	0.00578458693676132\\
0.0434782608695652	0.00530669592471422\\
0.0416666666666667	0.00490586267516502\\
0.04	0.00452930187562173\\
0.0384615384615385	0.00421078053275592\\
0.037037037037037	0.00391044231284576\\
0.0357142857142857	0.00365351743973552\\
0.0344827586206897	0.00341005304789233\\
0.0333333333333333	0.00319974823052027\\
0.032258064516129	0.00299951068880833\\
0.03125	0.00282512776196342\\
0.0303030303030303	0.00265868224935339\\
0.0294117647058824	0.0025126523487099\\
0.0285714285714286	0.00237270492776487\\
0.0277777777777778	0.00224908721602912\\
0.027027027027027	0.00213036917195608\\
0.0263157894736842	0.00202484393954039\\
0.0256410256410256	0.00192321426164277\\
0.025	0.00183241460379546\\
0.024390243902439	0.00174481148453709\\
0.0238095238095238	0.00166611946360706\\
0.0232558139534884	0.00159003447595074\\
0.0227272727272727	0.00152140613382123\\
0.0222222222222222	0.0014549433132921\\
0.0217391304347826	0.00139471409748604\\
0.0212765957446809	0.00133629419262071\\
0.0208333333333333	0.00128316972871001\\
0.0204081632653061	0.00123156296468963\\
0.02	0.0011844478932248\\
0.0196078431372549	0.00113862728024206\\
0.0192307692307692	0.00109666884527893\\
0.0188679245283019	0.00105580380716386\\
0.0185185185185185	0.00101825915106812\\
0.0181818181818182	0.000981663460055238\\
0.0178571428571429	0.000947949061870901\\
0.0175438596491228	0.000915043004164229\\
0.0172413793103448	0.000884644613455277\\
0.0169491525423729	0.000854956184887512\\
0.0166666666666667	0.000827460687221948\\
0.0163934426229508	0.000800575620511279\\
0.0161290322580645	0.000775620199330862\\
0.0158730158730159	0.000751205060001703\\
0.015625	0.000728488726591592\\
0.0153846153846154	0.00070624222301463\\
0.0151515151515152	0.000685505075427195\\
0.0149253731343284	0.000665184649042083\\
0.0147058823529412	0.000646201779481048\\
0.0144927536231884	0.000627585878225312\\
0.0142857142857143	0.000610167657971172\\
0.0140845070422535	0.000593075312906755\\
0.0138888888888889	0.0005770508530476\\
0.0136986301369863	0.000561317413101947\\
0.0135135135135135	0.000546545768899009\\
0.0133333333333333	0.000532032779989366\\
0.0131578947368421	0.000518382944606932\\
0.012987012987013	0.000504966683066876\\
0.0128205128205128	0.000492331401901858\\
0.0126582278481013	0.000479903950813432\\
0.0125	0.000468182017490881\\
0.0123456790123457	0.000456649583343349\\
0.0121951219512195	0.000445757679238601\\
0.0120481927710843	0.00043503479578233\\
0.0119047619047619	0.000424894421112154\\
0.0117647058823529	0.000414908959466054\\
0.0116279069767442	0.000405453925902965\\
0.0114942528735632	0.000396137751446646\\
0.0113636363636364	0.000387306954201305\\
0.0112359550561798	0.000378603590470039\\
0.0111111111111111	0.000370343461031819\\
0.010989010989011	0.000362198461227425\\
0.0108695652173913	0.000354461188051491\\
0.010752688172043	0.00034682943123987\\
0.0106382978723404	0.000339571125353544\\
0.0105263157894737	0.000332409065639805\\
0.0104166666666667	0.000325591955242288\\
0.0103092783505155	0.000318862838567746\\
0.0102040816326531	0.000312450762127159\\
0.0101010101010101	0.000306119746916078\\
0.01	0.000300082422998948\\
};
\addlegendentry{Com}

\addplot [dashed]
  table[row sep=crcr]{%
0.1	1e-1\\
0.01	1e-03\\
};

\end{axis}
\end{tikzpicture}%

%% file: errors-BSamND-PutOnAvg-set2-15D-r=0.06-yDomain-non-uniform-M=10to100.tikz
%
%
%
\begin{tikzpicture}

\begin{axis}[%
width=0.37\textwidth,
scale only axis,
xmode=log,
xmin=0.01,
xmax=0.1,
xminorticks=true,
xlabel={1/m},
ymode=log,
ymin=1e-5,
ymax=1e-0,
yminorticks=true,
legend cell align={left},
legend style={nodes={scale=0.75, transform shape}}
]
\addplot [color=color1, mark=square*, mark options={solid, color1, scale=0.75}]
  table[row sep=crcr]{%
0.1	0.0226864224967807\\
0.0909090909090909	0.0334712593014292\\
0.0833333333333333	0.028365564895537\\
0.0769230769230769	0.0179476676260832\\
0.0714285714285714	0.014602302126417\\
0.0666666666666667	0.0223117800626609\\
0.0625	0.0121630014350225\\
0.0588235294117647	0.010494877218246\\
0.0555555555555556	0.0111564314656944\\
0.0526315789473684	0.0160603527053678\\
0.05	0.00710793384884467\\
0.0476190476190476	0.00723113092890015\\
0.0454545454545455	0.00911598975254235\\
0.0434782608695652	0.0105039081039724\\
0.0416666666666667	0.00482244566086476\\
0.04	0.00561362098889617\\
0.0384615384615385	0.00812223515230581\\
0.037037037037037	0.00626377199884409\\
0.0357142857142857	0.00362454975143534\\
0.0344827586206897	0.00497025969242104\\
0.0333333333333333	0.00664506853922342\\
0.032258064516129	0.00356351696068224\\
0.03125	0.00302317587249123\\
0.0303030303030303	0.0046559443843921\\
0.0294117647058824	0.00480810776417728\\
0.0285714285714286	0.00127088061304192\\
0.0277777777777778	0.00282043056172898\\
0.027027027027027	0.00429360816782265\\
0.0263157894736842	0.00336979565519857\\
0.0256410256410256	0.0014935859561056\\
0.025	0.0027761648348652\\
0.024390243902439	0.0035806577108346\\
0.0238095238095238	0.00237326891807088\\
0.0232558139534884	0.00155537000094252\\
0.0227272727272727	0.00274445757773933\\
0.0222222222222222	0.00281943394614603\\
0.0217391304347826	0.00156542112721181\\
0.0212765957446809	0.00167140834713941\\
0.0208333333333333	0.00251252546910075\\
0.0204081632653061	0.00216929262718346\\
0.02	0.000808600550107119\\
0.0196078431372549	0.00167376475505843\\
0.0192307692307692	0.00223473002041796\\
0.0188679245283019	0.00169871417180611\\
0.0185185185185185	0.000441959789817314\\
0.0181818181818182	0.001784004487839\\
0.0178571428571429	0.00179916307622552\\
0.0175438596491228	0.00138532388682977\\
0.0172413793103448	0.000799420433519327\\
0.0169491525423729	0.00159496004349902\\
0.0166666666666667	0.00156021438364928\\
0.0163934426229508	0.00099735015015523\\
0.0161290322580645	0.00102193369552972\\
0.0158730158730159	0.00147976243767123\\
0.015625	0.00122238754855042\\
0.0153846153846154	0.000739467110657976\\
0.0151515151515152	0.00109690883544489\\
0.0149253731343284	0.00124866846613925\\
0.0147058823529412	0.00109011696498929\\
0.0144927536231884	0.000352394243565479\\
0.0142857142857143	0.00105870648581652\\
0.0140845070422535	0.00110793989964764\\
0.0138888888888889	0.000930957723959125\\
0.0136986301369863	0.000460359456357029\\
0.0135135135135135	0.000960496290477675\\
0.0133333333333333	0.000958487479123526\\
0.0131578947368421	0.00084542121628961\\
0.012987012987013	0.000570816733547641\\
0.0128205128205128	0.000904132631789789\\
0.0126582278481013	0.000888264198481092\\
0.0125	0.000654162330487651\\
0.0123456790123457	0.000648358580332076\\
0.0121951219512195	0.000830850154732277\\
0.0120481927710843	0.000807416543200157\\
0.0119047619047619	0.000402515356853428\\
0.0117647058823529	0.000646846516445754\\
0.0116279069767442	0.000751588427500893\\
0.0114942528735632	0.000769348965321548\\
0.0113636363636364	0.000204209142302014\\
0.0112359550561798	0.000624934568850044\\
0.0111111111111111	0.000720006315943156\\
0.010989010989011	0.000679493957969513\\
0.0108695652173913	0.000297875548664006\\
0.010752688172043	0.000583124470151519\\
0.0106382978723404	0.000688457108009022\\
0.0105263157894737	0.000586289550818186\\
0.0104166666666667	0.00033504711187482\\
0.0103092783505155	0.000602109782583238\\
0.0102040816326531	0.000671214259405417\\
0.0101010101010101	0.000360808408418567\\
0.01	0.000390268109017988\\
};
\addlegendentry{PCA (IT)}

\addplot [color=color2, mark=*, mark options={solid, color2, scale=0.75}]
  table[row sep=crcr]{%
0.1	0.0446853871056379\\
0.0909090909090909	0.0340693208052918\\
0.0833333333333333	0.0272774983588535\\
0.0769230769230769	0.026301108259373\\
0.0714285714285714	0.0236124487312959\\
0.0666666666666667	0.0187930381283514\\
0.0625	0.0174889805264242\\
0.0588235294117647	0.0165509389946752\\
0.0555555555555556	0.0139841403032615\\
0.0526315789473684	0.0124239527227388\\
0.05	0.0119797300969906\\
0.0476190476190476	0.0110254374061129\\
0.0454545454545455	0.0094613209511627\\
0.0434782608695652	0.00905203719244962\\
0.0416666666666667	0.00887329340674836\\
0.04	0.00759235545029391\\
0.0384615384615385	0.00711882295521415\\
0.037037037037037	0.00701829439686663\\
0.0357142857142857	0.00636148932620739\\
0.0344827586206897	0.00578804160605184\\
0.0333333333333333	0.00565833783468461\\
0.032258064516129	0.00547814621975862\\
0.03125	0.00484882068697501\\
0.0303030303030303	0.00465762306445439\\
0.0294117647058824	0.0046218959022768\\
0.0285714285714286	0.00417317307646092\\
0.0277777777777778	0.00392599427104989\\
0.027027027027027	0.00390521504364205\\
0.0263157894736842	0.00369181420694975\\
0.0256410256410256	0.00337895462516746\\
0.025	0.00332204180847073\\
0.024390243902439	0.00328648551181088\\
0.0238095238095238	0.00297232900433442\\
0.0232558139534884	0.00287167652195164\\
0.0227272727272727	0.00286911847077409\\
0.0222222222222222	0.00266271388825667\\
0.0217391304347826	0.00251557135281466\\
0.0212765957446809	0.00250903849314432\\
0.0208333333333333	0.00243742960372284\\
0.0204081632653061	0.00224485809248343\\
0.02	0.00220665199100001\\
0.0196078431372549	0.00220264800673409\\
0.0192307692307692	0.00203743575569471\\
0.0188679245283019	0.00195354320820251\\
0.0185185185185185	0.00196621747666081\\
0.0181818181818182	0.00187249927532562\\
0.0178571428571429	0.00175841671544807\\
0.0175438596491228	0.0017603424125241\\
0.0172413793103448	0.00173948035695815\\
0.0169491525423729	0.0016065411162054\\
0.0166666666666667	0.00157968647637508\\
0.0163934426229508	0.00158055705978044\\
0.0161290322580645	0.00149364090799442\\
0.0158730158730159	0.0014267891057349\\
0.015625	0.00143453643236868\\
0.0153846153846154	0.00139771520637044\\
0.0151515151515152	0.00130754550263812\\
0.0149253731343284	0.00130565800651117\\
0.0147058823529412	0.0013027162356698\\
0.0144927536231884	0.00121666925047847\\
0.0142857142857143	0.00119338471826613\\
0.0140845070422535	0.00119669127044464\\
0.0138888888888889	0.00114823873190684\\
0.0136986301369863	0.00109494564100676\\
0.0135135135135135	0.00109666754395721\\
0.0133333333333333	0.00108881558906448\\
0.0131578947368421	0.00101828677496374\\
0.012987012987013	0.00101070316388818\\
0.0128205128205128	0.0010163369321865\\
0.0126582278481013	0.000959119753398241\\
0.0125	0.000936265987886653\\
0.0123456790123457	0.000937950571307411\\
0.0121951219512195	0.000911976366565259\\
0.0120481927710843	0.000872374703724921\\
0.0119047619047619	0.000867954756807121\\
0.0117647058823529	0.00087086379595025\\
0.0116279069767442	0.000818847665242428\\
0.0114942528735632	0.000806119941038297\\
0.0113636363636364	0.000814518898792693\\
0.0112359550561798	0.000777851376607908\\
0.0111111111111111	0.000755163552450977\\
0.010989010989011	0.000758715803884202\\
0.0108695652173913	0.000744115546860602\\
0.010752688172043	0.000711953177380709\\
0.0106382978723404	0.000706833776573879\\
0.0105263157894737	0.000710068002279352\\
0.0104166666666667	0.000675945288551194\\
0.0103092783505155	0.000661111429497874\\
0.0102040816326531	0.000667783453232262\\
0.0101010101010101	0.000645994790742899\\
0.01	0.000624566328618625\\
};
\addlegendentry{Com (IT)}

\addplot [dashed]
  table[row sep=crcr]{%
0.1	1e-1\\
0.01	1e-03\\
};

\end{axis}
\end{tikzpicture}%

%% file: errors-BSeuND-PutOnAvg-set5-phi=0.0413-5D-r=0.04-yDomain-non-uniform-M=10to100.tikz
%
%
%
\begin{tikzpicture}

\begin{axis}[%
width=0.37\textwidth,
scale only axis,
xmode=log,
xmin=0.01,
xmax=0.1,
xminorticks=true,
ymode=log,
ymin=1e-4,
ymax=1e1,
yminorticks=true,
ylabel={error},
title={European basket option},
legend cell align={left},
legend style={nodes={scale=0.75, transform shape}}
]
\addplot [color=color1, mark=square*, mark options={solid, color1, scale=0.75}]
  table[row sep=crcr]{%
0.1	0.00232497495509065\\
0.0909090909090909	0.00132651740613987\\
0.0833333333333333	0.0150130180916115\\
0.0769230769230769	0.0159416032163087\\
0.0714285714285714	0.0252574741280451\\
0.0666666666666667	0.0245762917040846\\
0.0625	0.0282682780686603\\
0.0588235294117647	0.0259261086888714\\
0.0555555555555556	0.0265831414690645\\
0.0526315789473684	0.0237010778995757\\
0.05	0.0231542501423547\\
0.0476190476190476	0.020499984406019\\
0.0454545454545455	0.0195902469162288\\
0.0434782608695652	0.0174364999276726\\
0.0416666666666667	0.0165989757227774\\
0.04	0.0148960084897158\\
0.0384615384615385	0.0141692553654433\\
0.037037037037037	0.0128329580965652\\
0.0357142857142857	0.0122275664255032\\
0.0344827586206897	0.0111668446949569\\
0.0333333333333333	0.0106548129275676\\
0.032258064516129	0.00979326309006368\\
0.03125	0.00936067977366584\\
0.0303030303030303	0.00866318872094496\\
0.0294117647058824	0.00829151275523365\\
0.0285714285714286	0.00771219588673588\\
0.0277777777777778	0.00739784441392111\\
0.027027027027027	0.00691616230873926\\
0.0263157894736842	0.00664118316601581\\
0.0256410256410256	0.00623494203490793\\
0.025	0.00599991873908046\\
0.024390243902439	0.00565455022869088\\
0.0238095238095238	0.0054443753424831\\
0.0232558139534884	0.00514825554503773\\
0.0227272727272727	0.00496701350365747\\
0.0222222222222222	0.00470924687673424\\
0.0217391304347826	0.00454669989020218\\
0.0212765957446809	0.00432353767425475\\
0.0208333333333333	0.00417864273221902\\
0.0204081632653061	0.00398360426559385\\
0.02	0.00385430537253839\\
0.0196078431372549	0.003681698960456\\
0.0192307692307692	0.00356514185530799\\
0.0188679245283019	0.00341335959863542\\
0.0185185185185185	0.00330780459421298\\
0.0181818181818182	0.00317205579482227\\
0.0178571428571429	0.00307783499695269\\
0.0175438596491228	0.00295695578189203\\
0.0172413793103448	0.00287040721433662\\
0.0169491525423729	0.00276151802277269\\
0.0166666666666667	0.00268313896445882\\
0.0163934426229508	0.00258618870975802\\
0.0161290322580645	0.00251363565329221\\
0.0158730158730159	0.00242573162314663\\
0.015625	0.00236016103669456\\
0.0153846153846154	0.00228020360358272\\
0.0151515151515152	0.00221958600060645\\
0.0149253731343284	0.00214744879977324\\
0.0147058823529412	0.00209192846659967\\
0.0144927536231884	0.00202582025111475\\
0.0142857142857143	0.00197443615088133\\
0.0140845070422535	0.0019141111304517\\
0.0138888888888889	0.00186662011179628\\
0.0136986301369863	0.00181164703350589\\
0.0135135135135135	0.00176765188629169\\
0.0133333333333333	0.00171683117309485\\
0.0131578947368421	0.00167604965855084\\
0.012987012987013	0.00162958392769852\\
0.0128205128205128	0.001591548351076\\
0.0126582278481013	0.00154831508659292\\
0.0125	0.00151323512918466\\
0.0123456790123457	0.00147340882866231\\
0.0121951219512195	0.00144034764905321\\
0.0120481927710843	0.00140342890988876\\
0.0119047619047619	0.00137268753068476\\
0.0117647058823529	0.00133851560869758\\
0.0116279069767442	0.00130953712860382\\
0.0114942528735632	0.00127776336101526\\
0.0113636363636364	0.00125095444966128\\
0.0112359550561798	0.00122122914756773\\
0.0111111111111111	0.00119587346094185\\
0.010989010989011	0.00116828006753344\\
0.0108695652173913	0.00114447458740408\\
0.010752688172043	0.00111865367250452\\
0.0106382978723404	0.00109633049831714\\
0.0105263157894737	0.0010721769530484\\
0.0104166666666667	0.00105102183998085\\
0.0103092783505155	0.00102842119802737\\
0.0102040816326531	0.00100849391086832\\
0.0101010101010101	0.000987264465477367\\
0.01	0.000968523399402343\\
};
\addlegendentry{PCA}

\addplot [color=color2, mark=*, mark options={solid, color2, scale=0.75}]
  table[row sep=crcr]{%
0.1	0.0117515824491026\\
0.0909090909090909	0.00192483136158117\\
0.0833333333333333	0.00785843203532188\\
0.0769230769230769	0.0152781141480212\\
0.0714285714285714	0.0206707238803201\\
0.0666666666666667	0.0236747972796305\\
0.0625	0.0250601065591258\\
0.0588235294117647	0.0250234949530146\\
0.0555555555555556	0.024237781210509\\
0.0526315789473684	0.0228818783219253\\
0.05	0.021375758383769\\
0.0476190476190476	0.0197609078863774\\
0.0454545454545455	0.0182552917908971\\
0.0434782608695652	0.0168099251972436\\
0.0416666666666667	0.0155337096021793\\
0.04	0.0143477355388679\\
0.0384615384615385	0.0133146442431382\\
0.037037037037037	0.0123580254142119\\
0.0357142857142857	0.0115231324802529\\
0.0344827586206897	0.0107458939716629\\
0.0333333333333333	0.010063835658757\\
0.032258064516129	0.00942503323781096\\
0.03125	0.00886170571123834\\
0.0303030303030303	0.00833148600435507\\
0.0294117647058824	0.0078619657861978\\
0.0285714285714286	0.00741808303987135\\
0.0277777777777778	0.00702337562062993\\
0.027027027027027	0.0066485672081491\\
0.0263157894736842	0.00631377265905719\\
0.0256410256410256	0.00599444571100349\\
0.025	0.00570785871035007\\
0.024390243902439	0.00543338808599714\\
0.0238095238095238	0.00518591700416593\\
0.0232558139534884	0.00494807167812183\\
0.0227272727272727	0.00473269671925713\\
0.0222222222222222	0.0045251127089454\\
0.0217391304347826	0.00433639468588964\\
0.0212765957446809	0.00415409564027769\\
0.0208333333333333	0.00398775718267608\\
0.0204081632653061	0.00382679659276519\\
0.02	0.00367941753020418\\
0.0196078431372549	0.00353660914304044\\
0.0192307692307692	0.00340540908160136\\
0.0188679245283019	0.00327813978431024\\
0.0185185185185185	0.00316082966004316\\
0.0181818181818182	0.00304693851277626\\
0.0178571428571429	0.00294161764871248\\
0.0175438596491228	0.00283930294508394\\
0.0172413793103448	0.0027443834042149\\
0.0169491525423729	0.00265213480976634\\
0.0166666666666667	0.0025662834614657\\
0.0163934426229508	0.0024828292284873\\
0.0161290322580645	0.00240491930665421\\
0.0158730158730159	0.00232918210982724\\
0.015625	0.0022582575404666\\
0.0153846153846154	0.00218932034727359\\
0.0151515151515152	0.0021245652531956\\
0.0149253731343284	0.00206164310589863\\
0.0147058823529412	0.0020023578299071\\
0.0144927536231884	0.00194477556581774\\
0.0142857142857143	0.00189035705246177\\
0.0140845070422535	0.00183753129425845\\
0.0138888888888889	0.00178745754575793\\
0.0136986301369863	0.00173888212717621\\
0.0135135135135135	0.00169269942584371\\
0.0133333333333333	0.00164793346997349\\
0.0131578947368421	0.00160524594799405\\
0.012987012987013	0.00156390424406894\\
0.0128205128205128	0.00152436508786913\\
0.0126582278481013	0.00148610945565153\\
0.0125	0.00144941423853417\\
0.0123456790123457	0.00141394716626664\\
0.0121951219512195	0.00137982745021148\\
0.0120481927710843	0.00134688616437728\\
0.0119047619047619	0.00131510471856267\\
0.0117647058823529	0.0012844568595618\\
0.0116279069767442	0.00125480297957914\\
0.0114942528735632	0.00122624202550448\\
0.0113636363636364	0.00119852849790014\\
0.0112359550561798	0.00117187055389145\\
0.0111111111111111	0.00114593043432087\\
0.010989010989011	0.00112101131023579\\
0.0108695652173913	0.00109669531257772\\
0.010752688172043	0.00107336925575119\\
0.0106382978723404	0.00105054238008506\\
0.0105263157894737	0.00102867524144301\\
0.0104166666666667	0.00100721958683714\\
0.0103092783505155	0.000986694030048696\\
0.0102040816326531	0.000966500129285564\\
0.0101010101010101	0.000947209345712352\\
0.01	0.00092817948635826\\
};
\addlegendentry{Com}

\addplot [dashed]
  table[row sep=crcr]{%
0.1	1e-0\\
0.01	1e-2\\
};

\end{axis}
\end{tikzpicture}%

%% file: errors-BSamND-PutOnAvg-set5-phi=0.0413-5D-r=0.04-yDomain-non-uniform-M=10to100.tikz
%
%
%
\begin{tikzpicture}

\begin{axis}[%
width=0.37\textwidth,
scale only axis,
xmode=log,
xmin=0.01,
xmax=0.1,
xminorticks=true,
ymode=log,
ymin=1e-4,
ymax=1e1,
yminorticks=true,
title={American basket option},
legend cell align={left},
legend style={nodes={scale=0.75, transform shape}}
]
\addplot [color=color1, mark=square*, mark options={solid, color1, scale=0.75}]
  table[row sep=crcr]{%
0.1	0.0354246703716772\\
0.0909090909090909	0.0596003655276007\\
0.0833333333333333	0.0561018508880888\\
0.0769230769230769	0.0512899834870471\\
0.0714285714285714	0.0286028023142517\\
0.0666666666666667	0.0164757486906684\\
0.0625	0.0108517556400969\\
0.0588235294117647	0.0170935939411301\\
0.0555555555555556	0.0202184709717041\\
0.0526315789473684	0.0261198050262621\\
0.05	0.0250462826674234\\
0.0476190476190476	0.0229778433070642\\
0.0454545454545455	0.0160175399443148\\
0.0434782608695652	0.0111796745571784\\
0.0416666666666667	0.00769660083057921\\
0.04	0.00917224085494794\\
0.0384615384615385	0.010110499364945\\
0.037037037037037	0.0123401400527339\\
0.0357142857142857	0.0125934433457413\\
0.0344827586206897	0.0124092577038173\\
0.0333333333333333	0.0101345703327702\\
0.032258064516129	0.0082102724826143\\
0.03125	0.0060660212534831\\
0.0303030303030303	0.00573298644820497\\
0.0294117647058824	0.00591742061694056\\
0.0285714285714286	0.00689851593937973\\
0.0277777777777778	0.00725489060143403\\
0.027027027027027	0.00740590290696552\\
0.0263157894736842	0.00679864928295615\\
0.0256410256410256	0.00608359229913802\\
0.025	0.00464946291672241\\
0.024390243902439	0.0041882232116297\\
0.0238095238095238	0.00403386578549103\\
0.0232558139534884	0.00425028710699848\\
0.0227272727272727	0.00438609095645681\\
0.0222222222222222	0.00488592987959713\\
0.0217391304347826	0.0046246730771724\\
0.0212765957446809	0.00432262690330809\\
0.0208333333333333	0.00382807382379724\\
0.0204081632653061	0.00332998827906295\\
0.02	0.00278614782400588\\
0.0196078431372549	0.00289180600434413\\
0.0192307692307692	0.00308469622856755\\
0.0188679245283019	0.00316555658735673\\
0.0185185185185185	0.00325372660940459\\
0.0181818181818182	0.00323360621053226\\
0.0178571428571429	0.00291734838498137\\
0.0175438596491228	0.00264241072023452\\
0.0172413793103448	0.00236390049466628\\
0.0169491525423729	0.00221938206003891\\
0.0166666666666667	0.00215053232142637\\
0.0163934426229508	0.00223619143590703\\
0.0161290322580645	0.0023326055834989\\
0.0158730158730159	0.00235811346492198\\
0.015625	0.00226141639697275\\
0.0153846153846154	0.00214506020880556\\
0.0151515151515152	0.0019516322014681\\
0.0149253731343284	0.00182767333483369\\
0.0147058823529412	0.00166306390434023\\
0.0144927536231884	0.00169433859132084\\
0.0142857142857143	0.00171911864575236\\
0.0140845070422535	0.00174503190785202\\
0.0138888888888889	0.00173844077025009\\
0.0136986301369863	0.00175537159885053\\
0.0135135135135135	0.00160595849314937\\
0.0133333333333333	0.00149881940307317\\
0.0131578947368421	0.0014240517084243\\
0.012987012987013	0.00134970722759853\\
0.0128205128205128	0.00129589160914634\\
0.0126582278481013	0.00134262334994872\\
0.0125	0.00136335869262538\\
0.0123456790123457	0.00135698483751945\\
0.0121951219512195	0.00129524818997773\\
0.0120481927710843	0.001277226934155\\
0.0119047619047619	0.00118801319460538\\
0.0117647058823529	0.001104076840301\\
0.0116279069767442	0.00107315503835181\\
0.0114942528735632	0.00107310895812418\\
0.0113636363636364	0.00106808766575384\\
0.0112359550561798	0.00107320388295129\\
0.0111111111111111	0.00105940320110953\\
0.010989010989011	0.00106387415629428\\
0.0108695652173913	0.000998709196906589\\
0.010752688172043	0.000950710502314323\\
0.0106382978723404	0.000910559444960768\\
0.0105263157894737	0.000884207492347855\\
0.0104166666666667	0.00085376357592537\\
0.0103092783505155	0.000858225502138765\\
0.0102040816326531	0.00088342427238608\\
0.0101010101010101	0.000867512059501863\\
0.01	0.00084605082484579\\
};
\addlegendentry{PCA (IT)}

\addplot [color=color2, mark=*, mark options={solid, color2, scale=0.75}]
  table[row sep=crcr]{%
0.1	0.0507025028083277\\
0.0909090909090909	0.0625269467901646\\
0.0833333333333333	0.06203679782808\\
0.0769230769230769	0.050688356687365\\
0.0714285714285714	0.0312432675172722\\
0.0666666666666667	0.0166102585969963\\
0.0625	0.0144831983781746\\
0.0588235294117647	0.0186153195866492\\
0.0555555555555556	0.0227367542790979\\
0.0526315789473684	0.0268527589383698\\
0.05	0.0260072188047946\\
0.0476190476190476	0.022618076538194\\
0.0454545454545455	0.0158123670238943\\
0.0434782608695652	0.0111310293986975\\
0.0416666666666667	0.00880096767184391\\
0.04	0.00963980251128227\\
0.0384615384615385	0.0116070763841556\\
0.037037037037037	0.0131670792802883\\
0.0357142857142857	0.0132115347527844\\
0.0344827586206897	0.0122375445564895\\
0.0333333333333333	0.0102556133063221\\
0.032258064516129	0.00778770777956872\\
0.03125	0.00627323235426758\\
0.0303030303030303	0.00606457830001084\\
0.0294117647058824	0.00653780072545196\\
0.0285714285714286	0.00729507970824805\\
0.0277777777777778	0.00767066764347746\\
0.027027027027027	0.00739077029692758\\
0.0263157894736842	0.00682250978533361\\
0.0256410256410256	0.00576350753873633\\
0.025	0.00476226169818439\\
0.024390243902439	0.00429842162340854\\
0.0238095238095238	0.00430798865762583\\
0.0232558139534884	0.00440380910236904\\
0.0227272727272727	0.00476279764210474\\
0.0222222222222222	0.00492848080032893\\
0.0217391304347826	0.00459886991346714\\
0.0212765957446809	0.00419864894356614\\
0.0208333333333333	0.00382167032552516\\
0.0204081632653061	0.00324026848841008\\
0.02	0.00297664253206698\\
0.0196078431372549	0.00310603699687739\\
0.0192307692307692	0.00321125173380565\\
0.0188679245283019	0.0032529093321596\\
0.0185185185185185	0.00332439414285446\\
0.0181818181818182	0.00326603052521435\\
0.0178571428571429	0.00289086399187433\\
0.0175438596491228	0.00260121398974711\\
0.0172413793103448	0.00245208826098064\\
0.0169491525423729	0.00228874900929554\\
0.0166666666666667	0.00224532918339015\\
0.0163934426229508	0.00239376669986591\\
0.0161290322580645	0.00246792657835471\\
0.0158730158730159	0.00236211782167395\\
0.015625	0.00227981851592673\\
0.0153846153846154	0.00215510385903173\\
0.0151515151515152	0.00193059138085339\\
0.0149253731343284	0.00182331006589465\\
0.0147058823529412	0.00176304401473359\\
0.0144927536231884	0.00180768786268715\\
0.0142857142857143	0.00178449325450991\\
0.0140845070422535	0.00182302309477578\\
0.0138888888888889	0.00180704252443675\\
0.0136986301369863	0.00173777755867022\\
0.0135135135135135	0.0015717973524918\\
0.0133333333333333	0.00151096669387485\\
0.0131578947368421	0.00144094889787283\\
0.012987012987013	0.00140523143233118\\
0.0128205128205128	0.00137400235301932\\
0.0126582278481013	0.00142367546778033\\
0.0125	0.00140027115192431\\
0.0123456790123457	0.00137353221881575\\
0.0121951219512195	0.00131103219318263\\
0.0120481927710843	0.00128785564603184\\
0.0119047619047619	0.00116798207585767\\
0.0117647058823529	0.00113451843655454\\
0.0116279069767442	0.0011139206668993\\
0.0114942528735632	0.00112472171138123\\
0.0113636363636364	0.00109242247216379\\
0.0112359550561798	0.00110969140910377\\
0.0111111111111111	0.0010989225230329\\
0.010989010989011	0.00106105849992133\\
0.0108695652173913	0.000981354286890834\\
0.010752688172043	0.000975206846733911\\
0.0106382978723404	0.000936177074867928\\
0.0105263157894737	0.000896375053304865\\
0.0104166666666667	0.000892449119053396\\
0.0103092783505155	0.00092632432630424\\
0.0102040816326531	0.000887588499328018\\
0.0101010101010101	0.000883244921684678\\
0.01	0.000866270652695178\\
};
\addlegendentry{Com (IT)}

\addplot [dashed]
  table[row sep=crcr]{%
0.1	1e-0\\
0.01	1e-2\\
};

\end{axis}
\end{tikzpicture}%

%% file: errors-BSeuND-PutOnAvg-set5-phi=0.0413-10D-r=0.04-yDomain-non-uniform-M=10to100.tikz
%
%
%
\begin{tikzpicture}

\begin{axis}[%
width=0.37\textwidth,
scale only axis,
xmode=log,
xmin=0.01,
xmax=0.1,
xminorticks=true,
ymode=log,
ymin=1e-4,
ymax=1e1,
yminorticks=true,
ylabel={error},
legend cell align={left},
legend style={nodes={scale=0.75, transform shape}}
]
\addplot [color=color1, mark=square*, mark options={solid, color1, scale=0.75}]
  table[row sep=crcr]{%
0.1	0.327181118337579\\
0.0909090909090909	0.244015728592277\\
0.0833333333333333	0.210899412338392\\
0.0769230769230769	0.175797321506467\\
0.0714285714285714	0.162581128670864\\
0.0666666666666667	0.145195699186896\\
0.0625	0.13829646673158\\
0.0588235294117647	0.127486721554632\\
0.0555555555555556	0.122329282446765\\
0.0526315789473684	0.114138766548779\\
0.05	0.1091932419685\\
0.0476190476190476	0.101993185179772\\
0.0454545454545455	0.097068208114111\\
0.0434782608695652	0.0905293806942336\\
0.0416666666666667	0.0856433205377449\\
0.04	0.0797684078752141\\
0.0384615384615385	0.0752063221071317\\
0.037037037037037	0.0700376450855007\\
0.0357142857142857	0.0659336779000519\\
0.0344827586206897	0.0615002403970841\\
0.0333333333333333	0.0579290828172798\\
0.032258064516129	0.0541472449756757\\
0.03125	0.0511043384990248\\
0.0303030303030303	0.0479385393588796\\
0.0294117647058824	0.0453575110983824\\
0.0285714285714286	0.0426943820768049\\
0.0277777777777778	0.0405067733784517\\
0.027027027027027	0.0382622886094897\\
0.0263157894736842	0.0363946798919734\\
0.0256410256410256	0.0344715261530268\\
0.025	0.0328771160486436\\
0.024390243902439	0.0312287993356133\\
0.0238095238095238	0.0298427288231782\\
0.0232558139534884	0.0284209232597092\\
0.0227272727272727	0.0272177395785302\\
0.0222222222222222	0.0259809846343497\\
0.0217391304347826	0.0249206272536693\\
0.0212765957446809	0.0238335612292104\\
0.0208333333333333	0.0229028365982149\\
0.0204081632653061	0.0219472407151606\\
0.02	0.0211192317034836\\
0.0196078431372549	0.0202764893243987\\
0.0192307692307692	0.0195411870793514\\
0.0188679245283019	0.018790884121815\\
0.0185185185185185	0.0181315097707841\\
0.0181818181818182	0.0174604069526261\\
0.0178571428571429	0.0168704340723931\\
0.0175438596491228	0.0162678151604805\\
0.0172413793103448	0.015735405491105\\
0.0169491525423729	0.0151957951237289\\
0.0166666666666667	0.0147142761716221\\
0.0163934426229508	0.0142238972147268\\
0.0161290322580645	0.013786422085122\\
0.0158730158730159	0.0133429238350136\\
0.015625	0.01294321994577\\
0.0153846153846154	0.0125397045165538\\
0.0151515151515152	0.0121775072234023\\
0.0149253731343284	0.0118090855035735\\
0.0147058823529412	0.0114755675728215\\
0.0144927536231884	0.0111387396970972\\
0.0142857142857143	0.0108327059108078\\
0.0140845070422535	0.0105243390299989\\
0.0138888888888889	0.0102419216764478\\
0.0136986301369863	0.00995768708265032\\
0.0135135135135135	0.00969919578382061\\
0.0133333333333333	0.00943769932877103\\
0.0131578947368421	0.00919736758596379\\
0.012987012987013	0.0089562615789589\\
0.0128205128205128	0.00873336100941735\\
0.0126582278481013	0.00851027265890103\\
0.0125	0.00830332182562188\\
0.0123456790123457	0.00809651480655127\\
0.0121951219512195	0.00790510034857483\\
0.0120481927710843	0.00771241309211845\\
0.0119047619047619	0.00753314936956784\\
0.0117647058823529	0.00735508944725538\\
0.0116279069767442	0.00718785865968918\\
0.0114942528735632	0.00702117111944389\\
0.0113636363636364	0.00686503041411335\\
0.0112359550561798	0.0067099764698213\\
0.0111111111111111	0.00656370858746413\\
0.010989010989011	0.00641849805474415\\
0.0108695652173913	0.00628171285628554\\
0.010752688172043	0.00614560646123152\\
0.0106382978723404	0.00601703263799536\\
0.0105263157894737	0.00588931030217665\\
0.0104166666666667	0.00576866826342837\\
0.0103092783505155	0.00564909395323276\\
0.0102040816326531	0.0055353153723221\\
0.0101010101010101	0.00542295311104546\\
0.01	0.00531599851477083\\
};
\addlegendentry{PCA}

\addplot [color=color2, mark=*, mark options={solid, color2, scale=0.75}]
  table[row sep=crcr]{%
0.1	0.304661066988253\\
0.0909090909090909	0.23918601538961\\
0.0833333333333333	0.19931446260485\\
0.0769230769230769	0.172849374636309\\
0.0714285714285714	0.155442383681356\\
0.0666666666666667	0.142476546142772\\
0.0625	0.13288008625987\\
0.0588235294117647	0.124682110853602\\
0.0555555555555556	0.117781710898232\\
0.0526315789473684	0.11121906968242\\
0.05	0.105195964077916\\
0.0476190476190476	0.0992077624446974\\
0.0454545454545455	0.0935401384446877\\
0.0434782608695652	0.0879064284047342\\
0.0416666666666667	0.0825844250776679\\
0.04	0.0773842711663342\\
0.0384615384615385	0.0725359366649769\\
0.037037037037037	0.0678872022833552\\
0.0357142857142857	0.0636118202988687\\
0.0344827586206897	0.0595723286801757\\
0.0333333333333333	0.0558944722379753\\
0.032258064516129	0.052450802722154\\
0.03125	0.0493325592032878\\
0.0303030303030303	0.0464238599254649\\
0.0294117647058824	0.0437936731660944\\
0.0285714285714286	0.0413396722068633\\
0.0277777777777778	0.0391171817892548\\
0.027027027027027	0.0370379362930855\\
0.0263157894736842	0.0351487784232791\\
0.0256410256410256	0.0333745709157416\\
0.025	0.0317564026654331\\
0.024390243902439	0.0302305373927734\\
0.0238095238095238	0.0288336484650138\\
0.0232558139534884	0.0275114835778236\\
0.0227272727272727	0.0262969738187486\\
0.0222222222222222	0.0251436029118004\\
0.0217391304347826	0.0240809751061093\\
0.0212765957446809	0.0230688693587702\\
0.0208333333333333	0.0221339026266563\\
0.0204081632653061	0.0212410189458687\\
0.02	0.0204141684572274\\
0.0196078431372549	0.0196225874708968\\
0.0192307692307692	0.0188878659078142\\
0.0188679245283019	0.018182849480425\\
0.0185185185185185	0.0175270515150192\\
0.0181818181818182	0.01689638806684\\
0.0178571428571429	0.0163085423545279\\
0.0175438596491228	0.0157420623322118\\
0.0172413793103448	0.0152130164074968\\
0.0169491525423729	0.0147022199242564\\
0.0166666666666667	0.0142243111785181\\
0.0163934426229508	0.0137620666843237\\
0.0161290322580645	0.013328852987085\\
0.0158730158730159	0.0129091504316907\\
0.015625	0.0125151893207818\\
0.0153846153846154	0.0121329361926268\\
0.0151515151515152	0.0117736032313243\\
0.0149253731343284	0.011424459657567\\
0.0147058823529412	0.011095802665217\\
0.0144927536231884	0.0107760476249652\\
0.0142857142857143	0.0104746685497421\\
0.0140845070422535	0.0101810938185931\\
0.0138888888888889	0.00990405502200531\\
0.0136986301369863	0.00963388174010049\\
0.0135135135135135	0.00937862986120308\\
0.0133333333333333	0.00912943325878857\\
0.0131578947368421	0.00889374576383339\\
0.012987012987013	0.0086634116481612\\
0.0128205128205128	0.00844533698464112\\
0.0126582278481013	0.00823200851041328\\
0.0125	0.00802983398726198\\
0.0123456790123457	0.00783187705506094\\
0.0121951219512195	0.00764409227915053\\
0.0120481927710843	0.00746006462937565\\
0.0119047619047619	0.00728533454487845\\
0.0117647058823529	0.00711395870274067\\
0.0116279069767442	0.00695109983393216\\
0.0114942528735632	0.0067912419227909\\
0.0113636363636364	0.00663920205385082\\
0.0112359550561798	0.00648985235639366\\
0.0111111111111111	0.00634769290447323\\
0.010989010989011	0.00620795002261865\\
0.0108695652173913	0.0060748329852931\\
0.010752688172043	0.00594389046658961\\
0.0106382978723404	0.00581906431394508\\
0.0105263157894737	0.0056961983311119\\
0.0104166666666667	0.00557898837676341\\
0.0103092783505155	0.00546354824794726\\
0.0102040816326531	0.00535334661763365\\
0.0101010101010101	0.00524474585996693\\
0.01	0.00514100464309664\\
};
\addlegendentry{Com}

\addplot [dashed]
  table[row sep=crcr]{%
0.1	1e-0\\
0.01	1e-2\\
};

\end{axis}
\end{tikzpicture}%

%% file: errors-BSamND-PutOnAvg-set5-phi=0.0413-10D-r=0.04-yDomain-non-uniform-M=10to100.tikz
%
%
%
\begin{tikzpicture}

\begin{axis}[%
width=0.37\textwidth,
scale only axis,
xmode=log,
xmin=0.01,
xmax=0.1,
xminorticks=true,
ymode=log,
ymin=1e-4,
ymax=1e1,
yminorticks=true,
legend cell align={left},
legend style={nodes={scale=0.75, transform shape}}
]
\addplot [color=color1, mark=square*, mark options={solid, color1, scale=0.75}]
  table[row sep=crcr]{%
0.1	0.247535458487862\\
0.0909090909090909	0.221651562216804\\
0.0833333333333333	0.198207479395707\\
0.0769230769230769	0.156896594479939\\
0.0714285714285714	0.127136258993556\\
0.0666666666666667	0.0930154090720396\\
0.0625	0.0726199518233805\\
0.0588235294117647	0.059061116450934\\
0.0555555555555556	0.0541863280087167\\
0.0526315789473684	0.0539609019913154\\
0.05	0.0572389118714547\\
0.0476190476190476	0.0553297778264863\\
0.0454545454545455	0.0518743368481633\\
0.0434782608695652	0.0436446126811418\\
0.0416666666666667	0.037264009819868\\
0.04	0.0295191450110561\\
0.0384615384615385	0.0250331468404283\\
0.037037037037037	0.0210992433451072\\
0.0357142857142857	0.0199633545692421\\
0.0344827586206897	0.020116586762347\\
0.0333333333333333	0.0215476711448108\\
0.032258064516129	0.0225616354397236\\
0.03125	0.0229258171370894\\
0.0303030303030303	0.021346130202657\\
0.0294117647058824	0.0191078563791667\\
0.0285714285714286	0.0166573497562457\\
0.0277777777777778	0.0140375151112249\\
0.027027027027027	0.0120830683889981\\
0.0263157894736842	0.0110742152594181\\
0.0256410256410256	0.0105692181558297\\
0.025	0.0105104495384882\\
0.024390243902439	0.011085755412795\\
0.0238095238095238	0.0118012523859754\\
0.0232558139534884	0.0117546520745346\\
0.0227272727272727	0.0115618286790014\\
0.0222222222222222	0.0105795828545432\\
0.0217391304347826	0.00958658771978804\\
0.0212765957446809	0.00847275395160985\\
0.0208333333333333	0.00758084322132468\\
0.0204081632653061	0.00700595758208422\\
0.02	0.00669265979539446\\
0.0196078431372549	0.00657020356361038\\
0.0192307692307692	0.00662347978472866\\
0.0188679245283019	0.00696429576771074\\
0.0185185185185185	0.00711132419866622\\
0.0181818181818182	0.00705731826303158\\
0.0178571428571429	0.00689498826839419\\
0.0175438596491228	0.00634811843749183\\
0.0172413793103448	0.00584934980209084\\
0.0169491525423729	0.00528935782990914\\
0.0166666666666667	0.00493617722003847\\
0.0163934426229508	0.00469182937932366\\
0.0161290322580645	0.00456427286554018\\
0.0158730158730159	0.00453590657594738\\
0.015625	0.00460271930065836\\
0.0153846153846154	0.00474245510612192\\
0.0151515151515152	0.00475465084891624\\
0.0149253731343284	0.00468771696344916\\
0.0147058823529412	0.00445098207812578\\
0.0144927536231884	0.00426137180434694\\
0.0142857142857143	0.00395836236055125\\
0.0140845070422535	0.00365659995106782\\
0.0138888888888889	0.00344967674649688\\
0.0136986301369863	0.0033659586947028\\
0.0135135135135135	0.00329622616887448\\
0.0133333333333333	0.00329903098218232\\
0.0131578947368421	0.00332524236433329\\
0.012987012987013	0.00335144710785507\\
0.0128205128205128	0.00338345042845134\\
0.0126582278481013	0.00330647628270242\\
0.0125	0.00319495778069268\\
0.0123456790123457	0.00304681626593784\\
0.0121951219512195	0.00285199523714574\\
0.0120481927710843	0.00269338259003149\\
0.0119047619047619	0.002606494841908\\
0.0117647058823529	0.00257213273078882\\
0.0116279069767442	0.00247100546140899\\
0.0114942528735632	0.00249100377112299\\
0.0113636363636364	0.00250910497630308\\
0.0112359550561798	0.00249105999716726\\
0.0111111111111111	0.00256190792136124\\
0.010989010989011	0.00247928969587363\\
0.0108695652173913	0.00240579989442224\\
0.010752688172043	0.00230789346072235\\
0.0106382978723404	0.00220389700233881\\
0.0105263157894737	0.00210625044489987\\
0.0104166666666667	0.00199586153191156\\
0.0103092783505155	0.00199468479423004\\
0.0102040816326531	0.0019741737250385\\
0.0101010101010101	0.00196795571834407\\
0.01	0.00195819173789502\\
};
\addlegendentry{PCA (IT)}

\addplot [color=color2, mark=*, mark options={solid, color2, scale=0.75}]
  table[row sep=crcr]{%
0.1	0.22271975593709\\
0.0909090909090909	0.207926105483505\\
0.0833333333333333	0.178687804792567\\
0.0769230769230769	0.142824115178763\\
0.0714285714285714	0.111637331375071\\
0.0666666666666667	0.0833359560532738\\
0.0625	0.0643808755209232\\
0.0588235294117647	0.0549785437295487\\
0.0555555555555556	0.0514382246129319\\
0.0526315789473684	0.0525302502576519\\
0.05	0.0540118276531469\\
0.0476190476190476	0.052372515382908\\
0.0454545454545455	0.046788095037785\\
0.0434782608695652	0.0394691749402671\\
0.0416666666666667	0.0320059259045262\\
0.04	0.0264032156407943\\
0.0384615384615385	0.0219922868297804\\
0.037037037037037	0.0194375990024458\\
0.0357142857142857	0.0194926047616804\\
0.0344827586206897	0.0201265534675787\\
0.0333333333333333	0.0216275624931743\\
0.032258064516129	0.0217206596118533\\
0.03125	0.0211214066378105\\
0.0303030303030303	0.0191635905324308\\
0.0294117647058824	0.0166263292117801\\
0.0285714285714286	0.014468719110214\\
0.0277777777777778	0.0123234884211829\\
0.027027027027027	0.0111761947691011\\
0.0263157894736842	0.010320943060746\\
0.0256410256410256	0.0102788372006213\\
0.025	0.0105712123736144\\
0.024390243902439	0.0110389773290915\\
0.0238095238095238	0.0113933181534449\\
0.0232558139534884	0.0109245675149818\\
0.0227272727272727	0.0103088721001825\\
0.0222222222222222	0.00942213445887852\\
0.0217391304347826	0.00844988651729217\\
0.0212765957446809	0.00744437826183564\\
0.0208333333333333	0.0067967013196899\\
0.0204081632653061	0.00659287894819549\\
0.02	0.00637465430111206\\
0.0196078431372549	0.0065009521104038\\
0.0192307692307692	0.00668127776496341\\
0.0188679245283019	0.0068496273456784\\
0.0185185185185185	0.00679918841887428\\
0.0181818181818182	0.00637128876700999\\
0.0178571428571429	0.00598636712319767\\
0.0175438596491228	0.00557338874779845\\
0.0172413793103448	0.00514368000971999\\
0.0169491525423729	0.00470779146264633\\
0.0166666666666667	0.0045021871558788\\
0.0163934426229508	0.00447898943687974\\
0.0161290322580645	0.00438851067580259\\
0.0158730158730159	0.00444173280397031\\
0.015625	0.00450455794145554\\
0.0153846153846154	0.00455074486496798\\
0.0151515151515152	0.00443829487928227\\
0.0149253731343284	0.00417826545371369\\
0.0147058823529412	0.00397758462110964\\
0.0144927536231884	0.00371652103928355\\
0.0142857142857143	0.00352538837670657\\
0.0140845070422535	0.00330200309327999\\
0.0138888888888889	0.00322617108495038\\
0.0136986301369863	0.00322819313715605\\
0.0135135135135135	0.00321572588746299\\
0.0133333333333333	0.00322351276233412\\
0.0131578947368421	0.00322230695413239\\
0.012987012987013	0.00324056380319426\\
0.0128205128205128	0.00311222164351399\\
0.0126582278481013	0.00296740927763572\\
0.0125	0.00282804713098628\\
0.0123456790123457	0.00269159998419433\\
0.0121951219512195	0.00258882880418021\\
0.0120481927710843	0.00248222889750771\\
0.0119047619047619	0.00244064028750479\\
0.0117647058823529	0.00244767410775459\\
0.0116279069767442	0.00243931167267775\\
0.0114942528735632	0.00244886380361109\\
0.0113636363636364	0.00240650284326271\\
0.0112359550561798	0.00239735965201682\\
0.0111111111111111	0.00231959181542685\\
0.010989010989011	0.00223018825069587\\
0.0108695652173913	0.00212112877464676\\
0.010752688172043	0.00203883672415195\\
0.0106382978723404	0.00201594098530222\\
0.0105263157894737	0.00193613278208638\\
0.0104166666666667	0.0019096579826634\\
0.0103092783505155	0.00192181233422062\\
0.0102040816326531	0.00191747836572809\\
0.0101010101010101	0.0018930035717819\\
0.01	0.00185293487898264\\
};
\addlegendentry{Com (IT)}

\addplot [dashed]
  table[row sep=crcr]{%
0.1	1e-0\\
0.01	1e-2\\
};

\end{axis}
\end{tikzpicture}%

%% file: errors-BSeuND-PutOnAvg-set5-phi=0.0413-15D-r=0.04-yDomain-non-uniform-M=10to100.tikz
%
%
%
\begin{tikzpicture}

\begin{axis}[%
width=0.37\textwidth,
scale only axis,
xmode=log,
xmin=0.01,
xmax=0.1,
xminorticks=true,
xlabel={1/m},
ymode=log,
ymin=1e-4,
ymax=1e1,
yminorticks=true,
ylabel={error},
legend cell align={left},
legend style={nodes={scale=0.75, transform shape}}
]
\addplot [color=color1, mark=square*, mark options={solid, color1, scale=0.75}]
  table[row sep=crcr]{%
0.1	0.756614312300528\\
0.0909090909090909	0.548109880555185\\
0.0833333333333333	0.432718818632951\\
0.0769230769230769	0.342431418365701\\
0.0714285714285714	0.2931565783887\\
0.0666666666666667	0.25062812827489\\
0.0625	0.226956529961191\\
0.0588235294117647	0.20393180719244\\
0.0555555555555556	0.19038145819311\\
0.0526315789473684	0.175657842965105\\
0.05	0.166195477693137\\
0.0476190476190476	0.155384983233017\\
0.0454545454545455	0.147756067699159\\
0.0434782608695652	0.138940076039368\\
0.0416666666666667	0.132326993256099\\
0.04	0.124724605350648\\
0.0384615384615385	0.118557533536432\\
0.037037037037037	0.111730178081954\\
0.0357142857142857	0.106068393340603\\
0.0344827586206897	0.0999452865789401\\
0.0333333333333333	0.0947868475480895\\
0.032258064516129	0.0893248740588106\\
0.03125	0.0846747094060341\\
0.0303030303030303	0.0798565330868595\\
0.0294117647058824	0.0757102245390104\\
0.0285714285714286	0.0715079936078169\\
0.0277777777777778	0.0678666569738624\\
0.027027027027027	0.0642443202128664\\
0.0263157894736842	0.0610731785005694\\
0.0256410256410256	0.0579423051080852\\
0.025	0.0551923915818655\\
0.024390243902439	0.0524849789673727\\
0.0238095238095238	0.0501050579916669\\
0.0232558139534884	0.0477480849072123\\
0.0227272727272727	0.0456762703120202\\
0.0222222222222222	0.0436330256401174\\
0.0217391304347826	0.0418213147585611\\
0.0212765957446809	0.0400205798639774\\
0.0208333333333333	0.0384196099616521\\
0.0204081632653061	0.036839795820125\\
0.02	0.0354229704950644\\
0.0196078431372549	0.0340296115029481\\
0.0192307692307692	0.0327695255795817\\
0.0188679245283019	0.0315264813521097\\
0.0185185185185185	0.0303988208071448\\
0.0181818181818182	0.0292828969500611\\
0.0178571428571429	0.0282760891389611\\
0.0175438596491228	0.0272791240373778\\
0.0172413793103448	0.0263717860452797\\
0.0169491525423729	0.0254698854160065\\
0.0166666666666667	0.0246497611441647\\
0.0163934426229508	0.0238402413454768\\
0.0161290322580645	0.0230934117031314\\
0.0158730158730159	0.0223559922687482\\
0.015625	0.0216808839015599\\
0.0153846153846154	0.0210132386498092\\
0.0151515151515152	0.0203935269552105\\
0.0149253731343284	0.0197824778732496\\
0.0147058823529412	0.0192171826669636\\
0.0144927536231884	0.0186571682405479\\
0.0142857142857143	0.0181370920851638\\
0.0140845070422535	0.0176283805515531\\
0.0138888888888889	0.0171513091994662\\
0.0136986301369863	0.0166792942507321\\
0.0135135135135135	0.0162390343522691\\
0.0133333333333333	0.0158052413685823\\
0.0131578947368421	0.0153980824036442\\
0.012987012987013	0.0149974323385678\\
0.0128205128205128	0.0146224655642602\\
0.0126582278481013	0.0142514638892806\\
0.0125	0.0139017863460715\\
0.0123456790123457	0.013558595631288\\
0.0121951219512195	0.0132349651328294\\
0.0120481927710843	0.012914601417668\\
0.0119047619047619	0.0126127175908657\\
0.0117647058823529	0.0123161073080809\\
0.0116279069767442	0.0120345619550335\\
0.0114942528735632	0.0117571170705748\\
0.0113636363636364	0.0114941053256317\\
0.0112359550561798	0.0112349812981076\\
0.0111111111111111	0.010988265341533\\
0.010989010989011	0.0107461309102366\\
0.0108695652173913	0.0105168179981892\\
0.010752688172043	0.0102894607851418\\
0.0106382978723404	0.0100731657056965\\
0.0105263157894737	0.0098603591758355\\
0.0104166666666667	0.00965713312564631\\
0.0103092783505155	0.00945693246409185\\
0.0102040816326531	0.00926629011490299\\
0.0101010101010101	0.00907895184768037\\
0.01	0.00889903926044866\\
};
\addlegendentry{PCA}

\addplot [color=color2, mark=*, mark options={solid, color2, scale=0.75}]
  table[row sep=crcr]{%
0.1	0.709661391308028\\
0.0909090909090909	0.528991931567719\\
0.0833333333333333	0.412279892790469\\
0.0769230769230769	0.334676977721522\\
0.0714285714285714	0.282997681309539\\
0.0666666666666667	0.246606451942821\\
0.0625	0.220645368257978\\
0.0588235294117647	0.200722625753297\\
0.0555555555555556	0.185258735216687\\
0.0526315789473684	0.172336225981113\\
0.05	0.161519772826166\\
0.0476190476190476	0.151863901323535\\
0.0454545454545455	0.14331146877873\\
0.0434782608695652	0.135338983655862\\
0.0416666666666667	0.128016925325147\\
0.04	0.121033989913498\\
0.0384615384615385	0.114501877827271\\
0.037037037037037	0.108226153562752\\
0.0357142857142857	0.102322994458461\\
0.0344827586206897	0.0966626401933759\\
0.0333333333333333	0.0913476489390153\\
0.032258064516129	0.0862825878599836\\
0.03125	0.0815489348611251\\
0.0303030303030303	0.0770691420035412\\
0.0294117647058824	0.0729030290765103\\
0.0285714285714286	0.0689833435982603\\
0.0277777777777778	0.0653514770720367\\
0.027027027027027	0.0619478136003497\\
0.0263157894736842	0.0587998829696215\\
0.0256410256410256	0.0558551910624843\\
0.025	0.0531317612369673\\
0.024390243902439	0.0505842878360454\\
0.0238095238095238	0.0482246817954763\\
0.0232558139534884	0.0460147092327308\\
0.0227272727272727	0.043962515852785\\
0.0222222222222222	0.0420364840641039\\
0.0217391304347826	0.0402424290619816\\
0.0212765957446809	0.0385546180655449\\
0.0208333333333333	0.036977344647795\\
0.0204081632653061	0.0354899030168028\\
0.02	0.0340954835186498\\
0.0196078431372549	0.0327775530134247\\
0.0192307692307692	0.0315384113468955\\
0.0188679245283019	0.0303649257446725\\
0.0185185185185185	0.0292586280279519\\
0.0181818181818182	0.028209131010092\\
0.0178571428571429	0.0272172873692149\\
0.0175438596491228	0.0262749258335688\\
0.0172413793103448	0.0253823006952985\\
0.0169491525423729	0.0245330365554643\\
0.0166666666666667	0.0237268699710267\\
0.0163934426229508	0.0229588898528981\\
0.0161290322580645	0.0222283951018216\\
0.0158730158730159	0.0215316772873706\\
0.015625	0.0208676715331411\\
0.0153846153846154	0.0202336622924566\\
0.0151515151515152	0.0196282883422256\\
0.0149253731343284	0.0190496561813926\\
0.0147058823529412	0.0184961667797019\\
0.0144927536231884	0.017966607588237\\
0.0142857142857143	0.0174591915275908\\
0.0140845070422535	0.0169732693278686\\
0.0138888888888889	0.0165069077146498\\
0.0136986301369863	0.0160599214130208\\
0.0135135135135135	0.0156302658248499\\
0.0133333333333333	0.0152181374053377\\
0.0131578947368421	0.014821408553118\\
0.012987012987013	0.014440586818365\\
0.0128205128205128	0.0140734870971144\\
0.0126582278481013	0.0137208683403944\\
0.0125	0.013380508953091\\
0.0123456790123457	0.0130533697484321\\
0.0121951219512195	0.012737210247213\\
0.0120481927710843	0.0124331512361717\\
0.0119047619047619	0.0121389486525736\\
0.0117647058823529	0.01185584822416\\
0.0116279069767442	0.0115816151282253\\
0.0114942528735632	0.0113175905917959\\
0.0113636363636364	0.0110615599476898\\
0.0112359550561798	0.0108149355184359\\
0.0111111111111111	0.0105755316604093\\
0.010989010989011	0.0103448116202607\\
0.0108695652173913	0.0101206253220898\\
0.010752688172043	0.00990447023029972\\
0.0106382978723404	0.00969423897930888\\
0.0105263157894737	0.00949144778348909\\
0.0104166666666667	0.0092940366005454\\
0.0103092783505155	0.00910352946221593\\
0.0102040816326531	0.00891791618078663\\
0.0101010101010101	0.00873871893648648\\
0.01	0.00856398225203314\\
};
\addlegendentry{Com}

\addplot [dashed]
  table[row sep=crcr]{%
0.1	1e-0\\
0.01	1e-2\\
};

\end{axis}
\end{tikzpicture}%

%% file: errors-BSamND-PutOnAvg-set5-phi=0.0413-15D-r=0.04-yDomain-non-uniform-M=10to100.tikz
%
%
%
\begin{tikzpicture}

\begin{axis}[%
width=0.37\textwidth,
scale only axis,
xmode=log,
xmin=0.01,
xmax=0.1,
xminorticks=true,
xlabel={1/m},
ymode=log,
ymin=1e-4,
ymax=1e1,
yminorticks=true,
legend cell align={left},
legend style={nodes={scale=0.75, transform shape}}
]
\addplot [color=color1, mark=square*, mark options={solid, color1, scale=0.75}]
  table[row sep=crcr]{%
0.1	0.427669038187547\\
0.0909090909090909	0.326017209328466\\
0.0833333333333333	0.33503116162867\\
0.0769230769230769	0.310959401635637\\
0.0714285714285714	0.283845459443617\\
0.0666666666666667	0.243809896951287\\
0.0625	0.209817682996741\\
0.0588235294117647	0.171007672874712\\
0.0555555555555556	0.14149377202571\\
0.0526315789473684	0.117465450751935\\
0.05	0.100410476719032\\
0.0476190476190476	0.0892199757272234\\
0.0454545454545455	0.0839034355407264\\
0.0434782608695652	0.0811014161828947\\
0.0416666666666667	0.0817236277275946\\
0.04	0.0787931253765741\\
0.0384615384615385	0.0743143365527743\\
0.037037037037037	0.0666253322225607\\
0.0357142857142857	0.0593537966176072\\
0.0344827586206897	0.0519218946683413\\
0.0333333333333333	0.0459550027890501\\
0.032258064516129	0.0401169222667921\\
0.03125	0.0366063111253752\\
0.0303030303030303	0.0341923533873754\\
0.0294117647058824	0.0332330818763804\\
0.0285714285714286	0.0330341865670523\\
0.0277777777777778	0.0335707774310077\\
0.027027027027027	0.0336055012521932\\
0.0263157894736842	0.0332372014361209\\
0.0256410256410256	0.0317288352842855\\
0.025	0.029220376139925\\
0.024390243902439	0.0266923433773325\\
0.0238095238095238	0.0241931236505089\\
0.0232558139534884	0.0218694233849419\\
0.0227272727272727	0.0202622828921815\\
0.0222222222222222	0.0188998721679461\\
0.0217391304347826	0.0180401101067496\\
0.0212765957446809	0.017578358904256\\
0.0208333333333333	0.0175204962428381\\
0.0204081632653061	0.0178025849534205\\
0.02	0.0179283455421242\\
0.0196078431372549	0.0176779707256784\\
0.0192307692307692	0.0173358430466184\\
0.0188679245283019	0.01646632082271\\
0.0185185185185185	0.0153833769407896\\
0.0181818181818182	0.0143905277492955\\
0.0178571428571429	0.0132612972301391\\
0.0175438596491228	0.0124451939257781\\
0.0172413793103448	0.011746559786026\\
0.0169491525423729	0.011170044073447\\
0.0166666666666667	0.0109511452110951\\
0.0163934426229508	0.0108482704391584\\
0.0161290322580645	0.0108519940997294\\
0.0158730158730159	0.0109358082185516\\
0.015625	0.0108862002070484\\
0.0153846153846154	0.0108462744059619\\
0.0151515151515152	0.010589689618687\\
0.0149253731343284	0.0101562054710413\\
0.0147058823529412	0.00952046644678362\\
0.0144927536231884	0.00894409509392169\\
0.0142857142857143	0.00849444704764224\\
0.0140845070422535	0.00812693086064264\\
0.0138888888888889	0.00765946737277368\\
0.0136986301369863	0.00745822399167473\\
0.0135135135135135	0.00731582778522366\\
0.0133333333333333	0.00730315061173137\\
0.0131578947368421	0.00727622010534823\\
0.012987012987013	0.00730088601952517\\
0.0128205128205128	0.0073762227485048\\
0.0126582278481013	0.00718187562624273\\
0.0125	0.00706042591844103\\
0.0123456790123457	0.00677638736130248\\
0.0121951219512195	0.00650048656127034\\
0.0120481927710843	0.00618018693355893\\
0.0119047619047619	0.00590896623175396\\
0.0117647058823529	0.00565319096440042\\
0.0116279069767442	0.00550855916719151\\
0.0114942528735632	0.00532867108128343\\
0.0113636363636364	0.00532551799917513\\
0.0112359550561798	0.005242831104038\\
0.0111111111111111	0.00522756669481161\\
0.010989010989011	0.00523895089332171\\
0.0108695652173913	0.00519643080354903\\
0.010752688172043	0.00518994782158977\\
0.0106382978723404	0.0051167407613395\\
0.0105263157894737	0.00489514447431816\\
0.0104166666666667	0.00473436446208808\\
0.0103092783505155	0.0045779642881616\\
0.0102040816326531	0.00438974854339946\\
0.0101010101010101	0.00425922520628674\\
0.01	0.00410313697667952\\
};
\addlegendentry{PCA (IT)}

\addplot [color=color2, mark=*, mark options={solid, color2, scale=0.75}]
  table[row sep=crcr]{%
0.1	0.384578102170504\\
0.0909090909090909	0.311536303244043\\
0.0833333333333333	0.308132218278329\\
0.0769230769230769	0.288283959195887\\
0.0714285714285714	0.258969530919352\\
0.0666666666666667	0.222513585051525\\
0.0625	0.186837301720338\\
0.0588235294117647	0.154078976430649\\
0.0555555555555556	0.12691194525641\\
0.0526315789473684	0.107295969391634\\
0.05	0.0927917703201011\\
0.0476190476190476	0.0847744181308236\\
0.0454545454545455	0.0808852393027877\\
0.0434782608695652	0.0793743531274931\\
0.0416666666666667	0.0776965677515342\\
0.04	0.0727185128745376\\
0.0384615384615385	0.0665027273807954\\
0.037037037037037	0.0587721934181218\\
0.0357142857142857	0.0524371088756972\\
0.0344827586206897	0.0455734437379363\\
0.0333333333333333	0.0403876897204274\\
0.032258064516129	0.036649874538492\\
0.03125	0.0342587459144461\\
0.0303030303030303	0.033166781275785\\
0.0294117647058824	0.0328516488953792\\
0.0285714285714286	0.0330954976398488\\
0.0277777777777778	0.033005429893322\\
0.027027027027027	0.0320330180455688\\
0.0263157894736842	0.0302148348732008\\
0.0256410256410256	0.0278676806294555\\
0.025	0.0255698530628621\\
0.024390243902439	0.0234033646324399\\
0.0238095238095238	0.0211838672834563\\
0.0232558139534884	0.01960743482854\\
0.0227272727272727	0.0184629220013424\\
0.0222222222222222	0.0178288890928453\\
0.0217391304347826	0.0176112988654165\\
0.0212765957446809	0.0173711403904644\\
0.0208333333333333	0.0175294757984954\\
0.0204081632653061	0.017483950680516\\
0.02	0.0170018345046259\\
0.0196078431372549	0.0162377642639164\\
0.0192307692307692	0.0152995287566355\\
0.0188679245283019	0.0143451310141405\\
0.0185185185185185	0.0134286002637278\\
0.0181818181818182	0.0124124480501813\\
0.0178571428571429	0.0117567680179231\\
0.0175438596491228	0.0112910492681308\\
0.0172413793103448	0.0110267056797344\\
0.0169491525423729	0.0108424871596373\\
0.0166666666666667	0.0107400616998703\\
0.0163934426229508	0.0108085368252766\\
0.0161290322580645	0.0107318088609905\\
0.0158730158730159	0.0104908438705564\\
0.015625	0.0101495748663432\\
0.0153846153846154	0.00966844875025963\\
0.0151515151515152	0.00922723396231007\\
0.0149253731343284	0.00875299380850336\\
0.0147058823529412	0.0082599967876682\\
0.0144927536231884	0.00792653197203563\\
0.0142857142857143	0.00765471238339721\\
0.0140845070422535	0.00748863389730658\\
0.0138888888888889	0.00738146497982584\\
0.0136986301369863	0.00728030304022198\\
0.0135135135135135	0.00724386321282289\\
0.0133333333333333	0.00723956226672051\\
0.0131578947368421	0.00711771690924756\\
0.012987012987013	0.00693743117899892\\
0.0128205128205128	0.00668852081076388\\
0.0126582278481013	0.0064657705734863\\
0.0125	0.00620343254265521\\
0.0123456790123457	0.00591943348451984\\
0.0121951219512195	0.0057018833245035\\
0.0120481927710843	0.00553726709396152\\
0.0119047619047619	0.00542388644317704\\
0.0117647058823529	0.00532444895521067\\
0.0116279069767442	0.00526294318457232\\
0.0114942528735632	0.00523020185073442\\
0.0113636363636364	0.00521148482203948\\
0.0112359550561798	0.00513429073406613\\
0.0111111111111111	0.00506644234881648\\
0.010989010989011	0.004892857835463\\
0.0108695652173913	0.00477359353597961\\
0.010752688172043	0.0046163754699684\\
0.0106382978723404	0.00445528359622394\\
0.0105263157894737	0.0043102400935453\\
0.0104166666666667	0.00418874116548018\\
0.0103092783505155	0.00411296137353112\\
0.0102040816326531	0.00406374351273797\\
0.0101010101010101	0.00398150491032467\\
0.01	0.00394267232523404\\
};
\addlegendentry{Com (IT)}

\addplot [dashed]
  table[row sep=crcr]{%
0.1	1e-0\\
0.01	1e-2\\
};

\end{axis}
\end{tikzpicture}%

%% file: BlackScholesNDEuropean-K=40-T=0.5-vol=0.3-rho=0.8-r=0.05.tikz
%
%
%
\begin{tikzpicture}

\begin{axis}[%
width=0.37\textwidth,
scale only axis,
xmode=log,
xmin=0.01,
xmax=0.1,
xminorticks=true,
ymode=log,
ymin=1e-3,
ymax=1e2,
yminorticks=true,
ylabel={error},
title={European basket option},
legend cell align={left},
legend style={nodes={scale=0.75, transform shape}}
]
\addplot [color=color1, mark=square*, mark options={solid, color1, scale=0.75}]
  table[row sep=crcr]{%
0.1	0.447475078927405\\
0.0909090909090909	0.335680031854033\\
0.0833333333333333	0.266385838235297\\
0.0769230769230769	0.216209485350934\\
0.0714285714285714	0.179947333966255\\
0.0666666666666667	0.154277808023648\\
0.0625	0.133487145107259\\
0.0588235294117647	0.118916160329958\\
0.0555555555555556	0.106478730154304\\
0.0526315789473684	0.0964204533833222\\
0.05	0.088094389337277\\
0.0476190476190476	0.0805724049146042\\
0.0454545454545455	0.074514680676776\\
0.0434782608695652	0.0689713573919244\\
0.0416666666666667	0.0640367435790115\\
0.04	0.0598987380323557\\
0.0384615384615385	0.055871691221542\\
0.037037037037037	0.0524642653396761\\
0.0357142857142857	0.0490944942819631\\
0.0344827586206897	0.0461102659844608\\
0.0333333333333333	0.0433948331194802\\
0.032258064516129	0.0407386685017057\\
0.03125	0.0384721761849338\\
0.0303030303030303	0.0361712077907681\\
0.0294117647058824	0.0341583901308926\\
0.0285714285714286	0.03227952560378\\
0.0277777777777778	0.0304981252763064\\
0.027027027027027	0.0289594949223941\\
0.0263157894736842	0.0273812288308681\\
0.0256410256410256	0.026021157665423\\
0.025	0.0247330559387677\\
0.024390243902439	0.0235267433501001\\
0.0238095238095238	0.0224741416559571\\
0.0232558139534884	0.0213937205539647\\
0.0227272727272727	0.0204649206161003\\
0.0222222222222222	0.0195609172808515\\
0.0217391304347826	0.0187129560118189\\
0.0212765957446809	0.0179644168695985\\
0.0208333333333333	0.0172014574440231\\
0.0204081632653061	0.0165383202635843\\
0.02	0.0158730926273876\\
0.0196078431372549	0.0152623506016631\\
0.0192307692307692	0.0147021660090223\\
0.0188679245283019	0.0141401264989813\\
0.0185185185185185	0.0136463775343669\\
0.0181818181818182	0.013142435472461\\
0.0178571428571429	0.0126855097317584\\
0.0175438596491228	0.0122552882190625\\
0.0172413793103448	0.0118274087530086\\
0.0169491525423729	0.0114496517769886\\
0.0166666666666667	0.011057927592967\\
0.0163934426229508	0.0107058093260211\\
0.0161290322580645	0.0103686748616569\\
0.0158730158730159	0.0100354695926042\\
0.015625	0.00974087628244114\\
0.0153846153846154	0.0094291920720142\\
0.0151515151515152	0.00915432215675338\\
0.0149253731343284	0.00888280496112159\\
0.0147058823529412	0.0086211295439762\\
0.0144927536231884	0.00838596812498071\\
0.0142857142857143	0.0081369522831265\\
0.0140845070422535	0.00791618948168704\\
0.0138888888888889	0.007695962601872\\
0.0136986301369863	0.00748539156787054\\
0.0135135135135135	0.00729269490180862\\
0.0133333333333333	0.00709221807783944\\
0.0131578947368421	0.00691262413289184\\
0.012987012987013	0.00673095053225303\\
0.0128205128205128	0.00655997931522556\\
0.0126582278481013	0.00639914998386271\\
0.0125	0.00623496483150721\\
0.0123456790123457	0.00608807812358592\\
0.0121951219512195	0.00593605495622818\\
0.0120481927710843	0.00579461215532362\\
0.0119047619047619	0.00566069873920316\\
0.0117647058823529	0.00552439984519726\\
0.0116279069767442	0.00540252165377719\\
0.0114942528735632	0.00527367430951209\\
0.0113636363636364	0.00515607603078205\\
0.0112359550561798	0.00504164604552315\\
0.0111111111111111	0.00492797147569046\\
0.010989010989011	0.00482554049603046\\
0.0108695652173913	0.00471570196401672\\
0.010752688172043	0.00461670665225533\\
0.0106382978723404	0.00451880399328797\\
0.0105263157894737	0.00442243561888844\\
0.0104166666666667	0.00433522102079742\\
0.0103092783505155	0.0042414284182124\\
0.0102040816326531	0.00415711073720093\\
0.0101010101010101	0.00407228716067642\\
0.01	0.00399044158508222\\
};
\addlegendentry{PCA}

\addplot [color=color2, mark=*, mark options={solid, color2, scale=0.75}]
  table[row sep=crcr]{%
0.1	0.399945498912554\\
0.0909090909090909	0.304201852405824\\
0.0833333333333333	0.239007038788342\\
0.0769230769230769	0.194271903108161\\
0.0714285714285714	0.162584920145924\\
0.0666666666666667	0.139279412137959\\
0.0625	0.121832488752112\\
0.0588235294117647	0.108054794966998\\
0.0555555555555556	0.0970408391577919\\
0.0526315789473684	0.0879604014128912\\
0.05	0.0803079615862732\\
0.0476190476190476	0.0738598939882422\\
0.0454545454545455	0.068210954080536\\
0.0434782608695652	0.0633411754569817\\
0.0416666666666667	0.0589758523003434\\
0.04	0.0550461747580808\\
0.0384615384615385	0.0515002511449945\\
0.037037037037037	0.0482134667636993\\
0.0357142857142857	0.0452429747410203\\
0.0344827586206897	0.0424602997557662\\
0.0333333333333333	0.0399119586039527\\
0.032258064516129	0.0375424847129686\\
0.03125	0.0353396550576344\\
0.0303030303030303	0.0333163056201871\\
0.0294117647058824	0.0314229005799866\\
0.0285714285714286	0.0296931328260772\\
0.0277777777777778	0.0280817615887319\\
0.027027027027027	0.0265927942000737\\
0.0263157894736842	0.0252215201531172\\
0.0256410256410256	0.0239397406555364\\
0.025	0.022769724048163\\
0.024390243902439	0.0216734482061334\\
0.0238095238095238	0.0206606727932943\\
0.0232558139534884	0.0197185721528408\\
0.0227272727272727	0.0188357013686549\\
0.0222222222222222	0.0180203719529795\\
0.0217391304347826	0.0172490462563344\\
0.0212765957446809	0.0165346328417444\\
0.0208333333333333	0.0158603707979861\\
0.0204081632653061	0.0152262671622321\\
0.02	0.0146325578371611\\
0.0196078431372549	0.0140669243138918\\
0.0192307692307692	0.0135410224807622\\
0.0188679245283019	0.0130390141288359\\
0.0185185185185185	0.012565763550672\\
0.0181818181818182	0.0121182467868621\\
0.0178571428571429	0.0116909226497341\\
0.0175438596491228	0.0112901650130679\\
0.0172413793103448	0.0109055277408485\\
0.0169491525423729	0.0105430317330057\\
0.0166666666666667	0.0101970341828999\\
0.0163934426229508	0.00986674647849561\\
0.0161290322580645	0.00955448139430004\\
0.0158730158730159	0.00925347107128172\\
0.015625	0.00897006487182228\\
0.0153846153846154	0.00869729781609063\\
0.0151515151515152	0.00843716457349952\\
0.0149253731343284	0.00818918488717468\\
0.0147058823529412	0.007949868876584\\
0.0144927536231884	0.00772380240304837\\
0.0142857142857143	0.0075049957587332\\
0.0140845070422535	0.00729652090235966\\
0.0138888888888889	0.00709634334608644\\
0.0136986301369863	0.00690333190982351\\
0.0135135135135135	0.00671970543007827\\
0.0133333333333333	0.00654139405591003\\
0.0131578947368421	0.00637184866821361\\
0.012987012987013	0.00620781789273117\\
0.0128205128205128	0.00604991508742447\\
0.0126582278481013	0.00589866949801543\\
0.0125	0.00575143646707144\\
0.0123456790123457	0.00561167492171766\\
0.0121951219512195	0.00547554499692549\\
0.0120481927710843	0.00534478174517261\\
0.0119047619047619	0.00521869401840824\\
0.0117647058823529	0.00509616407446867\\
0.0116279069767442	0.00497915448734254\\
0.0114942528735632	0.00486492772989156\\
0.0113636363636364	0.00475544741034373\\
0.0112359550561798	0.00464920922193635\\
0.0111111111111111	0.00454620920667637\\
0.010989010989011	0.00444724620275583\\
0.0108695652173913	0.00435042859604273\\
0.010752688172043	0.00425789731039927\\
0.0106382978723404	0.0041675407646018\\
0.0105263157894737	0.00408015123144523\\
0.0104166666666667	0.00399568188721666\\
0.0103092783505155	0.00391308323233464\\
0.0102040816326531	0.00383401782722892\\
0.0101010101010101	0.00375650356975132\\
0.01	0.00368174748226835\\
};
\addlegendentry{Com}

\addplot [dashed]
  table[row sep=crcr]{%
0.1	1e0\\
0.01	1e-2\\
};

\end{axis}
\end{tikzpicture}%

%% file: BlackScholesNDAmerican-K=40-T=0.5-vol=0.3-rho=0.8-r=0.05.tikz
%
%
%
\begin{tikzpicture}

\begin{axis}[%
width=0.37\textwidth,
scale only axis,
xmode=log,
xmin=0.01,
xmax=0.1,
xminorticks=true,
ymode=log,
ymin=1e-3,
ymax=1e2,
yminorticks=true,
title={American basket option},
legend cell align={left},
legend style={nodes={scale=0.75, transform shape}}
]
\addplot [color=color1, mark=square*, mark options={solid, color1, scale=0.75}]
  table[row sep=crcr]{%
0.1	0.365861172294772\\
0.0909090909090909	0.249949524899886\\
0.0833333333333333	0.190676357117679\\
0.0769230769230769	0.163663145566604\\
0.0714285714285714	0.151783913461983\\
0.0666666666666667	0.145161045304559\\
0.0625	0.132910363132097\\
0.0588235294117647	0.118852473903275\\
0.0555555555555556	0.104695423443437\\
0.0526315789473684	0.0900514713219858\\
0.05	0.0769814818714503\\
0.0476190476190476	0.0661900564294555\\
0.0454545454545455	0.0572864982646086\\
0.0434782608695652	0.0507413116229385\\
0.0416666666666667	0.0441443546364617\\
0.04	0.0400926953463117\\
0.0384615384615385	0.037336573821146\\
0.037037037037037	0.0352201914457053\\
0.0357142857142857	0.033997109461775\\
0.0344827586206897	0.0332407878312377\\
0.0333333333333333	0.0316543653025514\\
0.032258064516129	0.0308667204996373\\
0.03125	0.0295420779204907\\
0.0303030303030303	0.0275259979572766\\
0.0294117647058824	0.0256987039257552\\
0.0285714285714286	0.0230277081362011\\
0.0277777777777778	0.0213796697792903\\
0.027027027027027	0.0194960599063805\\
0.0263157894736842	0.0176318669020494\\
0.0256410256410256	0.0162590858769773\\
0.025	0.0153521615355787\\
0.024390243902439	0.0144394220388184\\
0.0238095238095238	0.0138276992611326\\
0.0232558139534884	0.0134029521955235\\
0.0227272727272727	0.0131125159710308\\
0.0222222222222222	0.0129225057158084\\
0.0217391304347826	0.0126708347624289\\
0.0212765957446809	0.0123512250546698\\
0.0208333333333333	0.0121188717801806\\
0.0204081632653061	0.0117237649500703\\
0.02	0.0110289404087593\\
0.0196078431372549	0.0104709624437778\\
0.0192307692307692	0.00995071031935968\\
0.0188679245283019	0.00938297770764507\\
0.0185185185185185	0.00880821342745719\\
0.0181818181818182	0.0081817233223731\\
0.0178571428571429	0.00772709012053152\\
0.0175438596491228	0.00746334223261336\\
0.0172413793103448	0.0072885311311337\\
0.0169491525423729	0.00695708507898818\\
0.0166666666666667	0.00676355182339972\\
0.0163934426229508	0.00667102891125104\\
0.0161290322580645	0.00655435844943408\\
0.0158730158730159	0.00655118310901059\\
0.015625	0.00646720450664251\\
0.0153846153846154	0.00622840363232324\\
0.0151515151515152	0.00605400229672437\\
0.0149253731343284	0.00598606215707242\\
0.0147058823529412	0.00575345662912152\\
0.0144927536231884	0.00549960360047841\\
0.0142857142857143	0.00524939203807673\\
0.0140845070422535	0.00497047479213819\\
0.0138888888888889	0.00484944431966561\\
0.0136986301369863	0.00468024023571445\\
0.0135135135135135	0.00440838325047022\\
0.0133333333333333	0.00422749322299243\\
0.0131578947368421	0.00423940386775712\\
0.012987012987013	0.00413210322896163\\
0.0128205128205128	0.00408507684376769\\
0.0126582278481013	0.00395847014730943\\
0.0125	0.00384628664360864\\
0.0123456790123457	0.00388075268815591\\
0.0121951219512195	0.00382461138825896\\
0.0120481927710843	0.00379248866341175\\
0.0119047619047619	0.00367505778894994\\
0.0117647058823529	0.00359695387707237\\
0.0116279069767442	0.00343924965023934\\
0.0114942528735632	0.00337183154177723\\
0.0113636363636364	0.00331120596869638\\
0.0112359550561798	0.0031593830847978\\
0.0111111111111111	0.00303585908678272\\
0.010989010989011	0.0029354157863688\\
0.0108695652173913	0.00286812928397406\\
0.010752688172043	0.00281437214341462\\
0.0106382978723404	0.00279999017820742\\
0.0105263157894737	0.00266279865518193\\
0.0104166666666667	0.00264689028541021\\
0.0103092783505155	0.00262514876325337\\
0.0102040816326531	0.00260565357702536\\
0.0101010101010101	0.00258141144357804\\
0.01	0.00249911453939689\\
};
\addlegendentry{PCA (IT)}

\addplot [color=color2, mark=*, mark options={solid, color2, scale=0.75}]
  table[row sep=crcr]{%
0.1	0.315733012721794\\
0.0909090909090909	0.217766729159695\\
0.0833333333333333	0.167924989852245\\
0.0769230769230769	0.146483520427633\\
0.0714285714285714	0.138038414658555\\
0.0666666666666667	0.129337522286752\\
0.0625	0.116941456440051\\
0.0588235294117647	0.103461819207115\\
0.0555555555555556	0.0898386881629696\\
0.0526315789473684	0.0767459966885538\\
0.05	0.065498624892915\\
0.0476190476190476	0.0559701883811972\\
0.0454545454545455	0.0489928027810631\\
0.0434782608695652	0.0433194129519192\\
0.0416666666666667	0.0391362615267985\\
0.04	0.0360967519172162\\
0.0384615384615385	0.0341235052823547\\
0.037037037037037	0.0327576778652254\\
0.0357142857142857	0.0317386898773639\\
0.0344827586206897	0.0304977846507102\\
0.0333333333333333	0.0292692304881106\\
0.032258064516129	0.0274677863185895\\
0.03125	0.0256154302433815\\
0.0303030303030303	0.0234999724260305\\
0.0294117647058824	0.0213623974756478\\
0.0285714285714286	0.0194652767013404\\
0.0277777777777778	0.0178479814085817\\
0.027027027027027	0.0163860222074517\\
0.0263157894736842	0.0152555027468209\\
0.0256410256410256	0.0143817275423279\\
0.025	0.0137181979579664\\
0.024390243902439	0.0132530407194675\\
0.0238095238095238	0.012902931691535\\
0.0232558139534884	0.0126132926291112\\
0.0227272727272727	0.0123502958298989\\
0.0222222222222222	0.0119753565215319\\
0.0217391304347826	0.0116035071972265\\
0.0212765957446809	0.0111791562465768\\
0.0208333333333333	0.0105664729627231\\
0.0204081632653061	0.0100248758856782\\
0.02	0.00939880551199934\\
0.0196078431372549	0.00887871803320195\\
0.0192307692307692	0.0083677152724686\\
0.0188679245283019	0.00792742673600166\\
0.0185185185185185	0.00749072768565018\\
0.0181818181818182	0.00723457002679684\\
0.0178571428571429	0.00702450782622144\\
0.0175438596491228	0.00680903216426465\\
0.0172413793103448	0.00665028234151954\\
0.0169491525423729	0.0065058285259445\\
0.0166666666666667	0.00640486235661619\\
0.0163934426229508	0.00629461204371573\\
0.0161290322580645	0.00613251557010397\\
0.0158730158730159	0.00592934511160248\\
0.015625	0.00572286238727138\\
0.0153846153846154	0.00552695622186938\\
0.0151515151515152	0.00531167716181269\\
0.0149253731343284	0.0050686480690123\\
0.0147058823529412	0.00484100724345637\\
0.0144927536231884	0.00465126761739043\\
0.0142857142857143	0.00451627376810038\\
0.0140845070422535	0.00434037733997794\\
0.0138888888888889	0.00423118336781236\\
0.0136986301369863	0.00411380307825304\\
0.0135135135135135	0.00403286422216276\\
0.0133333333333333	0.00398532592317391\\
0.0131578947368421	0.00390909211815504\\
0.012987012987013	0.0038227326459026\\
0.0128205128205128	0.00375823791365626\\
0.0126582278481013	0.00368475357272136\\
0.0125	0.00358719267455143\\
0.0123456790123457	0.00350056272265498\\
0.0121951219512195	0.00339277803407079\\
0.0120481927710843	0.00327100913065159\\
0.0119047619047619	0.00317862196797503\\
0.0117647058823529	0.00308850849556652\\
0.0116279069767442	0.00298037266190221\\
0.0114942528735632	0.00291041018288407\\
0.0113636363636364	0.00282258515617695\\
0.0112359550561798	0.00276525940982797\\
0.0111111111111111	0.00271882782575883\\
0.010989010989011	0.00267259285685117\\
0.0108695652173913	0.00261682656601625\\
0.010752688172043	0.00257087140994194\\
0.0106382978723404	0.00254525799623906\\
0.0105263157894737	0.00250101348294951\\
0.0104166666666667	0.00245899392434445\\
0.0103092783505155	0.0023933802339382\\
0.0102040816326531	0.00234198070190139\\
0.0101010101010101	0.00229216681961386\\
0.01	0.00224067184206866\\
};
\addlegendentry{Com (IT)}

\addplot [dashed]
  table[row sep=crcr]{%
0.1	1e0\\
0.01	1e-2\\
};

\end{axis}
\end{tikzpicture}%

%% file: BlackScholesNDEuropean-K=40-T=1-vol=0.3-rho=0.8-r=0.05.tikz
%
%
%
\begin{tikzpicture}

\begin{axis}[%
width=0.37\textwidth,
scale only axis,
xmode=log,
xmin=0.01,
xmax=0.1,
xminorticks=true,
ymode=log,
ymin=1e-3,
ymax=1e2,
yminorticks=true,
ylabel={error},
legend cell align={left},
legend style={nodes={scale=0.75, transform shape}}
]
\addplot [color=color1, mark=square*, mark options={solid, color1, scale=0.75}]
  table[row sep=crcr]{%
0.1	1.27521590336948\\
0.0909090909090909	0.904470669230356\\
0.0833333333333333	0.651010443061034\\
0.0769230769230769	0.476788039427476\\
0.0714285714285714	0.361364072878092\\
0.0666666666666667	0.289672365768221\\
0.0625	0.243547496064864\\
0.0588235294117647	0.214108769619894\\
0.0555555555555556	0.196112724490177\\
0.0526315789473684	0.182873989107367\\
0.05	0.172237796716836\\
0.0476190476190476	0.164188557501024\\
0.0454545454545455	0.156050230266603\\
0.0434782608695652	0.148223784219104\\
0.0416666666666667	0.141202870661489\\
0.04	0.133638435386787\\
0.0384615384615385	0.126168696117011\\
0.037037037037037	0.119390285111178\\
0.0357142857142857	0.112518628487241\\
0.0344827586206897	0.105831015195841\\
0.0333333333333333	0.0997938499260673\\
0.032258064516129	0.0938334210056642\\
0.03125	0.0882864665044352\\
0.0303030303030303	0.0834265830071912\\
0.0294117647058824	0.0787136058433804\\
0.0285714285714286	0.074449079741048\\
0.0277777777777778	0.0707136744113033\\
0.027027027027027	0.0670915435543264\\
0.0263157894736842	0.063797375484449\\
0.0256410256410256	0.0609167649878302\\
0.025	0.0580669716776905\\
0.024390243902439	0.0554284598533119\\
0.0238095238095238	0.0530916665126213\\
0.0232558139534884	0.0507566539185813\\
0.0227272727272727	0.0485954789204097\\
0.0222222222222222	0.0466598340611082\\
0.0217391304347826	0.0447151191500206\\
0.0212765957446809	0.0429017825878137\\
0.0208333333333333	0.0412492180506963\\
0.0204081632653061	0.0395874759744848\\
0.02	0.038054814121633\\
0.0196078431372549	0.0366438753758747\\
0.0192307692307692	0.0352305914999773\\
0.0188679245283019	0.0339116084492401\\
0.0185185185185185	0.0327114440525138\\
0.0181818181818182	0.0314940811458229\\
0.0178571428571429	0.0303737611255057\\
0.0175438596491228	0.0293494747467093\\
0.0172413793103448	0.0283144733305223\\
0.0169491525423729	0.0273530714380952\\
0.0166666666666667	0.0264796176606046\\
0.0163934426229508	0.0255950034186512\\
0.0161290322580645	0.0247737788420315\\
0.0158730158730159	0.0240227359760778\\
0.015625	0.0232558349119989\\
0.0153846153846154	0.0225476613597904\\
0.0151515151515152	0.0218953271357236\\
0.0149253731343284	0.0212307531368632\\
0.0147058823529412	0.0206152963791393\\
0.0144927536231884	0.0200501829978901\\
0.0142857142857143	0.0194678941370139\\
0.0140845070422535	0.0189293940367614\\
0.0138888888888889	0.0184312472900245\\
0.0136986301369863	0.0179190526655759\\
0.0135135135135135	0.0174436976581944\\
0.0133333333333333	0.0170012249683298\\
0.0131578947368421	0.016546662848973\\
0.012987012987013	0.0161233558718079\\
0.0128205128205128	0.0157291331511047\\
0.0126582278481013	0.0153262378450725\\
0.0125	0.0149475968582893\\
0.0123456790123457	0.0145972228532454\\
0.0121951219512195	0.0142316152146469\\
0.0120481927710843	0.013891679501767\\
0.0119047619047619	0.0135769687098612\\
0.0117647058823529	0.0132488471816572\\
0.0116279069767442	0.0129420258892923\\
0.0114942528735632	0.0126587077300631\\
0.0113636363636364	0.0123620410220084\\
0.0112359550561798	0.012085197988009\\
0.0111111111111111	0.0118290382057529\\
0.010989010989011	0.0115599758795817\\
0.0108695652173913	0.0113101736381305\\
0.010752688172043	0.0110776411758975\\
0.0106382978723404	0.0108324025218058\\
0.0105263157894737	0.0106055680843751\\
0.0104166666666667	0.01039475425614\\
0.0103092783505155	0.0101722922337526\\
0.0102040816326531	0.00996563798258521\\
0.0101010101010101	0.00977282095386123\\
0.01	0.00956928278468361\\
};
\addlegendentry{PCA}

\addplot [color=color2, mark=*, mark options={solid, color2, scale=0.75}]
  table[row sep=crcr]{%
0.1	1.11944778133697\\
0.0909090909090909	0.805205166344805\\
0.0833333333333333	0.585637559856972\\
0.0769230769230769	0.437552186780856\\
0.0714285714285714	0.339588271269528\\
0.0666666666666667	0.275092598660522\\
0.0625	0.233240495271916\\
0.0588235294117647	0.205597797565158\\
0.0555555555555556	0.186607750814877\\
0.0526315789473684	0.173071645291593\\
0.05	0.162498036416546\\
0.0476190476190476	0.153465409830739\\
0.0454545454545455	0.145510853290755\\
0.0434782608695652	0.137911420265853\\
0.0416666666666667	0.130517937431947\\
0.04	0.123385111867566\\
0.0384615384615385	0.11635955690333\\
0.037037037037037	0.109565676057671\\
0.0357142857142857	0.103148413187578\\
0.0344827586206897	0.0970565331971995\\
0.0333333333333333	0.0913714974150297\\
0.032258064516129	0.0861356708956977\\
0.03125	0.0812731964899323\\
0.0303030303030303	0.0768148417928094\\
0.0294117647058824	0.0727281684404772\\
0.0285714285714286	0.0689354023712392\\
0.0277777777777778	0.0654491385841167\\
0.027027027027027	0.0622315779765259\\
0.0263157894736842	0.05923069673974\\
0.0256410256410256	0.0564626876562322\\
0.025	0.0538806138209802\\
0.024390243902439	0.0514503600527458\\
0.0238095238095238	0.049193299901594\\
0.0232558139534884	0.0470682294468903\\
0.0227272727272727	0.0450548573556562\\
0.0222222222222222	0.0431747590768747\\
0.0217391304347826	0.041390167628828\\
0.0212765957446809	0.0396933906433992\\
0.0208333333333333	0.0381068120683565\\
0.0204081632653061	0.0365972889331028\\
0.02	0.0351620613051864\\
0.0196078431372549	0.0338155489636289\\
0.0192307692307692	0.0325334591291595\\
0.0188679245283019	0.0313169325704799\\
0.0185185185185185	0.030174077238418\\
0.0181818181818182	0.0290849042318442\\
0.0178571428571429	0.0280518336217366\\
0.0175438596491228	0.0270787057545618\\
0.0172413793103448	0.0261509221506948\\
0.0169491525423729	0.0252718323679888\\
0.0166666666666667	0.0244407087626781\\
0.0163934426229508	0.023646192845229\\
0.0161290322580645	0.0228929670364142\\
0.0158730158730159	0.0221783895258731\\
0.015625	0.0214939241874132\\
0.0153846153846154	0.0208442414354568\\
0.0151515151515152	0.0202246846647469\\
0.0149253731343284	0.0196294870860338\\
0.0147058823529412	0.019064288315886\\
0.0144927536231884	0.0185231776640347\\
0.0142857142857143	0.018001743332106\\
0.0140845070422535	0.0175058568363262\\
0.0138888888888889	0.0170289968636581\\
0.0136986301369863	0.0165686396951958\\
0.0135135135135135	0.0161308892132341\\
0.0133333333333333	0.0157081764266804\\
0.0131578947368421	0.0152991136077656\\
0.012987012987013	0.0149098680092896\\
0.0128205128205128	0.0145329641303098\\
0.0126582278481013	0.0141684213675219\\
0.0125	0.013820473473336\\
0.0123456790123457	0.0134828157775111\\
0.0121951219512195	0.0131563014682845\\
0.0120481927710843	0.0128438793031123\\
0.0119047619047619	0.0125405134349208\\
0.0117647058823529	0.0122472761877575\\
0.0116279069767442	0.0119657517110072\\
0.0114942528735632	0.0116920442294761\\
0.0113636363636364	0.0114277583102398\\
0.0112359550561798	0.0111734829353969\\
0.0111111111111111	0.0109259892476263\\
0.010989010989011	0.0106870582102445\\
0.0108695652173913	0.0104565169726119\\
0.010752688172043	0.0102319800718558\\
0.0106382978723404	0.0100154868331037\\
0.0105263157894737	0.00980601686587157\\
0.0104166666666667	0.00960170215543066\\
0.0103092783505155	0.0094048215887792\\
0.0102040816326531	0.00921387716155042\\
0.0101010101010101	0.00902754523977745\\
0.01	0.00884812105256838\\
};
\addlegendentry{Com}

\addplot [dashed]
  table[row sep=crcr]{%
0.1	1e1\\
0.01	1e-1\\
};

\end{axis}
\end{tikzpicture}%

%% file: BlackScholesNDAmerican-K=40-T=1-vol=0.3-rho=0.8-r=0.05.tikz
%
%
%
\begin{tikzpicture}

\begin{axis}[%
width=0.37\textwidth,
scale only axis,
xmode=log,
xmin=0.01,
xmax=0.1,
xminorticks=true,
ymode=log,
ymin=1e-3,
ymax=1e2,
yminorticks=true,
legend cell align={left},
legend style={nodes={scale=0.75, transform shape}}
]
\addplot [color=color1, mark=square*, mark options={solid, color1, scale=0.75}]
  table[row sep=crcr]{%
0.1	1.1485513969066\\
0.0909090909090909	0.754429832642291\\
0.0833333333333333	0.485725811776904\\
0.0769230769230769	0.330419137702211\\
0.0714285714285714	0.24070643966861\\
0.0666666666666667	0.209500317433946\\
0.0625	0.207839542528788\\
0.0588235294117647	0.223698694287203\\
0.0555555555555556	0.247141185915879\\
0.0526315789473684	0.2621918692239\\
0.05	0.263635274626243\\
0.0476190476190476	0.255312038758226\\
0.0454545454545455	0.237324619198085\\
0.0434782608695652	0.218407962178539\\
0.0416666666666667	0.195537647810552\\
0.04	0.17500910822542\\
0.0384615384615385	0.152472089074566\\
0.037037037037037	0.132625876391952\\
0.0357142857142857	0.116407368762463\\
0.0344827586206897	0.10027012587443\\
0.0333333333333333	0.0870081341943889\\
0.032258064516129	0.076929880004065\\
0.03125	0.0676839887133767\\
0.0303030303030303	0.0622655077702365\\
0.0294117647058824	0.0573136551348785\\
0.0285714285714286	0.0548699859616999\\
0.0277777777777778	0.0533587524628931\\
0.027027027027027	0.0526002342622833\\
0.0263157894736842	0.0533949437291907\\
0.0256410256410256	0.0546056993359318\\
0.025	0.054391518999668\\
0.024390243902439	0.055252184452053\\
0.0238095238095238	0.0550816117749831\\
0.0232558139534884	0.0537816449558148\\
0.0227272727272727	0.0521142675808788\\
0.0222222222222222	0.0494419039168887\\
0.0217391304347826	0.0461973089722774\\
0.0212765957446809	0.0429161941989022\\
0.0208333333333333	0.0397712596827144\\
0.0204081632653061	0.0365787834338898\\
0.02	0.0332161249931699\\
0.0196078431372549	0.0305443195435045\\
0.0192307692307692	0.0285856870105494\\
0.0188679245283019	0.0262387915020232\\
0.0185185185185185	0.0245619372612227\\
0.0181818181818182	0.0236017634232972\\
0.0178571428571429	0.0222482758547224\\
0.0175438596491228	0.0215635969432242\\
0.0172413793103448	0.0212241232375625\\
0.0169491525423729	0.021204927772505\\
0.0166666666666667	0.0215039335409442\\
0.0163934426229508	0.0212427327208502\\
0.0161290322580645	0.0214152796311637\\
0.0158730158730159	0.0216075667773961\\
0.015625	0.0216671139021472\\
0.0153846153846154	0.0214856015337279\\
0.0151515151515152	0.0211964774075035\\
0.0149253731343284	0.0206939889658466\\
0.0147058823529412	0.0199689510504388\\
0.0144927536231884	0.0191415061254876\\
0.0142857142857143	0.0183182515760114\\
0.0140845070422535	0.0175969564399692\\
0.0138888888888889	0.0162237798361344\\
0.0136986301369863	0.0152767313628006\\
0.0135135135135135	0.0144080708780443\\
0.0133333333333333	0.0137736472957677\\
0.0131578947368421	0.013090336938995\\
0.012987012987013	0.0127818187483797\\
0.0128205128205128	0.0120866779311255\\
0.0126582278481013	0.01170394995638\\
0.0125	0.011375221138751\\
0.0123456790123457	0.0111524465633037\\
0.0121951219512195	0.0112784658743976\\
0.0120481927710843	0.011145045114576\\
0.0119047619047619	0.0110285534613244\\
0.0117647058823529	0.01106724360654\\
0.0116279069767442	0.0111345540637302\\
0.0114942528735632	0.0110434144220335\\
0.0113636363636364	0.0110882395071341\\
0.0112359550561798	0.0110846773664308\\
0.0111111111111111	0.0109456523824916\\
0.010989010989011	0.010696438615823\\
0.0108695652173913	0.0103902584300695\\
0.010752688172043	0.0101165330978157\\
0.0106382978723404	0.00983428592553537\\
0.0105263157894737	0.00934997689235395\\
0.0104166666666667	0.00905050852841427\\
0.0103092783505155	0.0086541363771957\\
0.0102040816326531	0.00824644316872725\\
0.0101010101010101	0.00784299752233864\\
0.01	0.00762709699210262\\
};
\addlegendentry{PCA (IT)}

\addplot [color=color2, mark=*, mark options={solid, color2, scale=0.75}]
  table[row sep=crcr]{%
0.1	0.978156360499254\\
0.0909090909090909	0.643201702579463\\
0.0833333333333333	0.418193661474361\\
0.0769230769230769	0.292553868960584\\
0.0714285714285714	0.229439896260611\\
0.0666666666666667	0.210263994055173\\
0.0625	0.216889199922253\\
0.0588235294117647	0.233501866851585\\
0.0555555555555556	0.245867439956932\\
0.0526315789473684	0.246101242642252\\
0.05	0.238240955249519\\
0.0476190476190476	0.222874667229214\\
0.0454545454545455	0.205137314100043\\
0.0434782608695652	0.184417910678896\\
0.0416666666666667	0.163771156776297\\
0.04	0.143816931613633\\
0.0384615384615385	0.12524829279333\\
0.037037037037037	0.108770842578674\\
0.0357142857142857	0.0941160884392973\\
0.0344827586206897	0.0819676950685668\\
0.0333333333333333	0.0724790084340503\\
0.032258064516129	0.0648076873989005\\
0.03125	0.0596677169497717\\
0.0303030303030303	0.0562945718874861\\
0.0294117647058824	0.0544843394805206\\
0.0285714285714286	0.0535437027355528\\
0.0277777777777778	0.0535179889411452\\
0.027027027027027	0.0539462988930994\\
0.0263157894736842	0.0540531106457696\\
0.0256410256410256	0.0540243471511648\\
0.025	0.0532069717466328\\
0.024390243902439	0.051418602316847\\
0.0238095238095238	0.0491680911751224\\
0.0232558139534884	0.0464243206826778\\
0.0227272727272727	0.043416902392746\\
0.0222222222222222	0.0401908560512059\\
0.0217391304347826	0.0372175023574712\\
0.0212765957446809	0.034252912773229\\
0.0208333333333333	0.0315357837156238\\
0.0204081632653061	0.0291109452711149\\
0.02	0.0270825095624536\\
0.0196078431372549	0.0253566419754172\\
0.0192307692307692	0.024041403383487\\
0.0188679245283019	0.023031325360261\\
0.0185185185185185	0.0223063473815053\\
0.0181818181818182	0.0218650401475449\\
0.0178571428571429	0.0216115456097947\\
0.0175438596491228	0.0216068616918941\\
0.0172413793103448	0.0215735026113046\\
0.0169491525423729	0.0214622062996686\\
0.0166666666666667	0.0214352737134051\\
0.0163934426229508	0.0212201312970821\\
0.0161290322580645	0.0208789156728946\\
0.0158730158730159	0.0203551075963269\\
0.015625	0.0196559919958235\\
0.0153846153846154	0.0188229374151927\\
0.0151515151515152	0.0178964134285353\\
0.0149253731343284	0.0170566614224281\\
0.0147058823529412	0.0161957481223016\\
0.0144927536231884	0.0153375786890422\\
0.0142857142857143	0.0145020287269348\\
0.0140845070422535	0.0138250233806998\\
0.0138888888888889	0.0131606352901272\\
0.0136986301369863	0.012684486764682\\
0.0135135135135135	0.0122817479980499\\
0.0133333333333333	0.0119530728686783\\
0.0131578947368421	0.0116949519383391\\
0.012987012987013	0.0114955135870645\\
0.0128205128205128	0.0114075510239928\\
0.0126582278481013	0.0113532490749257\\
0.0125	0.0113028081000399\\
0.0123456790123457	0.0112551787821351\\
0.0121951219512195	0.0111858503282409\\
0.0120481927710843	0.0110566184093228\\
0.0119047619047619	0.0108965112784416\\
0.0117647058823529	0.0107152005531539\\
0.0116279069767442	0.0104663299397831\\
0.0114942528735632	0.0101347645595284\\
0.0113636363636364	0.00978613995461597\\
0.0112359550561798	0.00943693751391272\\
0.0111111111111111	0.00907181701312165\\
0.010989010989011	0.00875066578278005\\
0.0108695652173913	0.00844153890735466\\
0.010752688172043	0.0081372527890462\\
0.0106382978723404	0.00786626977522253\\
0.0105263157894737	0.00762737028130633\\
0.0104166666666667	0.00744365058744112\\
0.0103092783505155	0.00730250619233797\\
0.0102040816326531	0.00716652062120193\\
0.0101010101010101	0.00708934705418773\\
0.01	0.00700538301269482\\
};
\addlegendentry{Com (IT)}

\addplot [dashed]
  table[row sep=crcr]{%
0.1	1e1\\
0.01	1e-1\\
};

\end{axis}
\end{tikzpicture}%

%% file: BlackScholesNDEuropean-K=40-T=2-vol=0.3-rho=0.8-r=0.05.tikz
%
%
%
\begin{tikzpicture}

\begin{axis}[%
width=0.37\textwidth,
scale only axis,
xmode=log,
xmin=0.01,
xmax=0.1,
xminorticks=true,
xlabel={1/m},
ymode=log,
ymin=1e-3,
ymax=1e2,
yminorticks=true,
ylabel={error},
legend cell align={left},
legend style={nodes={scale=0.75, transform shape}}
]
\addplot [color=color1, mark=square*, mark options={solid, color1, scale=0.75}]
  table[row sep=crcr]{%
0.1	3.72265996061185\\
0.0909090909090909	2.82751058465252\\
0.0833333333333333	2.09049516247921\\
0.0769230769230769	1.49605450953783\\
0.0714285714285714	1.02209282182759\\
0.0666666666666667	0.676892199049098\\
0.0625	0.424078152599591\\
0.0588235294117647	0.244932492950463\\
0.0555555555555556	0.123697977890313\\
0.0526315789473684	0.0472958323538695\\
0.05	0.00568432397823315\\
0.0476190476190476	0.0117952951353608\\
0.0454545454545455	0.0123557430311747\\
0.0434782608695652	0.00358097373338406\\
0.0416666666666667	0.0119753740069157\\
0.04	0.0316266624180388\\
0.0384615384615385	0.0526374391843039\\
0.037037037037037	0.0715106673741737\\
0.0357142857142857	0.0892198240215638\\
0.0344827586206897	0.104457556368195\\
0.0333333333333333	0.116781473896169\\
0.032258064516129	0.124722428396367\\
0.03125	0.12984606097438\\
0.0303030303030303	0.132366850418727\\
0.0294117647058824	0.133026898973743\\
0.0285714285714286	0.131046295171982\\
0.0277777777777778	0.127454048471181\\
0.027027027027027	0.122889452688432\\
0.0263157894736842	0.117854447100315\\
0.0256410256410256	0.11201906659344\\
0.025	0.105784225712138\\
0.024390243902439	0.0998344662037853\\
0.0238095238095238	0.0942617878996987\\
0.0232558139534884	0.0888578714766277\\
0.0227272727272727	0.0836161697251425\\
0.0222222222222222	0.0789107134221565\\
0.0217391304347826	0.0747890884315199\\
0.0212765957446809	0.0709826878558486\\
0.0208333333333333	0.0673277356435955\\
0.0204081632653061	0.0640596497695354\\
0.02	0.0611286050028177\\
0.0196078431372549	0.0584311280382446\\
0.0192307692307692	0.0556863666546734\\
0.0188679245283019	0.0532498917905784\\
0.0185185185185185	0.0511221083830158\\
0.0181818181818182	0.0492090599695496\\
0.0178571428571429	0.04726580439812\\
0.0175438596491228	0.0455540895508326\\
0.0172413793103448	0.0440601419469404\\
0.0169491525423729	0.0427381900432113\\
0.0166666666666667	0.0413843886126815\\
0.0163934426229508	0.0401578506347802\\
0.0161290322580645	0.0390502700634539\\
0.0158730158730159	0.0380710738188457\\
0.015625	0.0370664307926445\\
0.0153846153846154	0.036099853687281\\
0.0151515151515152	0.0352244360745342\\
0.0149253731343284	0.0344119626944455\\
0.0147058823529412	0.033580594850994\\
0.0144927536231884	0.0327323053249931\\
0.0142857142857143	0.031932617211706\\
0.0140845070422535	0.0311778479431757\\
0.0138888888888889	0.0304172458394048\\
0.0136986301369863	0.0296275252053544\\
0.0135135135135135	0.0288800903149191\\
0.0133333333333333	0.0281886456918725\\
0.0131578947368421	0.0275183555005842\\
0.012987012987013	0.0267947106116813\\
0.0128205128205128	0.0261284544020501\\
0.0126582278481013	0.0255012100126075\\
0.0125	0.0249001857234195\\
0.0123456790123457	0.0242529781517753\\
0.0121951219512195	0.0236378223283431\\
0.0120481927710843	0.0230715233042442\\
0.0119047619047619	0.0225359661582747\\
0.0117647058823529	0.0219645493717957\\
0.0116279069767442	0.0214320090560873\\
0.0114942528735632	0.0209340997047418\\
0.0113636363636364	0.0204731084059064\\
0.0112359550561798	0.0199882795919084\\
0.0111111111111111	0.0195187955135143\\
0.010989010989011	0.0190806153183924\\
0.0108695652173913	0.0186713637289797\\
0.010752688172043	0.0182518047931035\\
0.0106382978723404	0.0178411235663667\\
0.0105263157894737	0.0174615073044064\\
0.0104166666666667	0.0171053506888033\\
0.0103092783505155	0.0167524631597304\\
0.0102040816326531	0.0163965472789203\\
0.0101010101010101	0.0160647848509337\\
0.01	0.0157529070406897\\
};
\addlegendentry{PCA}

\addplot [color=color2, mark=*, mark options={solid, color2, scale=0.75}]
  table[row sep=crcr]{%
0.1	3.21777530043942\\
0.0909090909090909	2.41618022913473\\
0.0833333333333333	1.7716918076861\\
0.0769230769230769	1.26622135253974\\
0.0714285714285714	0.880995191394244\\
0.0666666666666667	0.595398403773436\\
0.0625	0.389823610143312\\
0.0588235294117647	0.247305038258721\\
0.0555555555555556	0.154035081764365\\
0.0526315789473684	0.0972426946167175\\
0.05	0.0667478498191647\\
0.0476190476190476	0.0545580784944502\\
0.0454545454545455	0.054924360881917\\
0.0434782608695652	0.062703536247608\\
0.0416666666666667	0.0740934503494834\\
0.04	0.0864405039734377\\
0.0384615384615385	0.0984233316869485\\
0.037037037037037	0.108867708358178\\
0.0357142857142857	0.117212630716206\\
0.0344827586206897	0.123286877199297\\
0.0333333333333333	0.127273301587863\\
0.032258064516129	0.12918256353418\\
0.03125	0.12927973190628\\
0.0303030303030303	0.127911495985532\\
0.0294117647058824	0.125494888700402\\
0.0285714285714286	0.122172632448373\\
0.0277777777777778	0.11818629149631\\
0.027027027027027	0.113761850397617\\
0.0263157894736842	0.109148830744131\\
0.0256410256410256	0.10438158580308\\
0.025	0.0995472581763606\\
0.024390243902439	0.0947069259416562\\
0.0238095238095238	0.0899699789433095\\
0.0232558139534884	0.0853489105493717\\
0.0227272727272727	0.0809078545024482\\
0.0222222222222222	0.0767022975142373\\
0.0217391304347826	0.0727977716561536\\
0.0212765957446809	0.0691563820333467\\
0.0208333333333333	0.0657841932823962\\
0.0204081632653061	0.0626849640984064\\
0.02	0.0598704714411546\\
0.0196078431372549	0.057287776968808\\
0.0192307692307692	0.0549226957408466\\
0.0188679245283019	0.0527662431067837\\
0.0185185185185185	0.0508143770956435\\
0.0181818181818182	0.0490158376245198\\
0.0178571428571429	0.0473500276591325\\
0.0175438596491228	0.0458040649836855\\
0.0172413793103448	0.0443679090046318\\
0.0169491525423729	0.0430075373887098\\
0.0166666666666667	0.0417167392641451\\
0.0163934426229508	0.0404979385978654\\
0.0161290322580645	0.039351717414454\\
0.0158730158730159	0.0382542211016563\\
0.015625	0.0372000522755158\\
0.0153846153846154	0.0361915065586729\\
0.0151515151515152	0.0352263074839998\\
0.0149253731343284	0.0342869255956044\\
0.0147058823529412	0.0333728668202999\\
0.0144927536231884	0.0324913028673599\\
0.0142857142857143	0.0316435866840425\\
0.0140845070422535	0.0308179389505456\\
0.0138888888888889	0.0300136760064582\\
0.0136986301369863	0.0292369029348922\\
0.0135135135135135	0.0284865648442736\\
0.0133333333333333	0.0277527727230718\\
0.0131578947368421	0.0270367913095901\\
0.012987012987013	0.0263464379255574\\
0.0128205128205128	0.0256820975217851\\
0.0126582278481013	0.0250366339745796\\
0.0125	0.024409908661462\\
0.0123456790123457	0.0238074458449908\\
0.0121951219512195	0.0232275027809088\\
0.0120481927710843	0.0226629660190119\\
0.0119047619047619	0.0221142120813829\\
0.0117647058823529	0.0215872203851077\\
0.0116279069767442	0.021081119279029\\
0.0114942528735632	0.0205908080994925\\
0.0113636363636364	0.0201159194970106\\
0.0112359550561798	0.0196606631904235\\
0.0111111111111111	0.0192227548405448\\
0.010989010989011	0.0187969343298242\\
0.0108695652173913	0.0183829289613504\\
0.010752688172043	0.0179854047131789\\
0.0106382978723404	0.0176028809387816\\
0.0105263157894737	0.0172316867964106\\
0.0104166666666667	0.0168714833152714\\
0.0103092783505155	0.0165257987427454\\
0.0102040816326531	0.0161922901047955\\
0.0101010101010101	0.0158671007988742\\
0.01	0.0155497051658617\\
};
\addlegendentry{Com}

\addplot [dashed]
  table[row sep=crcr]{%
0.1	1e1\\
0.01	1e-1\\
};

\end{axis}
\end{tikzpicture}%

%% file: BlackScholesNDAmerican-K=40-T=2-vol=0.3-rho=0.8-r=0.05.tikz
%
%
%
\begin{tikzpicture}

\begin{axis}[%
width=0.37\textwidth,
scale only axis,
xmode=log,
xmin=0.01,
xmax=0.1,
xminorticks=true,
xlabel={1/m},
ymode=log,
ymin=1e-3,
ymax=1e2,
yminorticks=true,
legend cell align={left},
legend style={nodes={scale=0.75, transform shape}}
]
\addplot [color=color1, mark=square*, mark options={solid, color1, scale=0.75}]
  table[row sep=crcr]{%
0.1	3.74865854016886\\
0.0909090909090909	2.75561149086716\\
0.0833333333333333	1.93517754201406\\
0.0769230769230769	1.27317519671952\\
0.0714285714285714	0.760893524540617\\
0.0666666666666667	0.405799369607267\\
0.0625	0.166233933622895\\
0.0588235294117647	0.00944087123641779\\
0.0555555555555556	0.076763750691665\\
0.0526315789473684	0.0940746413144939\\
0.05	0.0721131572959379\\
0.0476190476190476	0.0166931989409882\\
0.0454545454545455	0.0562576725285098\\
0.0434782608695652	0.136994503199672\\
0.0416666666666667	0.214006231167176\\
0.04	0.282000053091114\\
0.0384615384615385	0.337026361864659\\
0.037037037037037	0.364149654958194\\
0.0357142857142857	0.377381505413139\\
0.0344827586206897	0.381021817706362\\
0.0333333333333333	0.372029071826242\\
0.032258064516129	0.358318396701225\\
0.03125	0.335443368394166\\
0.0303030303030303	0.310973783093642\\
0.0294117647058824	0.283330371025935\\
0.0285714285714286	0.253633401173868\\
0.0277777777777778	0.223466554125357\\
0.027027027027027	0.193908382609234\\
0.0263157894736842	0.166355351300915\\
0.0256410256410256	0.141346616421655\\
0.025	0.116519799263971\\
0.024390243902439	0.0953702991345171\\
0.0238095238095238	0.0784276177746381\\
0.0232558139534884	0.0642961620941822\\
0.0227272727272727	0.0524702101480674\\
0.0222222222222222	0.0445844997247509\\
0.0217391304347826	0.037533104248741\\
0.0212765957446809	0.0355288826115601\\
0.0208333333333333	0.0344461856074414\\
0.0204081632653061	0.0361156252052233\\
0.02	0.0377735578658616\\
0.0196078431372549	0.0421027980994122\\
0.0192307692307692	0.0467999509339592\\
0.0188679245283019	0.0520577451098845\\
0.0185185185185185	0.0575943249226007\\
0.0181818181818182	0.0627733108906359\\
0.0178571428571429	0.0682585714219002\\
0.0175438596491228	0.0720565694265751\\
0.0172413793103448	0.0754244014925032\\
0.0169491525423729	0.0773186838578983\\
0.0166666666666667	0.0769329352600092\\
0.0163934426229508	0.0758058203520928\\
0.0161290322580645	0.0730949554539304\\
0.0158730158730159	0.069549197861333\\
0.015625	0.0653720016504309\\
0.0153846153846154	0.0607532148201493\\
0.0151515151515152	0.0567671989218397\\
0.0149253731343284	0.0520119659865923\\
0.0147058823529412	0.0470068495530853\\
0.0144927536231884	0.042748959322509\\
0.0142857142857143	0.0381962266454572\\
0.0140845070422535	0.0345745556371639\\
0.0138888888888889	0.0305301210132551\\
0.0136986301369863	0.0273707129413978\\
0.0135135135135135	0.0247319670318014\\
0.0133333333333333	0.0223209258999155\\
0.0131578947368421	0.0207040926439284\\
0.012987012987013	0.019250820803002\\
0.0128205128205128	0.0182624284048956\\
0.0126582278481013	0.0177990743764589\\
0.0125	0.0179678948826112\\
0.0123456790123457	0.0185321504776166\\
0.0121951219512195	0.0190995520143868\\
0.0120481927710843	0.0200844374810139\\
0.0119047619047619	0.0213508533307323\\
0.0117647058823529	0.0227039400989728\\
0.0116279069767442	0.0237977490621137\\
0.0114942528735632	0.0251598274019642\\
0.0113636363636364	0.0262901872685131\\
0.0112359550561798	0.0272601921533226\\
0.0111111111111111	0.0284794439737777\\
0.010989010989011	0.0290641751754617\\
0.0108695652173913	0.0290024086310083\\
0.010752688172043	0.028825878052734\\
0.0106382978723404	0.0281002682637972\\
0.0105263157894737	0.0272691958638536\\
0.0104166666666667	0.0260818955187787\\
0.0103092783505155	0.024745309900803\\
0.0102040816326531	0.0234350237081546\\
0.0101010101010101	0.0221464909311404\\
0.01	0.0205263635165114\\
};
\addlegendentry{PCA (IT)}

\addplot [color=color2, mark=*, mark options={solid, color2, scale=0.75}]
  table[row sep=crcr]{%
0.1	3.16513705595493\\
0.0909090909090909	2.27499106706284\\
0.0833333333333333	1.55574492159627\\
0.0769230769230769	0.997399350802563\\
0.0714285714285714	0.593244326865659\\
0.0666666666666667	0.314969928142516\\
0.0625	0.138035510837697\\
0.0588235294117647	0.0412759097748472\\
0.0555555555555556	0.0105031751908067\\
0.0526315789473684	0.0280788677151005\\
0.05	0.0772213117507663\\
0.0476190476190476	0.144479339077133\\
0.0454545454545455	0.216280880477199\\
0.0434782608695652	0.279346754485762\\
0.0416666666666667	0.324198056117623\\
0.04	0.35219962030834\\
0.0384615384615385	0.364718381941265\\
0.037037037037037	0.366202949752645\\
0.0357142857142857	0.357074221836785\\
0.0344827586206897	0.340833656509484\\
0.0333333333333333	0.318468685251539\\
0.032258064516129	0.293404152097565\\
0.03125	0.265696293106969\\
0.0303030303030303	0.237498519120229\\
0.0294117647058824	0.209794806048031\\
0.0285714285714286	0.183040564419357\\
0.0277777777777778	0.15796216976339\\
0.027027027027027	0.135428882840995\\
0.0263157894736842	0.115077985664201\\
0.0256410256410256	0.0976394756561527\\
0.025	0.0829647574726984\\
0.024390243902439	0.0711557726157919\\
0.0238095238095238	0.062139282171505\\
0.0232558139534884	0.0554891340244499\\
0.0227272727272727	0.051231984942758\\
0.0222222222222222	0.0491514191082958\\
0.0217391304347826	0.04935793258272\\
0.0212765957446809	0.051100781324692\\
0.0208333333333333	0.0541054964073178\\
0.0204081632653061	0.0581585539554945\\
0.02	0.0627352550576648\\
0.0196078431372549	0.0670858263132086\\
0.0192307692307692	0.0708518373891751\\
0.0188679245283019	0.0738394027061675\\
0.0185185185185185	0.0755330963327845\\
0.0181818181818182	0.0758672128962088\\
0.0178571428571429	0.0749458031414685\\
0.0175438596491228	0.0729743200730653\\
0.0172413793103448	0.0702960427068646\\
0.0169491525423729	0.0669213036741789\\
0.0166666666666667	0.0631492039814496\\
0.0163934426229508	0.0589899020321143\\
0.0161290322580645	0.0546689425076456\\
0.0158730158730159	0.0505014230274723\\
0.015625	0.0462574713796879\\
0.0153846153846154	0.0422200084357547\\
0.0151515151515152	0.0385007014043213\\
0.0149253731343284	0.0350108440495571\\
0.0147058823529412	0.0319801902398416\\
0.0144927536231884	0.029399184870404\\
0.0142857142857143	0.0271759623820822\\
0.0140845070422535	0.0254559150213955\\
0.0138888888888889	0.0242546065403806\\
0.0136986301369863	0.0235291792259522\\
0.0135135135135135	0.0231878845276903\\
0.0133333333333333	0.0231731466335168\\
0.0131578947368421	0.023568934379373\\
0.012987012987013	0.0242266864007394\\
0.0128205128205128	0.0250453519201921\\
0.0126582278481013	0.0260185155264026\\
0.0125	0.0269738870840923\\
0.0123456790123457	0.0278882436007155\\
0.0121951219512195	0.0286248470352515\\
0.0120481927710843	0.0290708477363832\\
0.0119047619047619	0.0292153986812451\\
0.0117647058823529	0.0290949662234414\\
0.0116279069767442	0.0287044622644883\\
0.0114942528735632	0.0280440289392052\\
0.0113636363636364	0.0271404121697021\\
0.0112359550561798	0.0261009902602023\\
0.0111111111111111	0.0249579972358456\\
0.010989010989011	0.023708667478731\\
0.0108695652173913	0.0224122181963855\\
0.010752688172043	0.021178877434032\\
0.0106382978723404	0.0199567937920948\\
0.0105263157894737	0.0188139117372206\\
0.0104166666666667	0.0177214224464377\\
0.0103092783505155	0.0167176901121877\\
0.0102040816326531	0.0157914772396266\\
0.0101010101010101	0.0150160954181509\\
0.01	0.0143829153803763\\
};
\addlegendentry{Com (IT)}

\addplot [dashed]
  table[row sep=crcr]{%
0.1	1e1\\
0.01	1e-1\\
};

\end{axis}
\end{tikzpicture}%